\documentclass[a4paper,oneside]{amsbook}
\usepackage[14pt]{extsizes}
\usepackage[cp1251]{inputenc}
\usepackage[russian]{babel}
\usepackage[matrix,arrow,tips]{xy}
\usepackage{amssymb,amsmath,amsthm,amscd}
\usepackage{type1cm}
\usepackage{enumerate,comment}
\usepackage{calrsfs}
\usepackage[normalem]{ulem}
\usepackage{pb-diagram}
\usepackage[pdftex,unicode]{hyperref}
\hypersetup{colorlinks=true, linkcolor=blue, filecolor=blue,
urlcolor=blue} 
\theoremstyle{plain}
 {
  \swapnumbers

\newtheorem{thm}{\sc\bfseries Tеорема}[section]
\newtheorem{prop}[thm]{\sc Предложение}
\newtheorem{lemma}[thm]{\sc Лемма}
\newtheorem{cor}[thm]{\sc Следствие}

 }
\theoremstyle{definition}
 {
  \swapnumbers
\newtheorem{dfn}[thm]{\sc Определение}

 }
\theoremstyle{remark}
 {
  \swapnumbers

\newtheorem*{remark*}{\sc Замечание}
\newtheorem{predo}[thm]{\sc Предостережение}
 }
\newcommand*{\abs}[1]{\left\vert#1\right\vert}

\newcommand*{\ehatotimes}{\mathbin{\hat{\underset{\eps}{\otimes}}}}

\newcommand*{\set}[1]{\left\{#1\right\}}
\newcommand*{\Set}[2]{\left\{#1\colon #2\right\}}
\newcommand*{\map}[3]{#1\colon #2\to #3}
\newcommand*{\cl}{\mathop{\rm cl}\nolimits}
\newcommand*{\ob}{\mathop{\rm ob}\nolimits}

\makeatletter
\renewcommand{\section}{\@startsection{section}{1}%
  \z@{.7\linespacing\@plus\linespacing}{.5\linespacing}%
  {\normalfont\bfseries\centering\S}}
\makeatother

\DeclareMathOperator{\Ker}{Ker} \DeclareMathOperator{\St}{St}

\renewcommand{\baselinestretch}{1.65}

\def\лл{\glqq}
\def\пл{\grqq}

\def\dl{\delta}

\def\lm{\lambda}

\def\vp{\varphi}
\def\vphi{\varphi}

\def\eps{\varepsilon}
\def\Cnum{\mathbb C}
\def\Cbb{\mathbb C}

\def\R{\mathbb R}

\def\Chi{\mathrm X}
\def\gel{\mathcal G}
\def\C{\mathcal C}
\def\Ccal{\mathcal C}
\def\Lcal{\mathcal L}
\def\Pcal{\mathcal P}
\def\les{\leqslant}
\def\ges{\geqslant}
\def\sbs{\subset}

\def\ul{\underline}

\def\ul{\underline}
\def\ol{\overline}

\def\eqq{\longleftrightarrow}
\def\Eq{\Leftrightarrow}

\def\опр{\stackrel{\text{опр}}{\eqq}}
\def\проп{$\langle\dots\rangle$}

\def\нач{$\boxed{\bot}$}
\def\шаг{$\boxed{\Rsh}$}

\begin{document}

\raggedbottom

{%
\bf
\renewcommand{\baselinestretch}{1}
\selectfont
\thispagestyle{empty}

\begin{center}
{\small МОСКОВСКИЙ ОРДЕНА ЛЕНИНА \\И ОРДЕНА ТРУДОВОГО КРАСНОГО ЗНАМЕНИ\\
ГОСУДАРСТВЕННЫЙ ПЕДАГОГИЧЕСКИЙ ИНСТИТУТ\\ имени В.~И.~ЛЕНИНА}

\medskip\noindent
{\small РЯЗАНСКИЙ ОРДЕНА "<ЗНАК ПОЧЁТА">\\
ГОСУДАРСТВЕННЫЙ ПЕДАГОГИЧЕСКИЙ ИНСТИТУТ\\ имени С.~А.~ЕСЕНИНА}
\end{center}
\vspace{-24pt}
\noindent\hrulefill

\begin{flushright}
На правах рукописи\\
УДК 519.4
\end{flushright}

\vspace{1cm}

\begin{center}
НАЗИЕВ Асланбек Хамидович

\bigskip
$K^*$-АЛГЕБРЫ И ДВОЙСТВЕННОСТЬ

\bigskip
01.01.06 --- математическая логика,
алгебра и теория чисел

\vspace{1cm}

ДИССЕРТАЦИЯ\\
на соискание учёной степени\\
кандидата физико-математических наук
\end{center}

\vspace{1cm}

\mbox{}\hfill\parbox{.52\textwidth}{%
Научный руководитель:\\
доктор физико-математических\\
наук, профессор\, {\sc\bfseries Д.~А.~Райков}}

\vspace{4cm}

\begin{center}
МОСКВА -- 1987
\end{center}

}

\eject

\tableofcontents

\chapter*{Введение}

Один из основных результатов диссертации --- теорема двойственности для категории \uwave{всех} --- как коммутативных, так и некоммутативных~--- отделимых локально компактных групп. Теорема двойственности для \uwave{коммутативных} отделимых локально компактных групп была получена Л.~С.~Понтрягиным в 1934 году~\cite{b35}, см. также~\cite{b36}. С тех пор неоднократно предпринимались попытки обобщить эту теорему на различные категории \uwave{некоммутативных} локально компактных групп. Первое обобщение такого рода было получено Т.~Таннакой~\cite{b47} и М.~Г.~Крейном~\cite{b18}. Т.~Таннака установил, что любая отделимая \uwave{компактная} группа с точностью до изоморфизма определяется своей представляющей алгеброй, а М.~Г.~Крейн дал для представляющих алгебр абстрактную характеризацию. У.~Стайнспринг~\cite{b43} обобщил теорему Таннаки-Крейна на унимодулярные группы. Дальнейшие обобщения были получены Г.~И.~Кацем, М.~Такесаки, Н.~Татсуумой, X.~Чу --- см. \cite{b13},  \cite{b44}, \cite{b48}, \cite{b45}, \cite{b46}, \cite{b58}, а также~\cite{b54}.

Во всех перечисленных работах существенным образом использовалась теория унитарных представлений групп. Попытки получения теоремы двойственности для некоммутативных локально компактных групп на этой основе предпринимались и другими авторами --- см., например, \cite{b39}, \cite{b50}, \cite{b53}, \cite{b59}. Однако теорема двойственности для категории всех отделимых локально компактных групп на этом пути до сих пор не получена. В 1981 году К.~Маккеннон даже высказал сомнение в том, что её вообще можно получить~--- см.~\cite{b22}.

Между тем в 1963 году оформилось другое направление исследований по теории двойственности для некоммутативных топологических групп. Начало ему было положено статьёй Дж.~Келли~\cite{b15}. В ней описывались новые объекты, названные автором когруппами, и давался набросок теории двойственности между отделимыми \uwave{компактными} группами и некоторыми когруппами.

В 1965 году Г.~Хохшильд во всех деталях реализовал идеи Дж.~Келли для \uwave{компактных} групп \cite{b55}. В 1970 году К.~X.~Хофман предложил в монографии \cite{b56} другую реализацию этих идей для \uwave{компактных} групп. Категорию, дуальную по К.~Х.~Хофману к категории отделимых компактных групп, образуют $C^*$-алгебры Хопфа, т.~е. $C^*$-алгебры, наделённые ещё структурой, определяемой так же, как и структура алгебры Хопфа, за исключением того, что вместо обычного тензорного произведения берётся найденное в \cite{b56} $C^*$-алгебраическое тензорное произведение. При этом К.~X.~Хофман доказал аналогичную теорему и для \uwave{компактных} полугрупп. В 1972 году С.~Санкаран и С.~Селесник опубликовали статью \cite{b40} с другим вариантом построений К.~X.~Хофмана для \uwave{компактных} полугрупп и групп.

Первые теоремы двойственности для всех отделимых локально компактных полугрупп и групп были опубликованы автором диссертации в 1973 году~--- см. \cite{b26}. Эти теоремы являются обобщениями теорем К.~X.~Хофмана. Однако доказательство их не свелось к простому перенесению рассуждений К.~X.~Хофмана в более общую ситуацию. Во-первых, К.~Х.~Хофман (также, как С.~Санкаран и С.~Селесник) в своих построениях использовал уже готовую двойственность для отделимых \uwave{компактных} пространств --- см., например, \cite{b5}, \cite{b10}, \cite{b33}. Поскольку для некомпактных отделимых \uwave{локально компактных} пространств такой двойственности ещё не было, автору пришлось строить её самостоятельно. Это, в свою очередь, потребовало дополнительных исследований по топологическим алгебрам. А, во-вторых, для определения тензорного произведения  $C^*$-алгебр К.~Х.~Хофман использовал теорию унитарных представлений этих алгебр, а С.~Санкаран и С.~Селесник --- тензорные произведения банаховых пространств и некоторые дополнительные ухищрения. Ни то, ни другое не представлялось возможным использовать в рассматриваемой нами более общей ситуации. Понадобилось отдельное изучение тензорных произведений топологических инволютивных алгебр.

В силу указанных причин построение теории двойственности для категорий всех отделимых локально компактных полугрупп и групп распалось на несколько этапов:

1) изучение топологических инволютивных алгебр;

2) построение двойственностей для пространств, более общих, чем компактные;

3) изучение тензорных произведений топологических инволютивных алгебр;

4) собственно получение теорем двойственности для всех отделимых локально компактных полугрупп и групп.

При этом нам удалось получить теоремы двойственности для категорий пространств и полугрупп даже более широких, чем локально компактные, а теоремы двойственности для локально компактного случая вывести в качестве простых следствий из указанных более общих теорем. Изложению результатов, полученных в процессе реализации намеченной выше программы, и посвящена диссертация.

Диссертация состоит из настоящего введения и пяти глав. В первой главе вводятся некоторые обозначения, сообщаются предварительные сведения о топологических инволютивных алгебрах и отмечаются некоторые особенности принятого в диссертации подхода к изучению топологических инволютивных алгебр.

Собственные результаты автора представлены в главах 2--5. Опишем содержание этих глав более подробно.

\medskip
\ul{Глава 2: $LC^*$-топологии на инволютивных алгебрах}

Алгеброй всюду в диссертации называется коммутативная, ассоциативная комплексная алгебра с единицей, морфизмом алгебр --- сохраняющий единицу гомоморфизм, характером алгебры --- морфизм из неё в поле $\Cbb$ комплексных чисел. Морфизмы инволютивных алгебр, сохраняющие инволюцию, называются эрмитовыми. Множество всех характеров алгебры $A$ обозначается через $\Chi(A)$, множество всех эрмитовых характеров инволютивной алгебры $A$ обозначается через $\Chi^*(A)$. Эти множества наделяются топологиями поточечной сходимости на $A$ и называются (будучи наделены ими) пространствами (всех, всех эрмитовых) характеров. Каждый элемент $a\in A$  порождает функцию $\hat{a}\colon \chi\mapsto \chi(a)$ из $\Chi(A)$ в $\Cbb$; сужение этой функции на множество $\Chi^*(A)$ обозначается тем же знаком. Все эти функции непрерывны и, более того, слабая топология, определяемая этими функциями на $\Chi(A)$ и $\Chi^*(A)$, совпадает с исходной. Отображение $a\mapsto\hat{a}$ из $A$ в алгебру $C(\Chi(A))$ всех непрерывных комплекснозначных функций на пространстве $\Chi(A)$ называется преобразованием Гельфанда алгебры $A$ и обозначается через $\gel_A$ или, если ясно о какой $A$ идёт речь, через $\gel$. То же~--- для инволютивных алгебр и отображений в~$C(\Chi^*(A))$.

Преднорма $p$ на инволютивной алгебре $A$ называется $C^*$-преднормой, если для всех $a\in A$ $p(aa^*) = p(a)^2$. Топологию на инволютивной алгебре, определяемую системой $C^*$-преднорм, мы называем $LC^*$-топологией, инволютивную алгебру, наделённую $LC^*$-топологией, --- $LC^*$-алгеброй, полную отделимую $LC^*$-алгебру --- $K^*$-алгеброй. Множество всех эрмитовых характеров алгебры $A$, непрерывных относительно $C^*$-преднормы~$p$ (соотв., $LC^*$-топологии $\tau$) на $A$, обозначается через $\Chi_p^*(A)$ (соотв., $\Chi_c^*(A,\ \tau))$. В случаях, когда ясно, о какой $A$ или $\tau$ идёт речь, пишем просто~$\Chi_p^*$, $\Chi_c^*(A)$.

Параграф 1 посвящен описанию всех $LC^*$-топологий на данной инволютивной алгебре. Здесь доказывается, что каждая $LC^*$-топология на инволютивной алгебре $A$ является топологией равномерной сходимости на некоторой совокупности компактных подмножеств пространства $X^*(A)$ и выясняется вопрос о существовании слабейшей и сильнейшей среди всех таких топологий. Отдельно рассмотрен случай, когда $A$ есть алгебра $C(T)$ всех непрерывных комплекснозначных функций на тихоновском (т.~е. отделимом вполне регулярном) пространстве $T$.

В §2 изучаются соотношения между $LC^*$-топологиями на инволютивных алгебрах и соответствующими им пространствами непрерывных эрмитовых характеров. Мы говорим, что топология $\tau$ на инволютивной алгебре $A$ согласуется с двойственностью между $A$ и непустым подмножеством $\Gamma$ множества $\Chi(A)$, если $\Chi_c(A,\ \tau)= \Gamma$, и доказываем, что каждая $LC^*$-топология на инволютивной алгебре $A$, согласующаяся с двойственностью между $A$ и $\Gamma$, есть топология равномерной сходимости на некоторой совокупности компактных подмножеств пространства $\Chi^*(A)$, содержащихся в $\Gamma$ и покрывающих $\Gamma$. Описываются слабейшая и сильнейшая среди всех таких топологий. Мы называем, далее, $C^*$-преднорму~$p$ на топологической инволютивной алгебре спектральной, если каждый  $p$-непрерывный эрмитов характер на~$A$ непрерывен относительно топологии алгебры~$A$. Мы доказываем, что преобразование Гельфанда $\mathcal{G}_A$ является гомеоморфизмом алгебры $A$ в алгебру $C_{co}(\Chi^*(A))$ (где через \glqq со\grqq\ обозначена компактно-открытая топология) тогда и только тогда, когда на $A$ каждая спектральная $C^*$-преднорма непрерывна и алгебра~$A$ отделима. Здесь же доказывается, что категории $LC^*$-алгебр, отделимых $LC^*$-алгебр и $K^*$-алгебр являются рефлексивными подкатегориями категории топологических инволютивных алгебр.

Прежде чем продолжить описание результатов главы 1, напомним, что изучая топологические инволютивные алгебры, мы преследуем вполне определённую цель: получение теорем двойственности для некоторых категорий топологических инволютивных алгебр и топологических пространств. Топологические пространства в доказательствах этих теорем будут возникать как пространства непрерывных эрмитовых характеров рассматриваемых алгебр. Поэтому нас интересуют в первую очередь такие результаты о топологических инволютивных алгебрах, в которых описываются соотношения между свойствами самих алгебр и пространств характеров на них. Для получения этих результатов мы применяем подход, который несколько отличается от общепринятого. Отличие это состоит в том, что мы не ограничиваемся, как это делают обычно, изучением одного только пространства непрерывных эрмитовых характеров, а рассматриваем соотношения между этим пространством и пространством всех эрмитовых характеров. Результаты следующих двух параграфов подтверждают плодотворность такого подхода.

В §3 изучаются полунепрерывные снизу $C^*$-преднормы на $LC^*$-ал\-геб\-рах (выбор именно этого свойства $C^*$-преднорм проясняется результатами следующего параграфа). Основной результат этого параграфа, теорема (3.1), гласит: $C^*$-преднорма $p$ на $LC^*$-алгебре $A$ полунепрерывна снизу тогда и только тогда, когда $\Chi_p^*(A)$ является замыканием в $\Chi^*(A)$ пересечения $\Chi_p^*(A)\cap\Chi_c^*(A)$. Отсюда, в частности, следует, что $C^*$-преднорма на $LC^*$-алгебре полунепрерывна снизу тогда и только тогда, когда она является верхней огибающей некоторого семейства минимальных непрерывных $C^*$-преднорм. Таким образом, после указания системы $C^*$-преднорм, задающей $LC^*$-топологию на инволютивной алгебре $A$, выделение непрерывных и полунепрерывных снизу $C^*$-преднорм на алгебре $A$ может быть осуществлено исключительно в терминах порядка на множестве $\mathcal{P}(A)$ всех $C^*$-преднорм на алгебре $A$.

В §4 изучаются $C^*$-бочечные алгебры. Так мы называем $LC^*$-алгебры, на которых все полунепрерывные снизу  $C^*$-преднормы непрерывны. Основным результатом этого параграфа является теорема (4.2): $LC^*$-алгебра $A$ $C^*$-бочечна тогда и только тогда, когда: а) пересечение множества $\Chi_c^*(A)$ с каждым компактом из $\Chi^*(A)$ замкнуто в $\Chi^*(A)$ и б) топология алгебры $A$ совпадает с топологией равномерной сходимости на всех компактах из пространства $\Chi_c^*(A)$. Эта теорема является значительным обобщением и даже усилением хорошо известной и важной теоремы Нахбина и Сироты (см. \cite{b34}, \cite{b42}, а также \cite{b1}) о бочечных пространствах непрерывных функций. Из неё следует, в частности, что категория $C^*$-бочечных $LC^*$-алгебр является корефлексивной подкатегорией в категории всех $LC^*$-алгебр и что для всех $LC^*$-топологий на инволютивной алгебре $A$, дающих одно и то же пространство непрерывных эрмитовых характеров, ассоциированные $C^*$-бочечные топологии на $A$ совпадают. Кроме того, она является основой для получения нескольких теорем двойственности.

\eject
\ul{Глава 3: Теоремы двойственности для компактологических}\par
\ul{и некоторых классов топологических пространств}

Напомним, что компактологическим пространством называется \cite{b6} пара $(S,\ \varkappa)$, образованная множеством $S$ и его покрытием $\varkappa$, удовлетворяющим условиям: 1) $\varkappa$ фильтруется вправо по отношению $\sbs$; 2) каждый $K\in\varkappa$ наделён компактной топологией; 3) если $K$, $L\in\varkappa$ и $K\sbs L$, то $K$ --- подпространство в $L$; 4) если $H \sbs K\in \varkappa$, то замыкание множества $H$ в $K$ с топологией, индуцированной из $K$, принадлежит~$\varkappa$. Морфизмом из компактологического пространства $(S,\ \varkappa)$ в компактологическое пространство $(S_1,\ \varkappa_1)$ называется всякое отображение $u\colon S\to S_1$, удовлетворяющее условию: для любого $K\in\varkappa$ существует $K_1\in\varkappa_1$ такой, что $u[K] \sbs K_1$ и сужение
$u$ на $K$ непрерывно.

С каждой топологией на множестве $S$ естественно связана компактология на том же множестве, состоящая из всех компактов относительно данной топологии. Такие компактологии и компактологические пространства с такими компактологиями называются топологически порождёнными. Компактологическое пространство называется регулярным, если морфизмы из него в $\Cbb$ различают точки. Все рассматриваемые в дальнейшем компактологические пространства предполагаются регулярными.

А.~Бухвальтер \cite[теорема 1.2.2]{b6} доказал, что компактологическое пространство регулярно тогда и только тогда, когда на нём существует вполне регулярная топология, индуцирующая на каждом члене компактологии его исходную топологию. В книге Дж.~Купера \cite{b19} утверждается больше (предложение А.2.2): каждое регулярное компактологическое пространство топологически порождено. В §1 приводится пример, показывающий, что это более сильное утверждение неверно, и доказывается, что компактологическое пространство топологически порождено тогда и только тогда, когда его компактология устойчива относительно образования компактных индуктивных пределов в категории отделимых топологических пространств.

В §2 даётся новое доказательство теоремы А.~Бухвальтера \cite{b6} о том, что категория, дуальная к категории регулярных компактологических пространств, эквивалентна категории $K^*$-алгебр, и выводится ряд (отсутствующих у А.~Бухвальтера) следствий относительно полноты и пополнений отделимых $LC^*$-алгебр.

В §3 доказываются теоремы двойственности для некоторых категорий  $K^*$-алгебр и топологических пространств (вообще говоря, более общих, чем локально компактные). Прежде чем привести здесь основные результаты, напомним некоторые определения. Тихоновское пространство $T$ называется: $k_{\Cbb}$-пространством, если каждая функция \mbox{$T\to\Cbb$}, непрерывная на всех компактах из $T$, непрерывна на всём $T$; $\mu$-про\-с\-тран\-ст\-вом, если каждое подмножество в $T$, на котором ограниченны все функции из $C(T)$, компактно; $Q$-пространством, если каждый эрмитов характер на алгебре $C_{co}(T)$ непрерывен. Говорят (см., например, \cite{b56}), что категории $\mathfrak{A}$ и $\mathfrak{B}$ образуют дуальную пару, если каждая из них эквивалентна дуальной к другой. Из обнаруженных в §3 дуальных пар выделим следующие:

1)	категория $k_{\Cbb}$-пространств и категория $K^*$-алгебр, на которых все спектральные $C^*$-преднормы непрерывны;

2)	категория $\mu k_\Cbb$-пространств и категория $C^*$-бочечных $K^*$-алгебр;

3) категория $Qk_\Cbb$-пространств и категория борнологических (см., например, \cite{b37}) $K^*$-алгебр.

Указываются некоторые следствия из полученных теорем двойственности, на которых мы здесь не будем останавливаться.

Наконец, в §4 доказываются теоремы двойственности для локально компактных пространств. $K^*$-алгебру $A$ мы называем $lc^*$-алгеброй (не $LC^*$-алгеброй!), если для любого непрерывного эрмитова характера $\chi$ на $A$ найдутся $a\in A$ и непрерывная $C^*$-преднорма $p$ на $A$ такие, что $\chi(a) = 1$ и $a\cdot N(p) = 0$. Основными результатами §4 (теоремы (4.2)---(4.5)) выделяются следующие дуальные пары категорий:

1) категория отделимых локально компактных пространств и категория $lc^*$-алгебр;

2)	категория локально компактных  $\mu$-пространств и категория  $C^*$-бочечных $lc^*$-алгебр;

3)	категория локально компактных  $Q$-пространств и категория борнологических $lc^*$-алгебр;

4)	категория  $\sigma$-компактных локально компактных пространств и категория метризуемых $lc^*$-алгебр.

\eject
\ul{Глава 4: Тензорные произведения $LC^*$-алгебр}

Важным этапом в построении двойственности для топологических полугрупп по методу Келли-Хофмана является определение тензорного произведения на соответствующих категориях $LC^*$-алгебр. Если $A_1$, $A_2$ --- $LC^*$-алгебры, то тензорное произведение $A_1\otimes A_2$ векторных пространств $A_1$, $A_2$ естественным образом превращается в инволютивную алгебру, умножение и инволюция в которой определяются условиями
\[
    (x_1 \otimes x_2)^* = x_1^*\otimes x_2^*
\]
для элементарных тензоров $a_1\otimes a_2$ и канонически продолжаются на остальные элементы из $A_1\otimes A_2$. Трудность состоит в выборе подходящей топологии. Теория локально выпуклых топологических векторных пространств предлагает целый спектр так называемых допустимых локально выпуклых топологий на тензорном произведении локально выпуклых пространств~--- см. \cite{b37}, \cite[Т. 2]{b16},~--- однако в случае $LC^*$-алгебр нам нужна не просто локально выпуклая, a $LC^*$-топология.

Для $C^*$-алгебр К.~Хофман в~\cite{b56} использовал найденное им $C^*$-алгебраическое тензорное произведение $\bar{\otimes}^*$. С.~Санкаран и С.~Селесник в~\cite{b40} пользовались инъективным тензорным произведением, но не на категории $C^*$-алгебр, а на категории банаховых пространств, и вводили $C^*$-алгебры Хопфа как такие алгебры Хопфа над указанной категорией, которые являются ещё и $C^*$-алгебрами. Как видим, даже в случае $C^*$-алгебр возможны, в принципе, совершенно различные подходы. Основные результаты главы 4 показывают, что это различие --- чисто внешнее.

В §1 определяется и изучается проективный способ топологизации тензорного произведения $LC^*$-алгебр. Напомним, что проективная локально выпуклая топология на тензорном произведении $E_1\otimes E_2$ локально выпуклых пространств $E_1$, $E_2$ определяется как сильнейшая из локально выпуклых топологий, относительно которых каноническое билинейное отображение
$$
  \otimes: E_1\times E_2 \to E_1\otimes E_2
$$
непрерывно. Эта топология обозначается через  $\pi$, наделённое ею пространство \mbox{$E_1\otimes E_2$}~--- через $E\underset{\pi}{\otimes} F$, его пополнение~--- через $E\mathbin{\hat{\underset{\pi}{\otimes}}} E_2$.

А.~Маллиос показал в \cite{b24}, что если $A_1$, $A_2$ --- $LMC$-алгебры, то и $A_1\otimes A_2$, $A_1\mathbin{\hat{\underset{\pi}{\otimes}}} A_2$~--- тоже. Вопрос~--- в том, должна ли $\pi$ быть $LC^*$-топологией, если $A_1$, $A_2$ --- $LC^*$-алгебры. В самом конце §2 показывается, что ответ~--- отрицательный. Тем не менее, для любых $LC^*$-алгебр $A_1$, $A_2$ на тензорном произведении $A_1\otimes A_2$ инволютивных алгебр $A_1$, $A_2$ существует сильнейшая $LC^*$-топология, относительно которой каноническое отображение $\otimes$ непрерывно. Эту топологию мы называем проективной $LC^*$-топологией и обозначаем через $\pi^*$ (в отличие от проективной локально выпуклой топологии $\pi$). Наделённую этой топологией алгебру $A_1\otimes A_2$ мы обозначаем через $A_1\underset{\ \pi^*}{\otimes} A_2$, а её отделимое пополнение~--- через $A_1\mathbin{\underset{\ \pi^*}{\bar\otimes}}A_2$.

Мы показываем, что $A_1\mathbin{\underset{\ \pi^*}{\otimes}} A_2$ является прямой суммой $A_1$ и $A_2$ в категории $LC^*$-алгебр, $A_1{\mathbin{\underset{\ \pi^*}{\bar\otimes}}} A_2$~--- прямой суммой $A_1$ и $A_2$ в категории  $K^*$-алгебр (следствие (l.9)).

В §2 изучаются инъективные тензорные произведения. Мы доказываем, что для любых отделимых $LC^*$-алгебр $A_1$, $A_2$ инъективная локально выпуклая топология~$\eps$ (см. \cite{b37}, \cite[Т. 2]{b16} или начало §2) на инволютивной алгебре $A_1\otimes A_2$ является отделимой $LC^*$-топологией, так что $A_1\mathbin{\underset{\varepsilon}{\otimes}} A_2$ является отделимой  $LC^*$-алгеброй, а её пополнение $A_1\ehatotimes A_2$ --- $K^*$-алгеброй. Из результатов главы 3 затем выводится, что $A_1\ehatotimes A_2$ является прямой суммой $A_1$ и $A_2$ в категории $K^*$-алгебр. Таким образом, для любых $K^*$-алгебр $A_1$, $A_2$ алгебры $A_1\mathbin{\underset{\ \pi^*}{\bar\otimes}}A_2$ и $A_1\ehatotimes A_2$ топологически эрмитово изоморфны и, значит, топологии $\pi$ и $\eps$ на алгебре $A_1\otimes A_2$ совпадают. Тем самым, на алгебре $A_1\otimes A_2$ существует только одна $LC^*$-топология $\tau\ges\eps$, относительно которой каноническое отображение $\otimes$ непрерывно,~--- это топология $\eps$. Используя пример В.~Дитриха \cite{b11}, мы показываем, что существуют $K^*$-алгебры $A_1$, $A_2$, для которых топология $\pi$ на $A_1\otimes A_2$ отлична от топологии $\pi^*$ и, значит, $\pi$  не является $LC^*$-топологией.

\medskip
\ul{Глава 5: Двойственность для полугрупп}

В этой главе мы оказываемся, наконец, готовыми к получению основных теорем двойственности для компактологических и некоторых классов топологических полугрупп и групп.

В §1 вводятся основные в этой главе определения \mbox{$K^*$-бигебры}, $K^*$-ал\-гебры Хопфа и их морфизмов. $K^*$-бигеброй (соотв., $K^*$-алгеброй Хоп\-фа) мы называем $K^*$-алгебру~$A$, наделённую непрерывным эрмитовым морфизмом $\mu\colon A \to A\ehatotimes A$ (и непрерывными эрмитовыми морфизмами $\eps\colon A \to \Cbb$, $\sigma\colon A \to A$), для которых коммутативна диаграмма
$$
  \begin{diagram}
  \node{A}\arrow{s,t}{\mu}\arrow{e,t}{\mu}
  \node{A\ehatotimes A} \arrow{s,b}{\mu\ehatotimes I_A}\\
  \node{A\ehatotimes A} \arrow{e,b}{I_A\ehatotimes\mu}
  \node{A\ehatotimes A\ehatotimes A}
  \end{diagram}
$$
(и диаграммы
$$
  \begin{diagram}
  \node[2]{A}\arrow{sw,t}{a\mapsto 1\otimes a}\arrow{s,b}{\mu}
  \arrow{se,t}{a\mapsto a\otimes 1}\\
  \node{\Cbb\ehatotimes A}
  \node{A\ehatotimes A} \arrow{e,b}{\eps\ehatotimes I_A}
  \arrow{w,b}{I_A\ehatotimes\eps}
  \node{A\ehatotimes\Cbb,}
  \end{diagram}
$$
\medskip
$$
  \begin{diagram}
  \node{A}\node{\Cbb}\arrow{w,t}{}
  \node{A}\arrow{w,t}{\eps}\arrow{s,t}{\mu}\arrow{e,t}{\eps}
  \node{\Cbb}\arrow{e}{}\node{A}\\
  \node{A\ehatotimes A}\arrow{n,t}{m}
  \node[2]{A\ehatotimes A}
  \arrow[2]{w,b}{I_A\ehatotimes\sigma}\arrow[2]{e,b}{I_A\ehatotimes\sigma}
  \node[2]{A\ehatotimes A,}\arrow{n,b}{m}
  \end{diagram}
$$
где $m$~--- линейное отображение, канонически соответствующее
умножению на $A$, а стрелка \glqq $\mathbb{C}\rightarrow A$\grqq\ изображает единственный морфизм из $\mathbb{C}$ в $A$). Морфизмы $\mu$,
$\varepsilon$ и~$\sigma$ называются, соответственно коумножением,
коединицей и кообращением. Последние два, если они существуют, определяются по $A$ и $\mu$ однозначно. По этой причине, называя $K^*$-алгебру Хопфа, мы будем обычно указывать лишь $K^*$-алгебру и коумножение.

Морфизмом из $K^*$-бигебры $(A_1,\ \mu_1)$ в $K^*$-бигебру $(A_2,
\mu_2)$ мы называем каждый непрерывный эрмитов морфизм $u\colon
A_1\rightarrow A_2$, согласующийся с коумножениями $\mu_1$ и
$\mu_2$, то есть такой, для которого диаграмма
$$
  \begin{diagram}
  \node{A_1}\arrow{s,t}{\mu_1}\arrow{e,t}{u}
  \node{A_2} \arrow{s,b}{\mu_2}\\
  \node{A_1\ehatotimes A_1} \arrow{e,b}{u_1\ehatotimes u_2}
  \node{A_2\ehatotimes A_2}
  \end{diagram}
$$
коммутативна. Если при этом $(A_1,\ \mu_1)$ и $(A_2,\ \mu_2)$~---
$K^*$-алгебры Хопфа, то $u$ согласуется также с коединицей и
кообращением, так что в отдельном понятии морфизма для $K^*$-алгебр Хопфа нет надобности.

В §2 доказываются теоремы двойственности для компактологических полугрупп и групп. Напомним, что группа над категорией $\Ccal$, обладающей конечными произведениями и, значит, финальным объектом, скажем, $F$ , может быть определена \cite{b3} как объект $G$, наделённый морфизмами $m\colon G\times G \to G$  (\glqq умножение\grqq), $e\colon F \to G$ (\glqq единица\grqq) и $s\colon G \to G$  (\glqq обращение\grqq), для которых коммутативны диаграммы
$$
  \begin{diagram}
  \node{G\times G\times G}\arrow{s,t}{m\times I_G}
  \arrow{e,t}{I_G\times e}
  \node{G\times G} \arrow{s,b}{m}\\
  \node{G\times G} \arrow{e,b}{m}
  \node{G,}
  \end{diagram}
$$
$$
  \begin{diagram}
  \node{F\times G}\arrow{e,t}{e\times I_G}
  \arrow{se,b}{pr_2}
  \node{G\times G} \arrow{s,b}{m}
  \node{G\times F} \arrow{sw,b}{pr_1}\arrow{w,t}{I_G\times e}\\
  \node[2]{G,}
  \end{diagram}
$$
$$
  \begin{diagram}
  \node{G\times G}\arrow[2]{e,t}{s\times I_G}
  \node[2]{G\times G} \arrow{s,t}{m}
  \node[2]{G\times G} \arrow[2]{w,b}{I_G\times s}\\
  \node{G}\arrow{n,t}{\Delta}\arrow{e,t}{}\node{F}\arrow{e,t}{e}
  \node{G}\node{F}\arrow{w,t}{e}
  \node{G.}\arrow{w,t}{}\arrow{n,b}{\Delta}
  \end{diagram}
$$
Опуская всё, что касается единицы и обращения, получаем понятие полугруппы над категорией. Единица и обращение, если они существуют, определяются по объекту и умножению на нём однозначно. Морфизм $u\colon G_1\to G_2$ категории $\Ccal$ называется морфизмом из полугруппы $(G_1,\ m_1)$ в полугруппу $(G_2,\ m_2)$, если он согласуется с умножениями $m_1$ и $M_2$, т.~е. если коммутативна диаграмма
$$
  \begin{diagram}
  \node{G_1}\arrow{e,t}{u}\node{G_2}\\
  \node{G_1\times G_1} \arrow{n,t}{m_1}\arrow{e,t}{u\times u}
  \node{G_2\times G_2.} \arrow{n,b}{m_2}
  \end{diagram}
$$
Если при этом $(G_1,\ m_1)$ и $(G_2,\ m_2)$ --- группы, то $u$ согласуется также с единицей и обращением.

Полугруппы, соотв., группы над категорией регулярных компактологических пространств мы называем, соответственно, компактологическими полугруппами и группами. Основной результат §2, теорема (2.5), гласит: категория, дуальная к категории компактологических полугрупп, соотв., групп, эквивалентна категории  $K^*$-бигебр, соотв.,  $K^*$-алгебр Хопфа.

В §3 доказываются теоремы двойственности для категорий полугрупп и групп над категориями топологических пространств, рассматривавшимися в главе~3. Мы приведём здесь лишь теорему, относящуюся к категории всех отделимых локально компактных полугрупп, соотв., групп. $K^*$-бигебру $(A,\ \mu)$ будем называть  $lc^*$-бигеброй, если $A$ --- $lc^*$-алгебра. Основной результат главы 5, теорема (3.8), гласит:\\

\parbox{.94\textwidth}{\em%
категория, дуальная к категории всех отделимых локально компактных полугрупп $($соотв., групп$)$ эквивалентна категории  $lc^*$-бигебр $($соотв.,  $lc^*$-алгебр Хопфа$)$.
}\\

Наконец, в §4 мы показываем, как выглядит двойственность Понтрягина в рамках построенной выше теории.

Заканчивая перечисление результатов диссертации, отметим, что в ней впервые:

описаны а) все вообще и б) все согласующиеся с данной двойственностью $LC^*$-топологии на коммутативных инволютивных алгебрах с единицей, отдельно рассмотрены сильнейшие и слабейшие среди всех таких топологий;

введены $C^*$-бочечные $LC^*$-алгебры, для которых доказаны обобщение теоремы Нахбина-Сироты о бочечных пространствах непрерывных функций, аналог теоремы Комуры и обобщение теоремы Бухвальтера-Шмета об ассоциированных бочечных пространствах;

найдены категории $K^*$-алгебр, дуальные к ряду категорий топологических пространств, включая категорию всех $k_\Cbb$-пространств, вполне регулярных $k$-пространств, отделимых локально компактных пространств и другие;

доказано, что проективная $LC^*$-топология тензорного произведения $A_1\otimes A_2$ отделимых коммутативных $LC^*$-алгебр с единицей $A_1$, $A_2$ (т.~е. сильнейшая $LC^*$-топология, относительно которой каноническое билинейное отображение $\otimes$ непрерывно), совпадает с обычной инъективной топологией $\eps$, в силу чего на $A_1\otimes A_2$ существует только одна $LC^*$-топология $\tau\ges\eps$, относительно которой $\otimes$ непрерывно, --- это топология~$\eps$;

получены теоремы двойственности для категории всех компактологических и ряда категорий топологических полугрупп и групп, включая категорию всех отделимых локально компактных полугрупп и групп.

Результаты диссертации докладывались и обсуждались на семинарах Д.~А.~Райкова по теории категорий в МГПИ им. В.~И.~Ленина (1971--1972~гг.) и по теории топологических векторных пространств в МГУ (1975~г.), на научных конференциях преподавателей Рязанского пединститута (1977--1982 гг.) , на IX Всесоюзном симпозиуме по теории групп (Москва, 1984~г.) и на семинаре Л.~Я.~Куликова по алгебре в МГПИ им. В.~И.~Ленина (1987~г.).

По теме диссертации автором опубликовано 7 научных работ \cite{b26}--\cite{b32}, выполненных без соавторов.

Диссертация является самостоятельным исследованием автора. Все полученные автором результаты являются новыми.


\eject

\thispagestyle{empty}
\mbox{}


\hfill\parbox{0.65\textwidth}{\large\textit{Автор с неизменной любовью и бесконечной благодарностью вспоминает своего научного руководителя доктора физико-математических наук, профессора Д.~А.~Райкова. Всё, что заслуживает внимания в этом исследовании, автор посвящает его светлой памяти.}}
\vfill

\mbox{}

\chapter{Объекты и методы исследования}

\subsection*{\S 1}
Наши обозначения и терминология в основном стандартны. Дополнение к множеству $X$ обозначается через $^cX$. Запись \glqq$f\colon X \to Y$\grqq\ означает, что $f$ есть функция из $X$ в $Y$, запись
$$
  f\colon\begin{cases}
     X\to Y;\\
     \mathcal{A}
  \end{cases}~\mbox{---}
$$
что $f$ есть функция из $X$ в $Y$, задаваемая правилом $\mathcal A$. Тождественная функция из множества $X$ в себя обозначается через
$I_X$, значение функции $f$ на элементе $x$~--- через $f(x)$ или
$\langle f,\ x\rangle$. Образ множества $X$ относительно функции $f$
обозначается через $f[X]$, множество всех функций из $X$ в $Y$~---
через $Y^X$. Композиция $^{\begin{CD}\cdot @>f>> \cdot @>g>> \cdot
\end{CD}}$ обозначается через $g\circ f$.

Предполагаются известными основные сведения из теории категорий и
функторов \cite{b3}, \cite{b23}, общей топологии \cite{b4}, \cite{b14} и теории
топологических векторных пространств \cite{b37}, \cite{b16}. Некоторые из них
напоминаются по ходу изложения. Основные определения и факты теории
топологических инволютивных алгебр напоминаются в этой главе.

\subsection*{\S 2} Относительно алгебр мы придерживаемся терминологии и обозначений из~\cite{b5}. Под алгеброй всюду в настоящей работе понимается коммутативная ассоциативная комплексная алгебра с единицей~$e$. Множество всех максимальных идеалов алгебры~$A$ обозначается через~$\mathcal{M}(A)$, множество всех идеалов коразмерности~1 (т.~е. идеалов~$I$, для которых факторалгебра~$A/I$ изоморфна полю~$\Cnum$ комплексных чисел),~--- через $\mathcal{M}_{\Cnum}(A)$. Алгебра $A$ называется полупростой, если пересечение всех элементов множества $\mathcal{M}(A)$  состоит лишь из нуля, $\bigcap \mathcal{M}(A) = \{0\}$, и $\Cnum$-полупростой, если аналогичным образом $\bigcap \mathcal{M}_{\Cnum}(A) =\{0\}$.

Морфизмами алгебр называются сохраняющие единицу гомоморфизмы,
характерами --- их морфизмы в $\Cnum$. Множество всех характеров
алгебры $A$ обозначается через $\mathrm{X}(A)$. Наделённое
топологией поточечной сходимости, оно становится тихоновским (т.~е.
вполне регулярным отделимым) топологическим пространством и
называется пространством характеров алгебры $A$.

Для любого морфизма $h$ из алгебры $A$ в алгебру $A'$ через
$\mathrm{X}(h)$ обозначается отображение $\chi'\mapsto \chi'\circ h$
из множества $\mathrm{X}'(A)$ в множество $\mathrm{X}(A)$. Это
отображение непрерывно относительно указанных топологий. Если $h$
--- тождественный морфизм алгебры $A$, то $\mathrm{X}(h)$ --- тождественное отображение на пространстве $\mathrm{X}(A)$: $\mathrm{X}(I_A) = I_{\mathrm{X}(A)}$. Если $h$ и $k$ --- морфизмы из $A$ в $A'$ и из $A'$ в $A''$ соответственно, то $\mathrm{X}(k\circ h) = \mathrm{X}(h)\circ \mathrm{X}(k)$. Таким образом, $\mathrm{X}$ ---
контравариантный функтор из категории алгебр в категорию тихоновских
пространств.

\subsection*{\S 3} Для любого топологического пространства $T$ множество $C(T)$ всех непрерывных комплекснозначных функций на $T$ является алгеброй относительно определяемых поточечно операций сложения, умножения и умножения на скаляр. Для любого непрерывного отображения $\map{u}{T}{T'}$ через $C(u)$ обозначается отображение $\vp' \mapsto \vp'\circ u$ из $C(T')$ в $C(T)$. Это --- морфизм алгебр. Если $u$ --- тождественное отображение пространства $T$, то $C(u)$ --- тождественный морфизм алгебры $C(T)$:  $C(I_T) = I_{C(T)}$. Если $u$ и $v$ --- непрерывные отображения из $T$ в $T'$ и из $T'$ в $T''$ соответственно, то $C(v\circ u) = C(u)\circ C(v)$. Тем самым $C$ --- контравариантный функтор из категории топологических пространств в категорию алгебр.

\subsection*{\S 4}
Пусть $T$ --- топологическое пространство. Для каждого $t\in T$ через $\dl_T(t)$ или просто $\dl(t)$ обозначается функция $\vp\mapsto \vp(t)$ из $C(T)$ в $\Cnum$, называемая преобразованием Дирака элемента $t$. Это --- характер алгебры $C(T)$. Отображение $t\mapsto \dl_T(t)$ из $T$ в $\Chi C(T)$ ($ = \Chi (C(T))$) называется преобразованием Дирака пространства $T$ (и обозначается через $\dl_T$). Оно инъективно тогда и только тогда, когда $C(T')$
различает точки пространства $T$, и является гомеоморфизмом на образ
тогда и только тогда, когда пространство $T$~--- тихоновское.

Для любого непрерывного отображения $\map{u}{T}{T'}$ диаграмма
$$
  \begin{CD}
   T@>u>>T'\\
  @V\dl_T VV @VV\dl_{T'}V\\
  \Chi C(T)@>\Chi C(u)>>T'
  \end{CD}
$$
коммутативна. Это означает, что соответствие $T \mapsto \dl_T$
является функторным морфизмом из тождественного функтора категории
топологических пространств в функтор~$\Chi C$.

\subsection*{\S 5}
Пусть $A$ --- алгебра. Для каждого $x\in A$ через $\gel_A(x)$,
или просто $\gel x$, обозначается функция $\chi \mapsto \chi(x)$ из
$\Chi(A)$ в $\Cbb$, называемая преобразованием Гельфанда
элемента $x$. Все такие функции на пространстве $\Chi(A)$
непрерывны. Более того, топологическое пространство $\Chi(A)$
таково, что его топология является слабейшей из всех топологий,
относительно которых каждая функция $\gel x$ с $x\in A$  непрерывна.
Отображение $\gel_A\colon x \mapsto \gel_A x$   является морфизмом
алгебры $A$ в алгебру $C\Chi(A)$ $(= C(\Chi(A)))$ и называется
преобразованием Гельфанда алгебры $A$. Алгебра $A$
$\Cbb$-полупроста тогда и только тогда, когда отображение $\gel_A$
инъективно. Для любого морфизма $h$ из алгебры $A$ в алгебру $A'$
диаграмма
$$
  \begin{CD}
   A @>h>> A'\\
   @V\gel_A VV @V\gel_{A'} VV\\
   C\Chi(A) @>C\Chi(A)>> C\Chi(A')
  \end{CD}
$$
коммутативна. Это означает, что соответствие $A \mapsto \gel_A$
является функторным морфизмом из тождественного функтора категории
алгебр в функтор $C\Chi$.

\subsection*{\S 6} Для любой алгебры $A$ и любого топологического пространства $T$ диаграммы
$$
  \xymatrix{
   & \Chi C\Chi(A)\ar[dr]^{\Chi(\gel_A)} & \\
   \Chi(A)\ar[rr]^{I_{\Chi(A)}}\ar[ur]^{\dl_{\Chi(A)}} && \Chi(A)\,,
   }\qquad
   \xymatrix{
   & C\Chi C(T)\ar[dr]^{C(\dl_T)} & \\
   C(T)\ar[rr]^{I_{C(T)}}\ar[ur]^{\gel_{C(T)}} && \Chi(A)
   }
$$
коммутативны. Иначе говоря, контравариантные функторы $\Chi$ и $C$
сопряжены справа.

\subsection*{\S 7} Преднорма $p$ на алгебре $A$ называется субмультипликативной, если она удовлетворяет условию
$$
  p(xy) \les p(x)p(y)\qquad \text{при всех } x, y\in A.\leqno(sm)
$$
Множество всех субмультипликативных преднорм на алгебре $A$
обозначается через ${\mathcal P}(A)$ . Алгебра, наделённая
субмультипликативной \mbox{(пред-)}нормой, называется (пред-)нормированной
алгеброй. Полная нормированная алгебра называется банаховой
алгеброй. Норму в нормированной алгебре всегда можно выбрать таким
образом, чтобы $\|e\| = 1$, поэтому обычно предполагают это условие
выполненным.

Пусть $p$ --- субмультипликативная преднорма на алгебре $A$,
$N(p)$~--- её ядро. Тогда $N(p)$ является идеалом в $A$ и условием
$p^{*} (x + N(p)) = p(x)$, $x\in A$, корректно определяется норма
$p^{*}$ на факторалгебре $A/N(p)$. Банахова алгебра, получаемая
пополнением алгебры $A/N(p)$ относительно указанной нормы,
называется $p$-фактором алгебры $A$ и обозначается через $A_p$.
Каноническое отображение $\pi_p\colon A\to A_p$ является морфизмом
алгебр, а также непрерывно и открыто.

\subsection*{\S 8} Топологической алгеброй называется множество, наделённое структурой алгебры и топологией так, что при этом операции сложения и умножения на скаляры совместно непрерывны, а операция умножения элементов алгебры --- раздельно непрерывна. Множество всех замкнутых максимальных идеалов топологической алгебры $A$ обозначается через ${\mathcal M}_c (A)$. Топологическая алгебра называется строго полупростой, если пересечение всех её замкнутых максимальных идеалов состоит лишь из нуля: ${\mathcal M}_c (A) = \{0\}$.

Множество всех непрерывных характеров алгебры $A$ обозначается через
$\Chi_c(A)$. Наделённое топологией поточечной сходимости, оно
становится тихоновским пространством, которое является
подпространством пространства $\Chi(A)$ и называется пространством
непрерывных характеров алгебры $A$.

Для любого непрерывного морфизма $h$ из топологической алгебры $A$ в
топологическую алгебру $A'$ отображение $\Chi (h)$ переводит $\Chi_c
(A')$ в $\Chi_c(A)$. Ограничение отображения $\Chi(h)$ множествами
$\Chi_c (A')$ и $\Chi_c(A)$ обозначается через $\Chi_c (h)$:
$$
  \Chi_c (h)\colon
\begin{cases}
  \Chi_c(A')\to \Chi_c(A);\\
  \chi'\mapsto \chi'\circ h.
\end{cases}
$$
Так же, как и для $\Chi$, $\Chi_c(I_A) = I_{\Chi_c(A)}$ и $\Chi_c
(k\circ h) = \Chi_c(h)\circ \Chi_c(k)$. Тем самым $\Chi_c$ есть
контравариантный функтор из категории топологических алгебр в
категорию тихоновских пространств.

\subsection*{\S 9}
Для любого топологического пространства $T$ алгебра $C(T)$, наделённая топологией равномерной сходимости на всех компактах из $T$, является топологической алгеброй и обозначается через $C_{co}(T)$. Для любого непрерывного отображения $u\colon T \to T'$ отображение $C(u)$ является непрерывным морфизмом из топологической алгебры $C_{co}(T')$ в топологическую алгебру $C_{co}(T)$. Рассматриваемое как морфизм топологических алгебр, оно иногда обозначается через $C_{co}(u)$. Тем самым получаем контравариантный функтор $C_{co}$ из категории топологических пространств в категорию топологических алгебр.

\subsection*{\S 10}
Пусть $T$ --- топологическое пространство. Для любого $t\in T$ характер $\dl_T(t)$ очевидным образом непрерывен относительно топологии на $C_{co}(T)$, так что $\dl_T$ отображает $T$ в $\Chi_c
C_{co}(T)$. Тем самым соответствие $T\mapsto \dl_T$ является также
функторным морфизмом из тождественного функтора категории
топологических пространств в функтор $\Chi_c C_{co}$.

\subsection*{\S 11}
Пусть $A$ --- топологическая алгебра. Для любого $x\in A$
сужение функции $\gel_A x$ на $X_c C(A)$ по-прежнему обозначается
через $\gel_A x$ и называется преобразованием Гельфанда элемента
$x$. Отображение $x\mapsto \gel_A x\colon A \to CX_c(A)$ по-прежнему
обозначается через $\gel_A$ и называется преобразованием Гельфанда
топологической алгебры~$A$. Подчеркнём, что отображение $\gel_A$ из
$A$ в $C_{co} \Chi_c(A)$ не обязательно непрерывно, поэтому
соответствие $A\mapsto \gel_A$ не является функторным морфизмом из
тождественного функтора категории топологических алгебр в функтор
$C_{co} \Chi_c$.

\subsection*{\S 12}
Пусть $A$ --- топологическая алгебра, $Q$ — множество всех её
обратимых элементов. Алгебра $A$ называется алгеброй с непрерывным
обращением, если отображение $x \mapsto x^{-1}\colon Q \to Q$
непрерывно. Подчеркнём, что алгебра с непрерывным обращением в
указанном смысле не является, вообще говоря, кольцом с непрерывным
обратным в смысле \cite{b33} (т.~е. $Q$-алгеброй в смысле \cite{b1}), ибо мы не
требуем, чтобы единица алгебры $A$ входила в $Q$ с некоторой
окрестностью.

\subsection*{\S 13}
Топология на алгебре, порождаемая каким-либо семейством
субмультипликативных преднорм, называется локально мультипликативно
выпуклой, локально $m$-выпуклой или  $LMC$-топологией \cite{b1}, \cite{b5}, \cite{b25}. Алгебра, наделённая $LMC$-топологией, называется $LMC$-алгеброй. Для любого топологического пространства $T$ топологическая алгебра $C_{co}(T)$ является $LMC$-алгеброй. Для
любой топологической алгебры $A$ через $\mathcal{P}_c(A)$
обозначается совокупность всех субмультипликативных непрерывных
преднорм на $A$.

\begin{thm}[(LMC) Основные свойства LMС-алгебр (1)]\label{LMC}\mbox{}\\
 Пусть $A$ --- $LMC$-алгебра. Тогда:
 \begin{enumerate}[$1^\circ$]
  \item $A$ --- топологическая алгебра.
  \item Умножение в $A$ совместно непрерывно.
  \item $A$ --- алгебра с непрерывным обращением.
  \item Если $A$ --- поле, то $A = \Cnum\cdot e$.
  \item Отображение $\chi \mapsto \Ker \chi$ является биекцией множества
  $\Chi_c(A)$ на множество $\mathcal{M}_c(A)$.
  \item Рассмотрим условия:

  \ \ а) $A$ полупроста;

  \ \ б) $A$ строго полупроста;

  \ \ в) непрерывные характеры различают точки $A$;

  \ \ г) отображение $\gel_A\colon A\to C\Chi_c(A)$ инъективно.

Условия б), в), г) равносильны и влекут условие а), а если алгебра
$A$ полна, то все четыре условия равносильны.
 \end{enumerate}
\end{thm}

\subsection*{\S 14}
Инволюция на алгебре $A$ --- это отображение $x\mapsto x^*$
из $A$ в $A$ такое, что
\begin{enumerate}[($*1)$\quad]
\item  $(x^*)^* = x$;
\item  $(x + y)^* = x^*+y^*$;
\item  $(x\cdot y)^* = y ^*\cdot x^*$;
\item  $(\lm\cdot x)^* = \ol{\lm}\cdot x^*$
\end{enumerate}
при всех $x,\ y\in A$, $\lm\in\Cnum$ (через $\ol{\lm}$ обозначается
число, сопряжённое к $\lm$). Наделённая инволюцией алгебра
называется инволютивной алгеброй, $*$-алгеброй или даже просто
алгеброй, если из контекста ясно, что она инволютивна. Поле $\Cnum$
является инволютивной алгеброй относительно стандартной инволюции
$\lm\mapsto\ol{\lm}$.

Элементы или подмножества инволютивной алгебры, устойчивые
относительно инволюции, называются самосопряжёнными. Самосопряжённые
идеалы называются также $*$-идеалами, самосопряжённые элементы ---
эрмитовыми. Множество всех максимальных самосопряжённых идеалов
инволютивной алгебры $A$ обозначается через $\mathcal{M}^* (A )$.
Алгебра $A$ называется $*$-полупростой, если $\bigcap\mathcal{M}^*
(A) = \{0\}$.

Морфизмы инволютивных алгебр --- это морфизмы алгебр, сохраняющие
инволюцию. Их называют также эрмитовыми морфизмами. В частности,
эрмитовы характеры инволютивной алгебры $A$ --- это те $\chi \in
\Chi(A)$, для которых $\chi(x^*) = \ol{\chi(x)}$ при всех $x\in A$.
Множество всех эрмитовых характеров алгебры $A$ обозначается через
$\Chi^*(A)$. Наделённое топологией поточечной сходимости, оно
становится тихоновским пространством, которое является
подпространством пространства $\Chi(A)$ и называется пространством
эрмитовых характеров алгебры $A$.

Для любого эрмитова морфизма $\map{h}{A}{A'}$ отображение $X(h)$
переводит $\Chi^*(A')$ в $\Chi^*(A)$. Ограничение отображения
$\Chi(h)$ множествами $\Chi^*(A')$ и $X^*(A)$ обозначается через
$X^*(h)$:
$$
  \Chi^*(h)\colon
  \begin{cases}
   \Chi^*(A') \to \Chi^*(A),\\
   \chi'\mapsto \chi'\circ h.
  \end{cases}
$$
Так же, как $\Chi(h)$ и ${\Chi}_c (h)$, отображение $\Chi^*(h)$
непрерывно и обладает свойствами: $\Chi^*(I_A) = I_{\Chi^*(A)}$,
$\Chi^*(k\circ h) = \Chi^*(h)\circ \Chi^*(k)$. Тем самым $\Chi^*$
является контравариантным функтором из категории инволютивных алгебр
в категорию тихоновских пространств.

\subsection*{\S 15}
Для любого топологического пространства $T$ на алгебре $C(T)$
имеется стандартная инволюция, определяемая условием $\vp^*(t) =
\ol{\vp(t)}$ для любых $\vp \in C(T)$, $t\in T$. Для каждого $t\in
T$ характер $\dl_T(t)$ алгебры $C(T)$ --- эрмитов, так что $\dl_T$
отображает $T$ в $\Chi^*C(T)$. Тем самым соответствие $T\mapsto
\dl_T$ является функторным морфизмом из тождественного функтора
категории топологических пространств в функтор $\Chi^*C$.

\subsection*{\S 16}
Пусть $A$ --- инволютивная алгебра, $x\in A$. Сужение функции
$\gel_A x$ на множество $\Chi^*(A)$ будет по-прежнему называться
преобразованием Гельфанда элемента $x$ и обозначаться через
$\gel_A$. Аналогично, отображение $x\mapsto \gel_A x$ из $A$ в
инволютивную алгебру $CX^*(A)$  будет называться преобразованием
Гельфанда инволютивной алгебры $A$ и обозначаться через $\gel_A$.
Для любой инволютивной алгебры $A$ так определённое преобразование
Гельфанда $A$ является эрмитовым морфизмом, а соответствие $A
\mapsto \gel_A$ --- функторным морфизмом из тождественного функтора
категории инволютивных алгебр в функтор $C\Chi^*$.

\begin{remark*}
В этом месте наш подход к изучению инволютивных алгебр отличается от
общепринятого. Обычно (см., например, \cite{b33}) рассматривают
пространство всех характеров (или, что равносильно, пространство
всех максимальных идеалов коразмерности l), отмечают, что
преобразование Гельфанда не является эрмитовым морфизмом, и для
преодоления этого затруднения вводят понятие симметрической
инволютивной алгебры, у которой, по определению, преобразование
Гельфанда является эрмитовым морфизмом. По нашему мнению это уводит
немного в сторону. Поскольку на инволютивной алгебре имеется по
сравнению с обычными алгебрами дополнительная структура~---
инволюция, постольку \glqq главными\grqq\ морфизмами инволютивных
алгебр следует считать те, которые сохраняют эту структуру, т.~е.
эрмитовы. Соответственно этому \glqq главными\grqq\ характерами
будут эрмитовы характеры, а \glqq главным\grqq\ структурным
пространством --- пространство эрмитовых характеров. Именно такой
подход позволяет наиболее естественным образом обобщать результаты о
$C^*$-алгебрах. Надобность в понятии симметрической инволютивной
алгебры при этом даже не возникает.
\end{remark*}

\subsection*{\S 17}
Для любой инволютивной алгебры $A$ и любого топологического пространства $T$ диаграммы
$$
  \xymatrix{
    & \Chi^* C\Chi^*(A)\ar[dr]^{\Chi^*(\gel_A)} & \\
    \Chi^*(A)\ar[ur]^{\dl_{\Chi^*(A)}}\ar[rr]^{I_{\Chi^*(A)}} &&
    \Chi^*(A),
  }\qquad
  \xymatrix{
    & C\Chi^*C(A)\ar[dr]^{C(\dl_T)} & \\
    C(T)\ar[ur]^{\gel_{C(T)}}\ar[rr]^{I_{C(T)}} &&
    C(T)
  }
$$
коммутативны. Иначе говоря, контравариантные функторы $\Chi^*$ и $C$
сопряжены справа.

\subsection*{§ 18}
Преднорма $p$ на инволютивной алгебре $A$ называется
$C^*$-пред\-нор\-мой, если она удовлетворяет условию
$$
  p(xx^*) = [p(x)]^2\qquad   \text{при всех } x\in A.
  \leqno{\mathbf{(C^*)}}
$$
Множество всех $C^*$-преднорм на инволютивной алгебре $A$
обозначается через $\mathcal{P}^*(A)$.

\begin{thm}[($C^*$-1) Основные свойства $C^*$-преднорм]\label{C*1}\hfill
\begin{enumerate}[$1^\circ$]
\item Пусть $p$ --- $C^*$-преднорма на инволютивной алгебре $A$. Тогда:

а) $p(e) = 1$;

б)  $p(x^*) = p(x)$      при всех $x\in A$;

в)  $p$ субмультипликативна: \cite{b41};

г)  $A_p$ --- $C^*$-алгебра (см. ниже), а морфизм $\pi_p\colon A \to
A_p$ --- эрмитов.

\item Верхняя огибающая любого семейства $C^*$-преднорм является
$C^*$-преднормой.
\end{enumerate}
\end{thm}

$C^*$-преднорма, являющаяся нормой, называется $C^*$-нормой. Полная
$C^*$-нормированная алгебра называется $C^*$-алгеброй.

\begin{thm}[($C^*$-2) Основные свойства коммутативных $C^*$-алгебр \cite{b5}]\label{C*2}
Пусть $A$ --- коммутативная $C^*$-алгебра с единицей. Тогда:
\begin{enumerate}[$1^\circ$]
\item Каждый характер алгебры $A$ непрерывен и имеет норму, равную единице.

\item Каждый характер алгебры $A$ эрмитов, так что $\Chi(A) =
\Chi^*(A)$.

\item Пространство $\Chi(A)$ компактно.

\item Преобразование Гельфанда является изометрическим изоморфизмом алгебры $A$ на алгебру $C\Chi(A)$, наделённую супремум-нормой. В частности, характеры алгебры $A$ различают её точки и для любого $x\in A$
$$
  \|x\| = \sup\{\chi(x)\colon \chi\in\Chi(A)\}.
$$
\end{enumerate}
\end{thm}

\subsection*{\S 19} Топологическая алгебра с непрерывной инволюцией называется
топологической инволютивной алгеброй. Множество всех замкнутых
максимальных $*$-идеалов топологической инволютивной алгебры $A$
обозначается через $\mathcal{M}_c^*(A)$. Топологическая инволютивная
алгебра $A$ называется строго $*$-полупростой, если $\bigcap
\mathcal{M}_c^*(A) = \{0\}$.

Множество всех непрерывных эрмитовых характеров топологической
инволютивной алгебры $A$ обозначается через $\Chi^*(A)$. Наделённое
топологией поточечной сходимости, оно становится тихоновским
пространством, которое является пересечением подпространств
$\Chi^*(A)$ и $\Chi_c(A)$ пространства $\Chi(A)$ и называется
пространством непрерывных эрмитовых характеров алгебры $A$.

Для любого непрерывного эрмитова морфизма $\map{h}{A}{A'}$
отображение $\Chi(h)$ переводит $\Chi_c^*(A')$ в $\Chi_c^*(A)$.
Ограничение отображения $\Chi(h)$ множествами $\Chi_c^*(A')$ и
$\Chi_c^*(A)$ обозначается через $\Chi_c^*(h)$:
$$
  \Chi_c^*(h)\colon
  \begin{cases}
  \Chi_c^*(A') \to \Chi_c^*(A);\\
  \chi'\mapsto \chi'\circ h.
  \end{cases}
$$
Так же, как для $\Chi$ и $\Chi^*$, для $\Chi_c^*$ имеем:
$\Chi_c^*(I_A) = I_{\Chi_c^*(A)}$, $\Chi_c^*(k\circ h) =
\Chi_c^*(h)\circ \Chi_c^*(k)$. Тем самым $\Chi_c^*$ является
контравариантным функтором из категории топологических инволютивных
алгебр в категорию тихоновских пространств.

\subsection*{\S 20} Топология на инволютивной алгебре, определяемая
каким-либо семейством $C^*$-преднорм, называется $LC^*$-топологией.
Инволютивная алгебра, наделённая $LC^*$-топологией, называется
$LC^*$-алгеброй (star-алгеброй в \cite{b1}). Полная отделимая
$LC^*$-алгебра называется $K^*$-алгеброй (локально $C^*$-алгеброй~--- в \cite{b12}).

\begin{thm}[$(LC^*)$ Основные свойства $LC^*$-алгебр]\label{LC*}
Пусть $A$ --- $LC^*$-алгебра. Тогда:
\begin{enumerate}[$1^\circ$]
\item $A$ есть топологическая инволютивная алгебра и $LMC$-алгебра
(коротко: $LMC$-$*$-алгебра).

\item Каждый непрерывный характер на $A$ --- эрмитов, а каждый
замкнутый максимальный идеал --- самосопряжён: $\Chi_c^*(A)=
\Chi_c(A)$, $\mathcal{M}_c^*(A) = \mathcal{M}_c(A)$.

\item Отображение $\chi\mapsto \Ker \chi$ является биекцией множества
$\Chi_c^*(A)$ на множество $\mathcal{M}_c^*(A)$.

\item Рассмотрим условия:

а$*)$ $A$ $*$-полупроста;

б$*)$ $A$ строго $*$-полупроста;

в$*)$ непрерывные эрмитовы характеры алгебры $A$ различают её точки;

г$*)$ отображение $\gel_A\colon A \to C\Chi_c^*(A)$ инъективно;

д$*)$ $A$ отделима.

Каждое из условий б$*)$, в$*)$, г$*)$, д$*)$ равносильно любому из
условий б), в), г) п. 6° предложения $(LMC)$ и влечёт условие~a$*)$.
Если $A$ полна, то все девять условий, перечисленных в п.~6°
предложения $(LMC)$ и здесь, --- равносильны. Выделим особо:

\item $LC^*$-алгебра строго $*$-полупроста тогда и только тогда, когда
она отделима. $K^*$-алгебра строго $*$-полупроста.
\end{enumerate}
\end{thm}

Утверждение 1° доказано в \cite{b1}. 3° в силу 1° следует из 2° и п.~5° предложения $(LMC)$. 4° следует из 1°, 2° и п.~6° предложения
$(LMC)$. 5° следует из 4°. Докажем~2°.

\begin{proof}
Пусть $\chi$ --- непрерывный характер $LC^*$-алгебры $A$. Существуют
$p_1,\ \ldots,\ p_s \in \mathcal{P}_c^*(A)$ и $\eps > 0$ такие, что
$|\chi(x)| \les 1$ на $U = \Set{x\in A}{\max\limits_i p_i(x) \les
\eps}$. Но верхняя огибающая семейства $C^*$-преднорм является
$C^*$-преднормой ($C^*$-1, п.~$2^\circ$). Поэтому $p = \max\limits_i
p_i(x) \in \mathcal{P}_c^*(A)$. Таким образом, $|\chi(x)| \les 1$ на
$U = \ol{V}_{p,\ \eps}$, т.~е. $\chi$ непрерывен относительно
топологии, определяемой одной только $p$. В силу $p$-непрерывности
$\chi$ имеем $N(p) \sbs \Ker\chi$. Поэтому существует характер
$\tilde\chi$ на $A/N(p)$ такой, что $\chi = \tilde\chi \circ C$, где
$C$ --- канонический морфизм из $A$ на $A/N(p)$. Норма $\tilde p$,
определяемая на $A/N(p)$ преднормой $p$, порождает на $A/N(p)$
топологию, финальную относительно отображения $C$, поэтому
$\tilde\chi$ непрерывно на алгебре $A/N(p)$, наделённой этой нормой.
Значит, существует непрерывный характер $\chi_p$ на $A_p$,
являющийся продолжением характера $\tilde{\chi}$ : $\tilde{\chi} =
\chi_p \circ i$, где $i$ --- тождественное вложение $A/N(p)$
в~$A_p$. Отсюда получаем, что $\chi = \chi_p\circ i\circ C =
\chi_p\circ\pi_p$, где $\pi_p$ --- канонический морфизм из $A$
в~$A_p$. Но $p$ --- $C^*$-преднорма, значит $A_p$ --- $C^*$-алгебра,
а на $C^*$-алгебре каждый характер~--- эрмитов. Значит, ${\chi}_p$~--- эрмитов. Но и морфизм ${\pi}_p$~--- тоже эрмитов. Тем самым и характер $\chi$, как композиция эрмитовых
морфизмов, эрмитов. Это и требовалось доказать.
\end{proof}

\begin{cor}
Пусть $A$ --- $LC^*$-алгебра, $B$~--- топологическая инволютивная
алгебра. Если непрерывные эрмитовы характеры алгебры $B$ различают
её точки, то каждый непрерывный морфизм из $A$ в $B$~--- эрмитов.
\end{cor}

\begin{proof}
Пусть $\map{u}{A}{B}$ --- непрерывный морфизм. Возьмём произвольный
$x\in A$ и покажем, что $u(x^*)= u(x)^*$.

Для любого $\chi \in \Chi_c^*(B)$ композиция $\chi\circ u$ есть
непрерывный характер алгебры $A$. В силу п. 2° предыдущего
предложения он --- эрмитов. Значит
$$
  \chi(u(x^*)) = (\chi\circ u)(x^*) =  \ol{(\chi\circ u)(x)} =
  \ol{\chi(u(x))}.
$$
Но и $\chi(u(x)^*) = \ol{\chi(u(x))}$ (в силу эрмитовости $\chi$).
Значит, для любого $\chi\in\Chi_c^*(B)$ имеем $\chi(u(x^*)) =
\chi(u(x)^*)$. Поскольку элементы множества $\Chi_c^*(B)$ различают
точки алгебры $B$, отсюда следует, что $u(x^*) = u(x)^*$, а это и
требовалось.
\end{proof}

\subsection*{\S 21} Для любого топологического пространства $T$
стандартная инволюция на $C_{co}(T)$ непрерывна, так что $C_{co}(T)$
этой инволюцией является топологической инволютивной алгеброй.

Более того, топология алгебры $C_{co}(T)$ определяется семейством
преднорм вида~$p_K$:
$$
   p_K(\vp) = \sup\{|\vp(t)|: t \in K\},
$$
где $K$ пробегает совокупность $\mathcal{K}(T)$ всех компактных
подмножеств пространства~$T$, а каждая такая преднорма является
$C^*$-преднормой. Тем самым для любого топологического пространства
$T$ алгебра $C_{co}(T)$ является $LC^*$-алгеброй, а $C_{co}$ ---
контравариантным функтором из категории топологических пространств в
категорию $LC^*$-алгебр.

Напомним, что тихоновское пространство~$T$ называется
$k_{\R}$-про\-ст\-ран\-ст\-вом, если каждая вещественнозначная функция
на~$T$, сужения которой на все компактные подмножества
пространства~$T$ непрерывны, непрерывна на всём~$T$. Это равносильно
тому, что аналогичное условие выполняется для комплекснозначных
функций. Поскольку практически все наши рассмотрения производятся
над полем комплексных чисел, нам будет удобнее называть такие
пространства $k_{\Cnum}$-пространствами.

Алгебра $C_{co}(T)$ полна тогда и только тогда, когда пространство
$T$ является $k_{\R}$-пространством \cite{b49}, \cite{b1}. Тем самым сужение
$C_{co}$ на категорию $k_{\R}$-пространств является контравариантным
функтором в категорию $K^*$-алгебр.

\chapter{$LC^*$-топологии на инволютивных алгебрах}

\section{$C^*$-преднормы на инволютивных алгебрах}

В этом параграфе даётся описание $C^*$-преднорм и $LC^*$-топологий
на инволютивных алгебрах в терминах пространств эрмитовых характеров
этих алгебр.

Для любой преднормы $p$ на инволютивной алгебре $A$ через
$\Chi_p^*(A)$, или просто $\Chi_p^*$, если ясно, о какой $A$ идёт
речь, мы обозначаем множество всех эрмитовых характеров алгебры $A$,
непрерывных относительно топологии, определяемой одной лишь
преднормой $p$. Следующее предложение, несмотря на простоту,
является основой для получения многих дальнейших результатов. В
приводимом здесь виде оно, насколько нам известно, в литературе не
встречалось, хотя близкие результаты, конечно, отмечались. Например,
в~\cite{b12} (лемма 4.1, пп.~2 и 3) доказано утверждение а) для
случая, когда $A$ есть полная $LC^*$-алгебра и $p$ --- непрерывная
$C^*$-преднорма на ней.

\begin{prop}
Пусть $p$ --- $C^*$-преднорма на инволютивной алгебре $A$. Тогда:
\begin{itemize}
\item[а$)$] $\Chi^*_p(A)$ --- компакт, гомеоморфный $\Chi^*(A_p)$;
\item[б$)$] $p(x) = \sup\Set{\abs{\chi}}{\chi\in\Chi^*_p(A)}$ для всех
  $x\in A$.
\end{itemize}
\end{prop}

\begin{proof}
а) Поскольку $A_p$ --- $C^*$-алгебра, пространство $\Chi^*(A_p)$
компактно (гл. 1, предложение ($C^*$-2)). Поэтому достаточно
доказать, что $\Chi^*_p$ гомеоморфно $\Chi^*(A_p)$. Поскольку
каноническое отображение $\pi_p\colon A \to A_p$ является эрмитовым
морфизмом, имеем непрерывное отображение $\Chi^*(\pi_p)\colon
\Chi^*(A_p) \to \Chi^*(A)$. Покажем, что оно является биекцией
$\Chi^*(A_p)$ на $\Chi^*_p$.

Сначала установим инъективность. Пусть $\chi$, $\chi'\in\Chi^*(A_p)$
и $\Chi^*(\pi_p)(\chi) = \Chi^*(\pi_p)(\chi)$. Тогда $\chi\circ\pi_p
= \chi'\circ\pi_p$, так что $\chi$ и $\chi'$ совпадают на $\pi_p(A)
= A/N(p)$. Но $A/N(p)$ плотно в $A_p$, а $\chi$ и $\chi'$
непрерывны. Значит, они совпадают на всей $A_p$, т.~е. $\chi =
\chi'$, что и требовалось.

Затем, рассуждая так же, как при доказательстве п.~$2^\circ$
предложения $(LC^*)$ гл.~1, получим сюръективность отображения
$\Chi^*(\pi_{p})$. Для завершения доказательства после этого
останется только вспомнить, что непрерывное биективное отображение
компактного пространства на отделимое является гомеоморфизмом.

б) Пусть $\tilde p$ --- $C^*$-норма на $A_{p}$, порождаемая
преднормой $p$ (то есть продолжение на $A_{p}$ нормы $p^\bullet$,
см. гл. 1, $\S\S3$ и 18). Тогда, в силу а) и п.~$4^\circ$
предложения ($C^*$-2), для любого $x\in A$ будем иметь:
\begin{multline*}
p(x)=\tilde p(\pi_{p}(x))=
     \sup\{|\chi_{p}(\pi_{p}(x))|\colon\chi_{p}\in\Chi^*(A_{p})\}=\\=
     \sup\{|(\chi_{p}\circ \pi_{p})(x)|\colon\chi_{p}\in\Chi^*(A_{p})\}=
     \sup\{|\chi(x)|:\chi\in \Chi^*_{p}.
\end{multline*}
Это и требовалось.
\end{proof}

\begin{cor} Пусть $p$, $q$ --- $C^*$-преднормы на инволютивной
алгебре~$A$. Тогда:
\begin{enumerate}[$1^\circ$]
\item следующие условия равносильны:
\begin{itemize}
\item[а$)$] $q$ $p$-непрерывна;
\item[б$)$] $\Chi^*_{q}(A)\subset\Chi^*_{p}(A)$;
\item[в$)$] $q\les p$.
\end{itemize}
\item $q=p \iff \Chi^*_{q}(A)=\Chi^*_{p}(A)$
\end{enumerate}
\end{cor}

\begin{proof}
Утверждение $1^\circ$ легко доказывается по схеме
$$
  \text{а)} \Rightarrow \text{б)} \Rightarrow \text{в)} \Rightarrow \text{a)},
$$
а $2^\circ$ тотчас же следует из~$1^\circ$.
\end{proof}

\begin{cor}
Пусть $\chi\in\Chi^*(A)$. Тогда $\chi\in\Chi^*_{p}\iff|\chi|\les p$.
\end{cor}

\begin{proof}
Применяем (1.2).
\end{proof}

В соответствии с обычным пониманием минимальности, $C^*$-преднорму
$p$ (класса~$\C$) на инволютивной алгебре $A$ будем называть
минимальной (в классе $\C$), если она нетривиальна и для любой
нетривиальной $C^*$-преднормы $q$ (класса $\C$) на $A$ из
неравенства $q\les p$ следует равенство $q = p$.

\begin{prop}
$C^*$-преднорма $p$ на инволютивной алгебре $A$ минимальна тогда и
только тогда, когда существует эрмитов характер $\chi$ на $A$ такой,
что $\Chi^*_{p}=\{\chi\}$ и, значит, $p(x) = |\chi(x)|$ при
всех~$x\in A$.
\end{prop}

\begin{proof} В силу п.~$2^\circ$ следствия (1.2) $C^*$-преднорма $p$
на инволютивной алгебре $A$ минимальна тогда и только тогда, когда
$\Chi^*_{p}$ --- минимальное относительно включения непустое
подмножество в $\Chi^*(A)$. Но минимальными непустыми множествами
являются одноэлементные множества и только они. Отсюда и следует
утверждение.
\end{proof}

В силу предложения (1.1) каждой $C^*$-преднорме $p$ на инволютивной
алгебре $A$ отвечает подмножество $S$ множества $\Chi^*(A)$,
связанное с $p$ соотношением
$$
  p(a)=\sup\{|\chi(a)|\colon \chi\in S\} \text{ для любого } a\in A.
$$
Теперь мы обратим это соотношение.

\begin{dfn}
Пусть $S$ --- непустое подмножество пространства $\Chi^*(A)$. Для
любого $a\in A$ и любого $\varphi\in C\Chi^*(A)$ положим
\begin{align*}
  p^S(a)&=\sup\{|\chi(a)|: \chi\in S\},\\
  p_{S}(\varphi)&=\sup\{|\varphi(\chi)|: \chi\in S\}.
\end{align*}
Множество $S$ будем называть $A$-ограниченным, если $p^S(a)<\infty$
для любого $a\in A$.
\end{dfn}

\begin{prop}
Пусть $S$ --- $A$-ограниченное подмножество в $\Chi^*(A)$. Тогда:
\begin{itemize}
\item[а$)$] $p^S$ и $p_{S}$ --- $C^*$-преднормы на $A$ и $C\Chi^*(A)$
соответственно;
\item[б$)$] $\Chi^*_{p^S}(A) = \cl_{\Chi^*(A)}S$.
\end{itemize}
\end{prop}

\begin{proof} a) Следует из п.~$1^\circ$ предложения
(1.1) и п.~$2^\circ$ предложения $(C^*\mbox{-}1)$ главы 1.

б) ($\Rightarrow$) Пусть $\chi\in \Chi^*_{p^S}$. Покажем, что
$\chi\in \cl_{\Chi^*(A)}S$. Допустим противное. Тогда, поскольку
пространство $\Chi^*(A)$ вполне регулярно, существует $\varphi\in
C\Chi^*(A)$ такая, что $\varphi [S] = \{0\}$, а $\varphi(\chi)=1$.
В~силу теоремы Вейерштрасса-Стоуна~\cite{b14} $\gel_{A}[A]$ плотна в
$C_{co}\Chi^*(A)$ как содержащая константы и различающая точки
самосопряжённая (см. гл. 1, $\S 16$, замечание) подалгебра
функциональной алгебры $C_{co}\Chi^*(A)$. Поэтому для любого
$\alpha>0$ существует $x\in A$ такое, что
$$
  p_{\Chi^*_{p^S}}(\varphi-\gel_{A}x)<\alpha
$$
и, следовательно, $p^S(x)<\alpha$, а $|\chi(x)|>1-\alpha$. Возьмём
$\varepsilon=1/3$. В силу сказанного, для любого $\delta>0$
существует $x\in A$ такое, что
$$
  p^S(x)<\min\{\delta, 1/3\}, \text{ а } |\chi(x)|>1 - \min\{\delta, 1/3\},
$$
то есть такое, что $p^S(x)<\varepsilon$, а $|\chi(x)|>\varepsilon$.
Но это противоречит $p^S$-непрерывности характера $\chi$.

($\Leftarrow$) Пусть $\chi\in \cl_{\Chi^*(A)}S$. Тогда, для любой
$\varphi\in C\Chi^*(A)$,
$$
  \varphi(\chi)\in \varphi[\cl_{\Chi^*(A)}S]\subset \cl_{\Cnum}\varphi[S].
$$
В частности, для любого $x\in A$
$$
  \gel_{A}(x)(\chi)\in \gel_{A}(x)[\cl_{\Chi^*(A)}S]\subset \cl_{\Cnum}\gel_{A}(x)[S],
$$
то есть $\chi(x)\in \cl_{\Cnum}\{\gamma(x)\colon \gamma\in S\}$. Тем
самым, при всех $x\in A$,
$$
  |\chi(x)|\les \sup\{|\gamma(x)|\colon\gamma\in S\}=p^S(x).
$$
Это и означает, что $\chi\in \Chi^*_{p^S}(A)$. Предложение доказано.
\end{proof}

\begin{cor}
Замыкание в $\Chi^*(A)$ любого $A$-ограниченного подмножества
компактно.
\end{cor}

\begin{proof}
Действительно, в силу (1.6) $\cl_{\Chi^*(A)}S=\Chi^*_{p^S}$, а
$\Chi^*_{p^S}$, в силу (1.1), --- компактно.
\end{proof}

Топологическое пространство, в котором замыкание каждого
ограниченного подмножества компактно, называется
$\mu$-пространством~\cite{b7} или $NS$-пространством~\cite{b1}, то
есть пространством Нахбина-Сироты. Поскольку все ограниченные
подмножества пространства $\Chi^*(A)$ содержатся среди
$A$-ограниченных, следствие (1.7) показывает, что для любой
инволютивной алгебры $A$ пространство $\Chi^*(A)$ является
\mbox{$\mu$-пространством}. Однако на самом деле имеет место
большее.

Тихоновское пространство $T$ называется
(вещественно-компактным~\cite{b9}, или)
$Q$-пространством~\cite{b57}, если каждый вещественный характер
алгебры $C_{\R}(T)$ всех вещественнозначных непрерывных функций на
$T$ фиксирован. Иными словами, для любого сохраняющего единицу
гомоморфизма $\chi$ $\R$-алгебры $C_{\R}(T)$ в $\R$ существует $t\in
T$ такое, что для всех $\varphi\in C_{\R}(T)$ выполняется равенство
$\varphi(t)=\chi(\varphi)$. Каждому тихоновскому пространству $T$
отвечает $Q$-пространство $\upsilon T$, обладающее свойствами:
\begin{itemize}
\item[а)] $T$ --- плотное подпространство в $\upsilon T$;
\item[б)] каждая непрерывная функция $\varphi\colon T\to \R$
 обладает непрерывным продолжением на $\upsilon T$,
\end{itemize}
и называемое хьюиттовским (или $Q$-) расширением пространства $T$.

Бинц доказал в~\cite{b2}, что для любой вещественной коммутативной
алгебры с единицей пространство всех её вещественных характеров
является $Q$-пространством. Наш следующий результат показывает, что
аналогично утверждение имеет место и для пространств эрмитовых
характеров инволютивных алгебр.

\begin{thm}
Для любой инволютивной алгебры $A$ пространство $\Chi^*(A)$ всех её
эрмитовых характеров является $Q$-пространством.
\end{thm}

\begin{remark*}
Эту теорему можно доказать примерно так же, как доказывается
аналогичное утверждение для вещественных алгебр в~\cite{b2}. Мы
приводим другое доказательство, которое, как нам кажется, лучше
проясняет природу рассматриваемого явления.
\end{remark*}

\begin{proof} Пусть $A$ --- инволютивная алгебра. Множество всех её
эрмитовых характеров с топологией поточечной сходимости является
подпространством пространства $\Cnum^A$. Пространство $\Cnum$
гомеоморфно $\R^2$, $\R$ --- $Q$-пространство и произведение любого
семейства $Q$-пространств само является $Q$-пространством. Значит,
$\Cnum^A$~--- $Q$-пространство. Кроме того, замкнутое подмножество
$Q$-пространства само является $Q$-пространством. Поэтому, если мы
убедимся, что множество $\Chi^*(A)$ замкнуто в $\Cnum^A$, то теорема
будет доказана. Но это немедленно следует из непрерывности сложения,
умножения и комплексного сопряжения на $\Cnum$. Пусть $\chi\in
\Cnum^A$ является пределом сети $\chi_{\alpha}$ элементов множества
$\Chi^*(A)$ в $\Cnum^A$. Тогда, для любого $z\in A$,
$\chi(z)=\lim\limits_{\alpha}\chi_{\alpha}(z)$. Поэтому для любых
$x,\ y\in A$ имеем $\chi(x+y) =
\lim\limits_{\alpha}\chi_{\alpha}(x+y) =
\lim\limits_{\alpha}[\chi_{\alpha}(x)+\chi_{\alpha}(y)] =
\lim\limits_{\alpha}\chi_{\alpha}(x)+\lim\limits_{\alpha}\chi_{\alpha}(y)
= \chi(x)+\chi(y)$. Совершенно аналогично проверяются равенства
$\chi(\lambda x) = \lambda\chi(x)$, $\chi(xy)=\chi(x)\chi(y)$ и
$\chi(x^*)=\ol{\chi(x)}$. Таким образом, $\chi\in \Chi^*(A)$ и тем
самым, $\Chi^*(A)$ замкнуто в $\Cnum^A$. Теорема доказана.
\end{proof}

Нам понадобится также описание хьюиттовского расширения $\upsilon T$
тихоновского пространства $T$. Хьюитт~\cite{b57} строил $\upsilon T$
как пространство всех вещественных характеров алгебры $C_{\R}(T)$
всех вещественнозначных непрерывных функций на пространстве $T$.
Поскольку мы работаем с комплексными алгебрами, для нас эта
конструкция неудобна. Следующее предложение дает другую, более
пригодную для наших целей, реализацию пространства $\upsilon T$.

\begin{thm}
Для любого тихоновского пространства $T$ пространство $\Chi^*C(T)$
является $Q$-расширением пространства $T$.
\end{thm}

\begin{proof} Пусть $T$ --- тихоновское пространство.
Тогда (гл.~1,\ \S4) преобразование Дирака $\delta_{T}$ является
гомеоморфизмом на образ. Отождествим $T$ с $\delta_{T}[T]$
посредством отображения $\delta_{T}$ и докажем, что $\Chi^*C(T)$
обладает перечисленными выше свойствами а) и б).

а) $T$ плотно в $\Chi^*C(T)$. Действительно, пусть $\chi\in
\Chi^*C(T)$ и $V$ --- произвольная окрестность точки $\chi$ в
пространстве $\Chi^*C(T)$. По определению топологии в этом
пространстве существуют $\varphi_{1},\ \ldots,\ \varphi_{s}\in C(T)$
и $\varepsilon > 0$ такие, что
$$
  W = \{\chi'\in \Chi^*C(T)\colon\max\limits_{i}%
  |\chi'(\varphi_{i})-\chi(\varphi_{i})| <
  \varepsilon\}\subset V
$$
Покажем, что найдется $t\in T$, для которого $\delta_{T}(t)\in W$.

Рассмотрим функцию $\varphi\colon T\to C$, определенную условием
$$
  \varphi(t)=\sum |\varphi_{i}(t)-\chi(\varphi_{i})|^2
  \text{ для любого } t\in T.
$$
Ясно, что $\varphi\in C(T)$. Кроме того,
\begin{align*}
 \chi(\varphi)
 &= \chi\left(\sum|\varphi_{i}-\chi(\varphi_{i})e|^2\right)=\\
 &= \chi\left(\sum(\varphi_{i}-\chi(\varphi_{i})e)
    (\varphi_{i}-\chi(\varphi_{i}))^*\right)=\\
 &=\sum(\chi(\varphi_{i})-\chi(\varphi_{i}))
   \cdot\ol{(\chi(\varphi_{i})-\chi(\varphi_{i}))} =\\
 &= 0.
\end{align*}

Значит, существует $t\in T$ такое, что $\varphi(t)=0$.
Действительно, в противном случае нашлась бы $\psi\in C(T)$ такая,
что $\varphi\cdot\psi=e$. Тогда бы $\chi(\varphi)\cdot\chi(\psi)=1$
и, значит, $\chi(\varphi)\neq 0$, что, как мы только что видели,
невозможно. Но равенство $\varphi(t)=0$ означает, что
$\sum|\varphi_{i}(t)-\chi(\varphi_{i})|^2=0$, откуда
$\varphi_{i}(t)=\chi(\varphi_{i})$ при всех $i$ или, что  --- то же
самое, $\delta_{T}(t)(\varphi_{i})=\chi(\varphi_{i})$ при всех $i$.
Следовательно,
$\max\limits_{i}|\delta_{T}(\varphi_{i})-\chi(\varphi_{i})|=0<\varepsilon$
при всех и $\delta_{T}(t)\in W$. Это (с учетом сделанного выше
соглашения) и означает, что $T$ плотно в~$\Chi^*C(T)$.

б) Каждая непрерывная вещественнозначная функция на $T$ обладает
непрерывным вещественнозначным продолжением на все $\Chi^*C(T)$.
Действительно, преобразование Гельфанда $\gel_{C(T)}$ каждой
$\varphi\in C_{\R}(T)$ ставит в соответствие функцию
$\gel_{C(T)}\varphi$ из $C\Chi^*C(T)$. Композиция этой функции с
$\delta_{T}$ есть $\varphi$ (гл. 1, \S6), так что она является
непрерывным продолжением $\varphi$ на все $\Chi^*C(T)$. При этом для
любого $\chi\in \Chi^*C(T)$ имеем:
$$
  \ol{\gel_{C(T)}(\varphi)(\chi)}=
  \ol{\chi(\varphi)}=\chi(\varphi^*)=
  \chi(\varphi)=\gel_{C(T)}(\varphi)(\chi),
$$
то есть $\gel_{C(T)}(\varphi)\in C_{\R}\Chi^*C(T)$. Это и
требовалось.

Теорема доказана.
\end{proof}

Теперь мы готовы получить описание всех $LC^*$-топологий на
инволютивных алгебрах в терминах пространств эрмитовых характеров
этих алгебр.

Пусть $\mathcal B$ --- непустое семейство $A$-ограниченных
подмножеств пространства $\Chi^*(A)$. Топологией равномерной
сходимости на членах семейства $\mathcal B$ называется топология на $A$, порождаемая всеми преднормами $p^S$ с $S\in \mathcal B$. В силу п.~$1^\circ$ предложения (1.6) это есть $LC^*$-топология, а в силу п.~$2^\circ$ предложения (1.1) та же топология является топологией равномерной сходимости на компактах, являющихся замыканиями членов семейства $\mathcal B$. Наконец, предложение~(1.1) показывает, что каждая $LC^*$-топология является топологией равномерной сходимости на членах подходящего семейства компактов из пространства $\Chi^*(A)$. Это приводит к п.~$1^\circ$ следующей теоремы. Остальные очевидны.

\begin{thm}
Пусть $A$ --- инволютивная алгебра. Тогда:

$1^\circ$ Топология на $A$ является $LC^*$-топологией тогда и только
тогда, когда она является топологией равномерной сходимости на
членах некоторой совокупности компактных подмножеств пространства
$\Chi^*(A)$.

$2^\circ$ Если на $A$ существует хотя бы два эрмитовых характера, то на ней нет слабейшей $LC^*$-топологии.

$3^\circ$ Если на $A$ существует хотя бы один характер, то на ней
имеет сильнейшая $LC^*$-топология --- топология равномерной
сходимости на всех компактных подмножеств пространства $\Chi^*(A)$.
\end{thm}

Применяя теорему (1.9), получаем описание всех $LC^*$-топологий на алгебрах непрерывных комплекснозначных функций.

\begin{thm} Пусть $T$~--- тихоновское пространство.
Тогда:

$1^\circ$ Каждая $C^*$-преднорма на $C(T)$ имеет вид $p_{K}$ для
некоторого компакта $K$ из пространства $T$.

$2^\circ$ Каждая $LC^*$-топология на $C(T)$ является топологией
равномерной сходимости на некоторой совокупности компактов из
пространства $\upsilon T$.

$3^\circ$ Если $|T|\geqslant2$, то на $C(T)$ не существует слабейшей $LC^*$-топологии.

$4^\circ$ Для любого $T$ на $C(T)$ существует сильнейшая
$LC^*$-топология~--- топология равномерной сходимости на всех
компактных подмножествах из $\upsilon T$.
\end{thm}

\section[Непрерывные $C^*$-преднормы]{Непрерывные $C^*$-преднормы на топологических инволютивных алгебрах}

B этом параграфе мы изучим расположение в пространстве $\Chi^*(A)$ множеств $\Chi^*_{p}$ для непрерывных $C^*$-преднорм $p$ на
топологических инволютивных алгебрах $A$ и зависимость между
топологией алгебры $A$ и пространством её непрерывных эрмитовых
характеров.

Начнём с двух предложений, которые, хотя и не встречаются в явном
виде в литературе, всё же могут считаться "<хорошо известными">.

\begin{prop}
Пусть $A$ --- топологическая инволютивная алгебра. Для того, чтобы $C^*$-преднорма $p$ на $A$ была непрерывна, необходимо и достаточно, чтобы множество $\Chi^*_{p}$ содержалось в $\Chi^*_c(A)$ и было равностепенно непрерывно.
\end{prop}

\begin{proof} Пусть $C^*$-преднорма $p$ на топологическом инволютивной алгебре $A$ непрерывна. Тогда и каждый характер из множества $\Chi^*_{p}$ непрерывен, так что $\Chi^*_{p}\subset \Chi^*_{c}(A)$. В силу непрерывности $p$ множество $\overline{V}_{p} = \{x\in A\colon p(x)\leq 1\} $ является окрестностью нуля в $A$. Поэтому поляра $\overline{V}_p^\circ$ множества $\overline{V}_{p}$ в сопряжённом к $A$ пространстве равностепенно непрерывна. Но, в силу следствия
(1.3), $\Chi^*_{p}\subset \overline{V}^\circ_{p}$. Значит, и
$\Chi^*_{p}$ равностепенно непрерывно.

Обратно, пусть $\Chi^*_{p}\subset \Chi^*_{c}(A)$ и $\Chi^*_{p}$
равностепенно непрерывно. Тогда существует непрерывная преднорма $q$ на $A$ такая, что $\Chi^*_{p}\subset \overline{V}^\circ_{q}$, откуда
$$
  V_{p}=(\Chi^*_{p})^\circ \supset \overline{V}^{\circ\circ}_{q}
  \supset \overline{V}_{q}.
$$
Таким образом, $\overline{V}_{p}$ --- окрестность нуля в $A$, а это и означает, что преднорма $p$ непрерывна. Предложение доказано.
\end{proof}

\begin{cor}
Каждая $LC^*$-топология $t$ на инволютивной алгебре $A$ есть топология равномерной сходимости на всех равностепенно непрерывных компактах из $\Chi^*_{c}(A,\ t)$.
\end{cor}

\begin{dfn}
Пусть $A$ --- инволютивная алгебра и пусть
$\Gamma\subset\Chi^*(A)$. Будем говорить, что топология $t$ на
алгебре $A$ согласуется с двойственностью между $A$ и $\Gamma$, если $\Chi^*_{c}(A,\ t)=\Gamma$.
\end{dfn}

\begin{thm}
Пусть $A$ --- инволютивная алгебра. Тогда:

$1^\circ$ Для любого непустого множества $\Gamma\subset\Chi^*(A)$ на алгебре $A$ существует по меньшей мере одна $LC^*$-топология,
согласующаяся с двойственностью между $A$ и~$\Gamma$.

$2^\circ$ Каждая $LC^*$-топология на $A$, согласующаяся с
двойственностью между $A$ и~$\Gamma$, есть топология равномерной
сходимости на некоторой совокупности компактных подмножеств
пространства $\Chi^*(A)$, содержащихся в $\Gamma$ и покрывающих
$\Gamma$.

$3^\circ$ Среди всех $LC^*$-топологий на $A$, согласующихся с
двойственностью между $A$ и $\Gamma$, имеются слабейшая и
сильнейшая. Слабейшей является топология равномерной сходимости на всех конечных подмножествах множества $\Gamma$, сильнейшей~---
топология равномерной сходимости на всех компактных подмножествах
пространства $\Chi^*(A)$, содержащихся в $\Gamma$.

$4^\circ$ Для того, чтобы топология, согласующаяся с двойственностью между $A$ и $\Gamma$, была отделима, необходимо и достаточно, чтобы элементы множества $\Gamma$ различали точки множества $A$.

$5^\circ$ Минимальные непрерывные $C^*$-преднормы на $A$ --- одни и те же для всех $LC^*$-топологий, согласующихся с данной
двойственностью.
\end{thm}

\begin{proof}
Состоит в непосредственном применении результатов предыдущего
параграфа и предложения (2.1).
\end{proof}

\begin{prop}
Пусть $\tau_{1}$, $\tau_{2}$ --- $LC^*$-топологии на инволютивной
алгебре~$A$, пусть $\tau_{1}$ сильнее $\tau_{2}$ и пусть
$\Chi^*_{c}(A,\ \tau_{1})=\Chi^*_{c}(A,\ \tau_{2})$. Тогда:

$1^\circ$ Пополнение алгебры $A$ относительно топологии $\tau_{1}$ является подалгеброй в пополнении алгебры $A$ относительно топологии $\tau_{2}$.

$2^\circ$ Если алгебра $A$ полна относительно топологии $\tau_{2}$, то она полна и относительно топологии $\tau_{1}$.
\end{prop}

\begin{proof} $1^\circ$ Известно~\cite[\S~18, 4.(4), c]{b16}, что
пополнение топологического векторного пространства $(E,\ \tau_{1})$ является подпространством в пополнении топологического векторного пространства $(E,\ \tau_{2})$, если топология $\tau_{1}$ обладает базисом окрестностей нуля, замкнутых относительно топологии $\tau_{2}$. Установим это для нашего случая.

Пусть $U$ --- окрестность нуля относительно топологии $\tau_{1}$.
Поскольку $\tau_{1}$ --- $LC^*$-топология, существует
$C^*$-преднорма $p$ такая, что $\overline{V}_{p,\
\varepsilon}=\{x\in A\colon p(x)\leq \varepsilon \} \subset U$. По предложению (2.1) существует компакт $K\subset \Chi^*_{c}(A,\
\tau_{1})$ такой, что $p = p^K$ и, значит, $\overline{V}_{p,\
\varepsilon}=\overline{V}_{p^K,\ \varepsilon}$. Это показывает, что $\overline{V}_{p,\ \varepsilon}$ замкнуто относительно топологии $\tau_{0}$ поточечной сходимости на множестве $\Chi^*_{c}(A,\ \tau_{1})$. Но по условию $\Chi^*_{c}(A,\ \tau_{1})=\Chi^*_{c}(A,\ \tau_{2})$, а $\tau_{0}$ --- слабейшая из топологий $\tau$ на $A$, для которых $\Chi^*_{c}(A,\ \tau)=\Chi^*_{c}(A,\ \tau_{1})$. Значит, $\tau_{2}$ сильнее $\tau_{0}$, так что $\overline{V}_{p,\ \varepsilon}$ замкнуто и относительно топологии $\tau_{2}$. Таким образом, произвольно взятая окрестность нуля относительно топологии $\tau_{1}$ содержит окрестность, замкнутую относительно топологии $\tau_{2}$. Это и требовалось установить.

\smallskip
$2^\circ$ Пусть $A_{1}$ --- пополнение $A$ относительно $\tau_{1}$, $A_{2}$ --- пополнение $A$ относительно топологии $\tau_{2}$. Тогда $A\subset A_{1}$ и, на основании п.~$1^\circ$, $A_{1}\subset A_{2}$. Но в силу полноты~$A$ относительно $\tau_{2}$ имеет место равенство $A = A_{2}$. Значит, также $A = A_{1}$, а это и означает, что $A$ полна относительно $\tau_{2}$. Предложение доказано.
\end{proof}

\begin{thm}
Категории всех $LC^*$-алгебр, всех отделимых $LC^*$-алгебр и всех
$K^*$-алгебр являются полными рефлексивными подкатегориями категории всех топологических инволютивных алгебр с непустыми множествами непрерывных эрмитовых характеров.
\end{thm}

\begin{proof} Докажем утверждение, относящееся к $K^*$-алгебрам;
остальное доказывается аналогично (и проще). Пусть $A$~---
топологическая инволютивная алгебра и
$\Chi^*_{c}(A)\neq\varnothing$. Тогда и
$\mathcal{P}^*_{c}(A)\neq\varnothing$. Положим
$N=\bigcap\{N(p)\colon p\in \mathcal{P}^*_{c}(A)\}$ и рассмотрим
факторалгебру $A/N$. Правилом $p^\bullet(a+N)=p(a)$ для каждой $p\in \mathcal{P}^*_{c}(A)$ корректно определяется $C^*$-преднорма $p^\bullet$ на инволютивной алгебре $A/N$. Наделим $A/N$ топологией, определяемой всеми преднормами $p^\bullet$ с $p\in
\mathcal{P}^*_{c}(A)$, и обозначим через $\St(A)$ пополнение алгебры $A/N$ относительно этой топологии. Имеется каноническое морфизм $j\colon A\rightarrow \St(A)$ с плотным образом --- композиция канонических морфизмов $C\colon A\to A/N$ и $A/N\to \widehat{A/N}$. Покажем, что он обладает нужным универсальным свойством.

Пусть $\varphi\colon A\rightarrow B$ --- непрерывный эрмитов морфизм из алгебры $A$ в $K^*$-алгебру $B$. Тогда для любой $q\in
\mathcal{P}^*_{c}(B)$ композиция $q'=q\circ\varphi$ есть непрерывная $C^*$-преднорма на $A$, так что $\set{q'\colon q\in
\mathcal{P}^*_c(B)}\sbs \mathcal{P}^*_c(A)$.  Поэтому $N\subset
\Ker\varphi=\cap\{N(q')\colon q\in \mathcal{P}^*_{c}(B)\}$, и
существует эрмитов морфизм $\varphi_{1}$ из $A/N$ в $B$ такой, что $\varphi_{1}\circ C=\varphi$. При этом для любой $q\in
\mathcal{P}^*_{c}(B)$ и любого $a+N\in A/N$ имеем
$$
  (q\circ\varphi_{1})(a+N) = (q\circ\varphi_{1})(C(a)) =
  (q\circ\varphi_{1}\circ C)(a) = (q\circ \varphi)(a),
$$
то есть $q\circ\vphi_{1}=(q')^\bullet$. Тем самым, для любой $q\in \mathcal{P}^*_{c}(B)$ преднорма $q\circ\varphi_{1}$ непрерывна.
Значит, и $\varphi$ непрерывен.

Теперь, поскольку $B$ полна, морфизм $\varphi_{1}$ обладает
каноническим продолжением $\overline{\varphi}\colon
\St(A)\rightarrow B$. Из определений $\varphi_{1}$ и
$\overline{\varphi}$ ясно, что $\overline{\varphi}\circ j = \varphi$. Если $\psi$ --- еще один эрмитов морфизм, для которого $\psi\circ j = \varphi$, то $\psi\circ j=\overline{\varphi}\circ j$, и в силу плотности образа $j$ в $\St(A)$ и отделимости $\St(A)$ будем иметь $\psi=\overline{\varphi}$. Предложение доказано.
\end{proof}

\begin{dfn}
Преднорму $p$ на топологической инволютивной алгебре $A$ будем
называть спектральной, если каждый $p$-непрерывный эрмитов характер на $A$ непрерывен относительно топологии алгебры $A$.
\end{dfn}

\begin{prop}
Пусть $A$ --- $LC^*$-алгебра. Для того, чтобы преобразование
Гельфанда $\gel_{A}\colon A\to C_{co}(\Chi^*_{c}(A))$ было
гомеоморфизмом на образ, необходимо и достаточно, чтобы алгебра~$A$ была отделима, а каждая спектральная $C^*$-преднорма на ней~--- непрерывна.
\end{prop}

\begin{proof} "<Необходимость">. Пусть $\gel_{A}$ ---
гомеоморфизм на образ. Тогда, в частности, $\gel_{A}$ инъективно.
Значит, алгебра $A$ строго полупроста и потому отделима (гл.~1,
\S20, предложение ($LC^*$), п.~$4^\circ$). Пусть $p$ ---
спектральная $C^*$-преднорма на $A$. Тогда $\Chi^*_{p}(A)$ ---
компакт в $\Chi^*_{c}(A)$, поэтому $C^*$-преднорма
$p_{\Chi^*_{p}(A)}$ на $C_{co}(\Chi^*_{c}(A))$ непрерывна. Но для
любого $a\in A$ имеем:
\begin{align*}
p(a)&=\sup\{|\chi(a)|\colon \chi\in \Chi^*_{p}(A)\}=\\
    &=\sup\{|\gel_{A}(a)(\chi)|\colon \chi\in \Chi^*_{p}(A)\}=\\
    &=p_{\Chi^*_{p}(A)}(\gel_{A}(a))=(p_{\Chi^*_{p}(A)}\circ \gel_{A})(a).
\end{align*}
Значит, $p$ является композицией $p_{\Chi^*_{p}(A)}$ и $\gel_{A}$ и потому непрерывна.

"<Достаточность">. Пусть $A$ отделима и каждая спектральная
$C^*$-преднорма на ней непрерывна. В силу первого из названных
условий $A$ строго $*$-полупроста. Поэтому $\gel_{A}$ инъективно, то есть является биекцией на образ. В силу второго условия топология алгебры $A$ есть топология равномерной сходимости на
\underline{всех} компактах из пространства $\Chi^*_{c}(A)$. Но и
топология алгебры $C_{co}\Chi^*_{c}(A)$ тоже есть топология
равномерной сходимости на всех компактах из пространства
$\Chi^*_{c}(A)$. Значит, $\gel_{A}$ есть гомеоморфизм на образ, что и требовалось.
\end{proof}

\begin{thm}
Пусть $T$ --- тихоновское пространство. Тогда:

$1^\circ$ Каждая $LC^*$-топология на $C(T)$, согласующаяся с
двойственностью между $C(T)$ и $T$, рассматриваемым как подмножество в $\Chi^*(C(T))$, есть топология равномерной сходимости на некоторой совокупности компактов из $T$, покрывающей~$T$.

$2^\circ$ Среди всех $LC^*$-топологий на $C(T)$, согласующихся с
двойственностью между~$C(T)$ и~$T$, имеются слабейшая и сильнейшая. Слабейшей является топология поточечной сходимости на множестве~$T$, сильнейшей --- компактно-открытая топология. В частности, ком\-пакт\-но-открытая топология на $C(T)$ является сильнейшей среди всех $LC^*$-топологий на $C(T)$, относительно которых непрерывны все фиксированные $($то есть вида $\delta_{T}(t)$\textup{)} характеры и только они.

$3^\circ$ $C^*$-преднорма $p$ на $C_{co}(T)$ непрерывна тогда и
только тогда, когда\linebreak \mbox{$\Chi^*_{p}(C(T))\subset T$},
или, что равносильно, когда \mbox{$\{t\in T\colon
|\delta_{T}(t)|\leq p\}$} --- компактное подмножество в~$T$. В
частности, все спектральные $C^*$-преднормы на $C_{co}(T)$
непрерывны.
\end{thm}

\begin{proof}
Состоит в непосредственном применении теоремы (2.4).
\end{proof}

\begin{prop}
Для любого тихоновского пространства $T$ следующие утверждения
равносильны:

а) $T$ --- $Q$-пространство;

б) каждая $C^*$-преднорма на $C_{co}(T)$ непрерывна.
\end{prop}

\begin{proof} Пусть каждая $C^*$-преднорма на
$C_{co}(T)$ непрерывна. Тогда (предложение (1.4)) и все эрмитовы
характеры на ней непрерывны. Но все непрерывные характеры на алгебре $C_{co}(T)$ фиксированны~\cite[cor. 2.3]{b25}. Значит, имеет место равенство $\Chi^*(C(T))=T$, и $T$ --- $Q$-пространство.

Обратно, пусть $T$ --- $Q$-пространство. Тогда имеет место равенство $\Chi^*(C(T))=T$ и, в силу п.~$2^\circ$ теоремы (2.9),
компактно-открытая топология на алгебре $C(T)$ есть сильнейшая
$LC^*$-топология на ней. Значит, все $C^*$-преднормы на $C_{co}(T)$ непрерывны, и предложение доказано.
\end{proof}

\section[Полунепрерывные снизу $C^*$-преднормы]{Полунепрерывные снизу $C^*$-преднормы на топологических инволютивных алгебрах}

В настоящем параграфе мы рассмотрим вопросы, связанные с
полунепрерывностью снизу $C^*$-преднорм на $LC^*$-алгебрах. Начнём с теоремы, характеризующей полунепрерывные снизу $C^*$-преднормы в терминах пространств эрмитовых характеров.

\begin{thm}
Пусть $A$ --- $LC^*$-алгебра. Для того, чтобы $C^*$-преднорма $p$ на $A$ была полунепрерывна снизу, необходимо и достаточно, чтобы
выполнялось равенство
$$
\Chi^*_{p}(A)=\cl_{\Chi^*(A)}[\Chi^*_{p}(A)\cap \Chi^*_{c}(A)].
$$
\end{thm}

\begin{proof}  "<Необходимость">: методом приведения к
абсурду. Положим, для краткости ${\Chi^*_{p}}'=\cl_{\Chi^*(A)}
[\Chi^*_{p}(A)\cap \Chi^*_{c}(A)]$. И допустим, что, вопреки
утверждению, ${\Chi^*_{p}}' \neq \Chi^*_{p}$. Тогда, поскольку
$\Chi^*_{p}$ замкнуто в $\Chi^*(A)$ и, значит, ${\Chi^*_{p}}'\subset \Chi^*_{p}$, непременно $\Chi^*_{p}\setminus {\Chi^*_{p}}'\ne\varnothing$. Возьмем $\chi\in \Chi^*_{p}\setminus {\Chi^*_{p}}'$ и покажем, что существуют $a\in A$ и $\alpha > 1$
такие, что $\chi(a)=\alpha$ и $0<\gamma(a)<\alpha-1$ при всех
$\gamma\in {\Chi^*_{p}}'$.

Пусть $\gamma\in {\Chi^*_{p}}'$. Тогда $\gamma\neq \chi$ и потому
существует $a'_{\gamma}$ в $A$ такое, что  $\gamma(a'_{\gamma})\neq \chi(a'_{\gamma})$. Полагая
$a_{\gamma}=a'_{\gamma}-\chi(a'_{\gamma})$, получаем, что
$\chi(a_{\gamma})=0$ и $\gamma(a_{\gamma})\neq0$. Умножая, если
нужно, $a_{\gamma}$ на $a^*_{\gamma}$, можно считать, что
$\chi(a_{\gamma})=0$, $\gamma(a_{\gamma})>0$ и
$\delta(a_{\gamma})\geqslant0$ при всех $\delta\in \Chi^*(A)$.
Умножая, если нужно, $a_{\gamma}$ на подходящее число, можно
считать, что $\gamma(a_{\gamma})>1$. Таким образом, для любого
$\gamma\in {\Chi^*_{p}}'$ существует $a_{\gamma}\in A$ такое, что
$\chi(a_{\gamma}) = 0$, $\gamma(a_{\gamma})>1$ и
$\delta(a_{\gamma})\geqslant0$ при всех $\delta\in \Chi^*(A)$. Но
тогда и $\delta(a_{\gamma})>1$ при всех $\delta$ из некоторой
окрестности точки $\gamma$. Поскольку множество ${\Chi^*_{p}}'$
компактно, существует конечное число таких окрестностей точек
$\gamma_{i}$, $i=1,\ \ldots,\ n$, покрывающих в
совокупности~${\Chi^*_{p}}'$. Положим $a' = \sum\limits_{i=1}^n
a_{\gamma_i}$. Тогда $\chi(a')=0$ и $\gamma(a')>1$ при всех
$\gamma\in {\Chi^*_{p}}'$. Положим $\alpha = \max\{\gamma(a')\colon
\gamma\in {\Chi^*_{p}}'\}$ и $a = -a'+\alpha$. Тогда
$\chi(a)=\alpha$ и, при всех $\gamma\in {\Chi^*_{p}}'$, $\gamma(a) = \gamma(\alpha)-\gamma(a') = \alpha-\gamma(a')$. Тем самым $0 <
\gamma(a) < \alpha-1$ при всех $\gamma\in {\Chi^*_{p}}'$, что и
требовалось.

Теперь покажем, что для любого непустого компакта $K$ из
$\Chi^*_{c}(A)$ и любого $\varepsilon>0$ существует $a_{(K,\
\varepsilon)}\in A$, для которого выполняются неравенства
$p_K(a-a_{K,\ \varepsilon})<\varepsilon$ и $p(a_{(K,\
\varepsilon)})<\alpha-1+\varepsilon$.

Пусть $K$ --- непустой компакт из $\Chi^*_{c}(A)$ и $\varepsilon>0$. Без ограничения общности можно предполагать, что $\varepsilon<1$. Положим $H_{p,\ \varepsilon} = \{\gamma\in \Chi^*_{p}\colon \gamma(a)\geq \alpha-1+\varepsilon/ 2\}$. Тогда $H_{p,\ \varepsilon}$~--- замкнутое подмножество компактного множества $\Chi^*_{p}$ и потому само компактно. Кроме того, $H_{p,\ \varepsilon}\cap K=\varnothing$, ибо если бы нашёлся $\gamma\in H_{p,\ \varepsilon}\cap K$, то мы бы имели, с одной стороны, $\gamma(a)\geqslant\alpha-1+\frac{\varepsilon}{2}$ и, с другой стороны, $\gamma(a)<\alpha-1$, ибо $\gamma\in \Chi^*_{p}\cap \Chi^*_{c}(A)\subset {\Chi^*_{p}}'$, что несовместимо. Наконец, $\chi(a)=\alpha>\alpha-1+\varepsilon/2$, так что $\chi\in H_{p,\ \varepsilon}$ и, тем самым, $H_{p,\ \varepsilon}\neq\varnothing$. Таким образом, $H_{p,\ \varepsilon}$ и $K$~--- непустые непересекающиеся компактные подмножества в $\Chi^*(A)$. Значит, существует $\varphi\in C(\Chi^*(A))$ такая, что $0\leqslant\varphi\leqslant1$, $\varphi=1$ на $K$ и $\varphi=0$ на $H_{p,\ \varepsilon}$. Ясно, что при этом
$\gel_{A}(\varphi)=\gel_{A}(a)$ на $K$, $\gel_{A}(a)\varphi=0$ на
$H_{p,\ \varepsilon}$ и $\gel_{A}(a)\varphi<\alpha-1+\varepsilon/2$,
ибо, по определению $H_{p,\ \varepsilon}$,
$\gel_{A}(a)<\alpha-1+\varepsilon/2$ на $\Chi^*_{p}\backslash H_{p,\ \varepsilon}$. Теперь, на основании теоремы Вейерштрасса-Стоуна, алгебра $\gel_{A}[A]$ плотна в алгебре $C_{co}(\Chi^*(A))$ как самосопряжённая подалгебра, содержащая постоянные функции и различающая точки множества $\Chi^*(A)$. Поэтому существует $a_{K,\ \varepsilon}\in A$ такой, что
$$
  |\gel_{A}(a_{(K,\ \varepsilon)})-\gel_{A}(a)\varphi|<\frac{\varepsilon}{2}
$$
на $K\cup\Chi^*_{p}$. Следовательно,
$$
  |\gel_{A}(a_{(K,\ \varepsilon)})-\gel_{A}(a)|<\frac{\varepsilon}{2}<\varepsilon
$$
на $K$ и
$$
|\gel_{A}(a_{(K,\varepsilon)})|<\alpha-1+\varepsilon
$$
на $\Chi^*_{p}$. Отсюда и получаем требуемое.

Теперь на множестве всех упорядоченных пар $(K,\ \varepsilon)$ с
компактным $K\subset\Chi^*_{c}(A)$ и положительным $\varepsilon$
введём отношение порядка, полагая, по определению,
$$
 (K,\ \varepsilon)\leqslant(K',\ \varepsilon')\leftrightarrow K\subset
 K'\; \&\; \varepsilon'\leqslant\varepsilon.
$$
На основании только что доказанного получим, что
$$
  \gel_{A}(a_{(K,\ \varepsilon)})\rightarrow \gel_{A}(a)
$$
равномерно на \underline{всех} компактах из $\Chi^*_{c}(A)$.
А~поскольку $A$ --- $LC^*$-алгебра, её топология есть топология
равномерной сходимости на \underline{некоторой} совокупности
компактов из $\Chi^*_{c}(A)$ (предложение (2.4)). Значит, $a_{(K,\ \varepsilon)}\rightarrow a$ относительно топологии алгебры $A$. При этом
$$
  \varliminf p(a_{(K,\ \varepsilon)}) =
  \varliminf p_{\Chi^*_{p}}(\gel_{A}(a_{(K,\ \varepsilon)})) =
  \varliminf (\alpha-1+\varepsilon)=\alpha-1.
$$
Но $\alpha=\chi(a)$ для $\chi\in \Chi^*_{p}$. Значит,
$$
  \alpha\leqslant p(a)=\sup\{|\gamma(a)|\colon \gamma\in \Chi^*_{p}\}.
$$
Получается, что $(a_{(K,\ \varepsilon)})\rightarrow a$, но
$\varliminf p(a_{(K,\ \varepsilon)}) \neq p(a)$. Однако это
несовместимо с полунепрерывностью снизу преднормы $p$ на
алгебре~$A$. Получаем противоречие, которое и доказывает
"<необходимость">.

"<Достаточность">. Пусть
$$
  \Chi^*_{p}=\cl_{\Chi^*(A)}[\Chi^*_{p}\cap \Chi^*_{c}(A)].
$$
Тогда для любого $a\in A$ имеем
$$
  p(a)=\sup\{|\chi(a)|\colon \chi\in \Chi^*_{p}\cap\Chi^*_{c}(A)\}.
$$
В силу этого $p$ оказывается верхней огибающей семейства непрерывных функций
$$
\chi\mapsto |\chi(a)|,\ \chi\in \Chi^*_{p}\cap\Chi^*_{c}(A).
$$
Значит, она полунепрерывна снизу, как и утверждалось. Доказательство теоремы завершено.
\end{proof}

Хорошо известно, что каждая полунепрерывная снизу вещественнозначная функция на тихоновском пространстве является верхней огибающей некоторого семейства непрерывных функций. Доказанная только что теорема даёт для $C^*$-преднорм на $LC^*$-алгебрах больше.

\begin{cor}
Каждая полунепрерывная снизу $C^*$-преднорма на $LC^*$-алгебре
является верхней огибающей некоторого семейства минимальных
непрерывных $C^*$-преднорм.
\end{cor}

\begin{proof}
Это тотчас же следует из приведённого выше доказательства
"<достаточности">\ для теоремы (4.1) и предложения (1.4).
\end{proof}

Теперь вспомним (теорема (2.4), п.~$5^\circ$), что минимальные
непрерывные $C^*$-пред\-нор\-мы --- одни и те же для всех
$LC^*$-топологий, согласующихся с данной двойственностью. Тем самым получаем
\begin{cor}
Полунепрерывные снизу $C^*$-преднормы --- одни и те же для всех
$LC^*$-топологий, согласующихся с данной двойственностью.
\end{cor}

Для алгебры $C_{co}(T)$ теорему (3.1) можно немного дополнить.

\begin{prop}
Пусть $T$ --- тихоновское пространство, $p$~--- $C^*$-преднорма на $C(T)$. Тогда следующие утверждения равносильны.

а) Преднорма $p$ полунепрерывна снизу на $C_{co}(T)$.

б) Ядро преднормы $p$ замкнуто в $C_{co}(T)$.

в) Преднорма $p$ есть верхняя огибающая некоторого семейства
минимальных непрерывных $C^*$-преднормы.

г) $\Chi^*_{p}(C(T))=\cl_{\upsilon T}\{t\in T\colon
|\delta_{T}(t)|\leqslant p\}$.
\end{prop}

\begin{proof} Равносильность утверждений а) и г) следует из теоремы~(3.1) и следствия (9.1). Что а) $\Rightarrow$ б) и в) $\Rightarrow$ a) --- ясно. Остается показать, что б) $\Rightarrow$ в).

Пусть ядро $N(p)$ преднормы $p$ на $C(T)$ замкнуто в $C_{co}(T)$.
Тогда $N(p)$~--- замкнутый самосопряжённый идеал в $C_{co}(T)$.
Поэтому~\cite[теорема 2.1]{b25} существует замкнутое подмножество
$Z$ пространства $T$ такое, что для любой функции $\varphi\in C(T)$ соотношение $p(\varphi)=0$ равносильно соотношению
$\varphi[Z]=\{0\}$.

Теперь заметим, что $Z\subset\Chi^*_{p}$. Действительно, в противном случае нашлись бы $z\in Z\backslash \Chi^*_{p}$ и $\varphi\in C(T) = C(\upsilon T)$ такие, что $\varphi[\Chi^*_{p}]=\{0\}$ и $\varphi(z)=1$. Тогда бы получилось, что $p(\varphi)=0$, а $\varphi[Z]\neq\{0\}$, но это невозможно.

Кроме того,
\begin{align*}
\varphi[Z] =\{0\} &\Leftrightarrow p(\varphi) = 0\\
                  &\Leftrightarrow \varphi\in \Ker\pi_{p}\\
                  &\Leftrightarrow \pi_{p}(\varphi)=0\\
                  &\Leftrightarrow (\forall \chi\in\Chi^*(C(\Chi)_{p}))
                  (\chi(\pi_{p}(\varphi))=0)\\
                  &\Leftrightarrow (\forall\chi\in
                  \Chi^*_{p})(\chi(\varphi)=0),\\
\intertext{то есть}
 \varphi[Z] = \{0\} &\Leftrightarrow \varphi[\Chi^*_{p}]=\{0\}.
\end{align*}
Значит, $Z$ плотно в $\Chi^*_{p}$. Поэтому
$$
  p(\varphi)=\sup\{|\chi(\varphi)|\colon \chi\in
  \Chi^*_{p}\}=\sup\{|\chi(\varphi)|\colon \chi\in Z\}.
$$
Но все $\chi\in Z$~--- непрерывные характеры, так что все функции
$\varphi\mapsto|\chi(\varphi)|$~--- непрерывные $C^*$-преднормы,
причём~--- минимальные (предложение(1.4)). Тем самым $p$~--- верхняя огибающая некоторого семейства минимальных непрерывных
$C^*$-преднорм, что и требовалось.
\end{proof}

\section{$C^*$-бочечные $LC^*$-алгебры}

Напомним, что локально выпуклое пространство бочечно тогда и только тогда, когда каждая полунепрерывная снизу преднорма на нем
непрерывна. По аналогии с этим мы вводим следующее

\begin{dfn}
$LC^*$-алгебру будем называть \textit{$C^*$-бочечной}, если каждая полунепрерывная снизу $C^*$-преднорма на ней непрерывна.
\end{dfn}

Напомним также, что подмножество $S$ пространства $T$ называется
ограниченным, если каждая функция из $C(T)$ ограниченна на $S$.
Тихоновское пространство~$T$ называется $\mu$-пространством, если в нём замыкание каждого ограниченного подмножества компактно
(см.~\cite{b7}, а также~\cite{b1}, где такие пространства называются $NS$-про\-странствами, то есть пространствами Нахбина-Сироты,~--- по имени математиков, впервые обративших внимание на эти пространства). Известная теорема Нахбина и Сироты (\cite{b34},~\cite{b42}; см. также~\cite{b1}, теорема 2.5-1) гласит, что если $T$~--- тихоновское пространство, то пространство $C_{co}(T)$ бочечно тогда и только тогда, когда $T$ является $\mu$-пространством. Наш следующий результат является обобщением этой теоремы на произвольные $LC^*$-алгебры. Чтобы сформулировать его, напомним, что подмножество $S$ топологического пространства $T$ называется $k$-замкнутым, если пересечение~$S$ с каждым компактом $K$ из $T$ замкнуто в $K$.

\begin{thm}
$LC^*$-алгебра $A$ $C^*$-бочечна тогда и только тогда, когда
выполняются следующие два условия:

а) $\Chi^*_{c}(A)$ $k$-замкнуто в $\Chi^*(A)$;

б) топология алгебры $A$ есть топология равномерной сходимости на
всех компактах из $\Chi^*_{c}(A)$.
\end{thm}

\begin{proof} "<Только тогда">. Пусть $A$ --- $C^*$-бочечная
$LC^*$-алгебра. Покажем, что выполняется условие а). Пусть $K$~--- компакт из $\Chi^*(A)$ и $S=K\cap\Chi^*_{c}(A)$. Тогда $S$ ---
$A$-ограниченное подмножество в пространстве $\Chi^*(A)$. Поэтому
$p^S$ (см. определение (1.5)) является (предложение (1.6))
$C^*$-преднормой на $A$ и притом --- полунепрерывной снизу, ибо
$$
  p^S=\sup\{p^{\{t\}}\colon t\in S\},
$$
а каждая преднорма $p^{\{t\}}$ с $t\in \Chi^*_{c}$ --- непрерывна. В силу теоремы (3.1)
$$
  \Chi^*_{p^S}(A)=\cl_{\Chi^*(A)}[\Chi^*_{p^S}(A)\cap\Chi^*_{c}(A)].
$$
А в силу $C^*$-бочечности алгебры $A$ преднорма $p^S$ непрерывна и потому $\cl_{\Chi^*(A)}S=\Chi^*_{p^S}(A)\subset\Chi^*_{c}(A)$
(предложение (1.6), б) и предложение (2.1)). С другой стороны,
$S\subset K$, а $K$ замкнуто в $\Chi^*(A)$, поэтому
$\cl_{\Chi^*(A)}S\subset K$. Таким образом, $\cl_{\Chi^*(A)}S\subset K\cap \Chi^*_{c}(A) = S$, так что $S$ замкнуто в $\Chi^*(A)$ и, тем более, в $K$, что и требовалось.

Теперь покажем, что выполняется условие б). Для любого компакта $K$
из $\Chi^*_{c}(A)$ $C^*$-преднорма $p^K$ на $A$ полунепрерывна снизу (как верхняя огибающая семейства непрерывных $C^*$-преднорм
$p^{\{t\}}$, $t\in K$) и потому непрерывна (в силу $C^*$-бочечности алгебры $A$). Тем самым, топология алгебры $A$ сильнее топологии равномерной сходимости на всех компактах из $\Chi^*_{c}(A)$. Но последняя~--- сильнейшая среди всех $LC^*$-топологий на $A$, согласующихся с двойственностью между $A$ и $\Chi^*_{c}(A)$. Значит, топология алгебры $A$ совпадает с топологией равномерной сходимости на всех компактах из $\Chi^*_{c}$, и условие б) выполняется.

Обратно, пусть выполняются условия а) и б). Покажем, что алгебра~$A$ $C^*$-бочечна. Пусть $p$ --- полунепрерывная снизу $C^*$-преднорма на $A$. Тогда $\Chi^*_{p}(A)$ --- компактное подмножество в $\Chi^*(A)$, причём, в силу теоремы (3.1),
$$
  \Chi^*_{p}(A)=\cl_{\Chi^*(A)}[\Chi^*_{p}(A)\cap \Chi^*_{c}(A)].
$$
Но $\Chi^*_{c}(A)$ $k$-замкнуто в $\Chi^*(A)$, поэтому
$$
  \cl_{\Chi^*(A)}[\Chi^*_{p}(A)\cap\Chi^*_{c}(A)] = \Chi^*_{p}(A)\cap\Chi^*_{c}(A).
$$
Тем самым
$$
  \Chi^*_{p}(A)=\Chi^*_{p}(A)\cap\Chi^*_{c}(A),
$$
то есть $\Chi^*_{p}(A)\subset \Chi^*_{c}(A)$. Таким образом,
$\Chi^*_{p}(A)$~--- компакт из $\Chi^*_{p}(A)$ и, в силу~б),
преднорма~$p$ непрерывна. Это и требовалось.
\end{proof}

Чтобы увидеть, что теорема (4.2) является обобщением теоремы Нахбина и Сироты, нам понадобится следующая теорема.

\begin{thm}
Пусть $T$~--- тихоновское пространство. Тогда следующие утверждения равносильны:

а) $T$ --- $\mu$-пространство;

б) $T$ $k$-замкнуто в $\upsilon T$;

в) алгебра $C_{co}(T)$ $C^*$-бочечна;

г) пространство $C_{co}(T)$ бочечно.
\end{thm}

\begin{proof}
Равносильность а) и г) --- это упомянутая выше теорема Нахбина и Сироты. Равносильность б) и в) получается из теоремы (4.2)
отождествлением пространства~$T$ с $\Chi^*_{c}(C_{co}(T))$, а
пространства $\upsilon T$ --- с $\Chi^*(C(T))$ (теорема (1.9)).
Установим равносильность а) и б).

а) $\Rightarrow$ б) Пусть $T$ --- $\mu$-пространство. Покажем, что
$T$ $k$-замкнуто в $\upsilon T$. Пусть $K$ --- компакт в $\upsilon
T$. Тогда $K\cap T$ --- замкнутое ограниченное подмножество
пространства~$T$. Поскольку $T$ --- $\mu$-пространство, $K\cap T$
компактно. Значит, оно замкнуто в $K$. Тем самым, $T$ $k$-замкнуто в
$\upsilon T$, что и требовалось.

б) $\Rightarrow$ а) Пусть $T$ $k$-замкнуто в $\upsilon T$. Докажем,
что $T$ --- $\mu$-пространство. Пусть $Z$~--- замкнутое ограниченное
подмножество в $T$. Тогда $\cl_{\upsilon T}Z$ --- компакт в
$\upsilon T$~\cite[8.E.1]{b9}. В силу $k$-замкнутости $T$ в
$\upsilon T$ пересечение $T\cap \cl_{\upsilon T}Z$ замкнуто в
$\cl_{\upsilon T}Z$ и потому --- компактно. Но $T\cap \cl_{\upsilon
T} Z = \cl_T Z = Z$ (в силу замкнутости $Z$). Таким образом, $Z$
компактно, что и требовалось.
\end{proof}

Теперь мы видим, что теорема (4.2) действительно является обобщением
теоремы Нахбина и Сироты о бочечных пространствах непрерывных
функций. Однако даже для случая $A = C_{co}(T)$ теорема (4.2) даёт
больше, чем теорема Нахбина и Сироты.

Действительно, в силу отмеченной в теореме (4.3) равносильности
бочечности и $C^*$-бочечности для алгебр непрерывных функций,
теорема Нахбина и Сироты "<в одну сторону">\ состоит в том, что если
алгебра непрерывных функций с компактно-открытой топологией
$C^*$-бочечна, то $T$ является $\mu$-пространством. Наша же теорема
в этом случае гласит, что если алгебра непрерывных функций с
какой-то $LC^*$-топологией $C^*$-бочечна, то $T$ является
$\mu$-пространством и топология алгебры является компактно-открытой
(так что одно из условий не только не используется, но даже выводится из другого).

В заключение параграфа рассмотрим вопрос об ассоциированных
$C^*$-бочечных $LC^*$-алгебрах. Комура доказал~\cite{b17}, что для
любого локально выпуклого пространства $(E,\ \tau)$ на пространстве
$E$ существует бочечная локально выпуклая топология~$\tau_t$,
которая, во-первых, сильнее топологии $\tau$ и, во-вторых, слабее
любой бочечной локально выпуклой топологии на $E$, более сильной,
чем $\tau$. Это равносильно тому, что тождественное отображение
$(E,\ \tau_t)\to (E,\ \tau)$ непрерывно и каждое непрерывное
линейное отображение $(F,\ \sigma)\to (E,\ \tau)$ из бочечного
пространства $(F,\ \sigma)$ в пространство $(E,\ \tau)$ непрерывно и
как отображение из $(F,\ \sigma)$ в $(E,\ \tau_t)$. Таким образом,
категория бочечных локально выпуклых пространств является
корефлексивной подкатегорией категории всех локально выпуклых
пространств. Следующая теорема показывает, что категория
$C^*$-бочечных $LC^*$-алгебр играет аналогичную роль в категории
всех $LC^*$-алгебр.

\begin{thm}
Категория $C^*$-бочечных $LC^*$-алгебр является корефлексивной
подкатегорией категории всех $LC^*$-алгебр. В частности, на каждой
$LC^*$-алгебре существует ровно одна $C^*$-бочечная
$LC^*$-то\-по\-ло\-гия, более сильная, чем исходная, и более слабая, чем любая $C^*$-бочечная $LC^*$-топология, более сильная, чем исходная.
\end{thm}

\begin{proof} Пусть $A$ --- инволютивная алгебра,
$\tau$ --- $LC^*$-топология на $A$ и $\Gamma$ --- $k$-замыкание
множества $\Chi^*_{c}(A,\ \tau)$ в $\Chi^*(A)$. Обозначим через
$\tau^*_{t}$ топологию равномерной сходимости на всех компактах из
$\Chi^*(A)$, содержащихся в $\Gamma$. Тогда $\Chi^*_{c}(A,\
\tau^*_{t})=\Gamma$, так что $\Chi^*_{c}(A,\ \tau^*_{t})$
$k$-замкнуто в $\Chi^*(A)$ и топология алгебры $(A,\ \tau^*_{t})$
есть топология равномерной сходимости на всех компактах пространства
$\Chi^*_{c}(A,\ \tau^*_{t})$. В силу теоремы (4.2) это означает, что
алгебра $(A,\ \tau^*_{t})$ $C^*$-бочечна. Кроме того, ясно, что
топология $\tau^*$ сильнее исходной.

Пусть $B$ --- $C^*$-бочечная $LC^*$-алгебра, $u\colon B \rightarrow
(A,\ \tau)$ --- непрерывный эрмитов морфизм. Покажем, что $u$
непрерывен и как отображение из $B$ в $(A,\ \tau^*_{t})$, или, что
то же самое, для любой непрерывной $C^*$-преднормы $p$ на $(A,\
\tau^*_{t})$ преднорма $p\circ u$ на $B$ непрерывна.

Пусть $p$ --- непрерывная $C^*$-преднорма на $(A,\ \tau^*_{t})$. Для
любого $b\in B$ имеем:
$$
  (p\circ u)(b)=p(u(b))=\sup\{|\chi(u(b))|\colon  \chi\in \Chi^*_{p}(A)\}=
$$
$$
  =\sup\{|(\chi\circ u)(b)|\colon  \chi\in \Chi^*_{p}(A)\}=
$$
$$
  =\sup\{|\chi'(b)|\colon \chi\in \Chi^*(u)[\Chi^*_{p}(A)]\}.
$$
Иначе говоря,
$$
  p\circ u=p^{\Chi^*(u)[\Chi^*_{p}(A)]}.
$$

Далее, в силу непрерывности $p$,
$$
  \Chi^*_p(A)\subset \Chi^*_c(A, \tau^*_t)=\Gamma,
$$
откуда
$$
  \Chi^*(u)[\Chi^*_p(A)]\subset \Chi^*(u)[\Gamma].
$$

C другой стороны, в силу непрерывности $u$
$$
  \Chi^*(u)[\Chi^*_c(A,\ \tau)]\subset\Chi^*_c(B).
$$
Но $\Chi^*(u)$ есть непрерывное отображение из $\Chi^*(A)$ в
$\Chi^*(B)$ и, значит, из $k\Chi^*(A)$ в $k\Chi^*(B)$, поэтому
$$
  \Chi^*(u)[\cl_{k\Chi^*(A)}\Chi^*_c(A,\tau)] \subset
  \cl_{k\Chi^*(B)}\Chi^*_c(B).
$$
При этом, по определению $\Gamma$,
$$
  \cl_{k\Chi^*(A)}\Chi^*_c(A,\tau)=\Gamma,
$$
а в силу бочечности $B$
$$
  \cl_{k\Chi^*(B)}\Chi^*_c(B)=\Chi^*_c(B).
$$
Значит,
$$
  \Chi^*(u)[\Gamma]\subset\Chi^*_c(B)
$$
и, тем самым,
$$
  \Chi^*(u)[\Chi^*_p(A)]\subset\Chi^*_c(B).
$$

Теперь, $\Chi^*_p(A)$ --- компакт, а $\Chi^*(u)$, как уже
говорилось, непрерывно. Значит, $\Chi^*(u)[\Chi^*_p(A)]$ --- компакт
в $\Chi^*_c(B)$. Однако, $B$ --- $C^*$-бочечна, и её топологией
является топология равномерной сходимости на всех компактах из
$\Chi^*_c(B)$. Поэтому преднорма $p\circ u$ на $B$ непрерывна. Это и
означает, что морфизм $u$ непрерывен и как отображение из $B$ в
$(A,\ \tau^*_t)$. Тем самым доказано, что категория $C^*$-бочечных
$LC^*$-алгебр является корефлексивной подкатегорией категории всех
$LC^*$-алгебр. Остальные утверждения теоремы теперь очевидны.
\end{proof}

\begin{dfn}
Топологию $\tau^*_t$, построенную при доказательстве теоремы (4.4),
будем называть {\em $C^*$-бочечным усилением} топологии $\tau$
алгебры $A$. Алгебру  $A$ с топологией $\tau^*_t$ будем обозначать
через $b^*A$ и называть ассоциированной $C^*$-бочечной
$LC^*$-алгеброй $LC^*$-алгебры.
\end{dfn}

Следующие три результата выводятся из теоремы (4.4) стандартными
теоретико-категорными рассуждениями.

\begin{prop}
Пусть $A_1$, $A_2$ --- $LC^*$-алгебры, $u$-непрерывный эрмитов
морфизм из $A_1$ в $A_2$. Тогда $u$ является непрерывным морфизмом и
из $b^*A_1$ в $b^*A_2$.
\end{prop}

\begin{prop}
$LC^*$-алгебра $A$ $C^*$-бочечна тогда и только тогда, когда
\mbox{$b^*A = A$}.
\end{prop}

\begin{prop}
Индуктивный предел $C^*$-бочечных $LC^*$-алгебр является
\mbox{$C^*$-бо}\-чечной $LC^*$-алгеброй.
\end{prop}

И последнее в этом параграфе предложение. Бухвальтер и
Шмет~\cite{b8} доказали, что для любого тихоновского пространства
$T$ бочечные усиления топологии поточечной сходимости и топологии
компактной сходимости на $C(T)$ совпадают. Из теоремы (4.4) и её
доказательства непосредственно следует намного более общее

\begin{prop}
Пусть $A$ --- инволютивная алгебра, $\Gamma$ --- непустое
подмножество в $\Chi^*(A)$. Тогда все $LC^*$-топологии на $A$,
согласующиеся с двойственностью между $A$ и $\Gamma$, имеют одно и
то же $C^*$-бочечное усиление.
\end{prop}

\chapter[Теоремы двойственности для некоторых классов пространств]{Теоремы двойственности для компактологических
         и некоторых классов топологических пространств}

\section{Топологически порожденные компактологические
пространства}

В этом параграфе напоминаются основные сведения о компактологических пространствах, опровергается одно из утверждений, содержащихся в книге Дж.~Б.~Купера~\cite{b19}, и даётся характеризация топологически порожденных компактологических пространств.

Компактологическим пространством называется~\cite{b6} множество $S$ с классом $\varkappa$ его подмножеств, называемых компактами, для которых выполняются следующие условия:
\begin{itemize}
\item[а)] $\varkappa$ --- фильтрующееся вправо (по отношению $\subset$) покрытие множества $S$;
\item[б)] каждый компакт $K\in \varkappa$ так наделён отделимой компактной топологией, что если $K\in \varkappa$, $L\in \varkappa$ и $K\subset L$, то вложение $K\rightarrow L$ непрерывно;
\item[в)] если $K\subset L$ и $L\in \varkappa$, то замыкание в $L$ множества $K$ тоже принадлежит $\varkappa$.
\end{itemize}

В силу б) если $K\in \varkappa$, $L\in \varkappa$ и $K\subset L$, то $K$ --- подпространство в $L$, а в силу с) все компактные
подпространства любого компакта $K\in \varkappa$ тоже принадлежат
$\varkappa$.

Пусть $(S,\ \varkappa)$, $(S_1,\ \varkappa_1)$ ---
компактологические пространства. Отображение $k\colon S\rightarrow S_1$ называется компактологическим морфизмом из $(S,\ \varkappa)$ в $(S_1,\ \varkappa_1)$, если для любого $K$ из $\varkappa$ найдется такой $K_1$ из $\varkappa_1$, что $k[K]\subset K_1$ и сужение $k|_K\colon  K\rightarrow K_1$ непрерывно; или, что равносильно, если для любого $K\in \varkappa$ множество $k[K]$, наделённое топологией образа, принадлежит $\varkappa_1$.

На каждом топологическом пространстве $(T,\ \tau)$ существует
естественная компактология, $\varkappa_{\tau}$, состоящая, по
определению, из всех подмножеств пространства $T$, компактных
относительно топологии $\tau$. В частности, такая компактология
существует на множестве $\mathbb{C}$ комплексных чисел.
Компактологические пространства, получаемые таким образом из
подходящих топологических пространств, называются топологически
порожденными. Ясно, что морфизмы из компактологического пространства $(S,\ \varkappa)$ в компактологическое пространство, порождённое топологическим пространством $T$,~--- это просто функции из $S$ в $T$, непрерывные на всех компактах из $\varkappa$.

Компактологическое пространство $(S,\ \varkappa)$ называется
регулярным, если морфизмы из $(S,\ \varkappa)$ в $\mathbb{C}$
различают точки множества $S$. В работе А.~Бухвальтера~\cite{b6}
доказано (теорема I.2.2), что компактологическое пространство $(S,\ \varkappa)$ регулярно тогда и только тогда, когда на множестве $S$ существует отделимая вполне регулярная топология, сужение которой на любой компакт $(K,\ \tau_K)\in \varkappa$ совпадает с $\tau_K$. В книге Дж.~Б.~Купера~\cite{b19} утверждается больше (теорема A.2.2): компактологическое пространство $(S,\ \varkappa)$ регулярно тогда и
только тогда, когда на множестве $S$ существует отделимая вполне
регулярная топология, порождающая компактологию $\varkappa$. Это
более сильное утверждение --- неверно.
\begin{thm}
Существуют регулярные компактологии, не порождаемые никакими
топологиями.
\end{thm}
\begin{proof}  Пусть $\mathcal{C}$~--- множество всех сходящихся
последовательностей действительных чисел. Для любой $f\in
\mathcal{C}$ положим $R(f) = \{f_n|n\in \mathbb{N}\}\cup\{\lim f\}$. Каждое $R(f)$ наделим топологией, индуцированной из $\mathbb{R}$, превращая его тем самым в отделимое компактное топологическое пространство. Обозначим через $\varkappa$ совокупность всех замкнутых подмножеств любых конечных объединений множеств $R(f)$ и наделим каждое такое подмножество топологией, индуцированной из $\R$. Получим компактологию на множестве $\R$, причём состоящую из множеств, не более чем счётных. Покажем, что эта компактология не порождается никакой топологией на $\R$.

Пусть $t$ --- какая-нибудь топология на $\R$, которая на каждом
множестве $R(f)$ совпадает с обычной евклидовой топологией $e$ на
$\R$. Тогда $t$ слабее $e$. Действительно, пусть $Z$~---
$t$-замкнутое подмножество в $\R$. Покажем, что оно  и $e$-замкнуто. Пусть $z_0$ --- $e$-предельная точка множества $Z$. Тогда существует последовательность $(z_n)$ элементов множества $Z$, сходящаяся к $z_0$ относительно топологии $e$. Поскольку топология~$t$ индуцирует на $R(z)$ ту же топологию, что и $e$, последовательность $z$ сходится к~$z_0$ и относительно топологии $t$. Но множество $Z$ $t$-замкнуто. Значит, $z_n\in Z$. Тем самым множество $Z$ оказывается также и $e$-замкнутым. Таким образом, каждое $t$-замкнутое множество является также и $e$-замкнутым, а это и означает, что топология $t$ действительно слабее топологии $e$.

В силу доказанного каждое $e$-компактное множество является также и $t$-компакт\-ным. В частности, $t$-компактным оказывается отрезок $[0,\ 1]$. Однако он несчётен, а компактология $\varkappa$ состоит лишь из счётных множеств. Поэтому отрезок $[0,\ 1]$ не принадлежит $\varkappa$. Значит, $\varkappa$ не порождается топологией $t$. В силу произвольности выбора $t$ компактология $\varkappa$ не порождается никакой топологией, что и требовалось доказать.
\end{proof}

В связи с доказанной только что теоремой естественно возникает
вопрос, какие компактологии топологически порождены. Ответ: те и
только те, которые устойчивы относительно образования компактных
индуктивных пределов в категории отделимых топологических
пространств. Оставшаяся часть параграфа посвящена уточнению и
обоснованию этого ответа.

\begin{thm}
Пусть $(S,\ \varkappa)$ --- регулярное компактологическое
пространство,~$t_\varkappa$ --- топология на $S$, относительно
которой замкнуты те и только те подмножества множества $S$,
пересечение которых с каждым $K\in \varkappa$ замкнуто в~$K$,
и~$\widehat{\varkappa}$~--- совокупность всех компактных
подпространств пространства $(S,\ t_\varkappa)$. Тогда:

$1^\circ$ $\varkappa\subset\widehat{\varkappa};$

$2^\circ$ $t_\varkappa$ --- отделимая $k$-топология;

$3^\circ$ пространство $(S,\ t_\varkappa)$ является индуктивным
пределом пространств $K\in \varkappa$ в категории отделимых
топологических пространств;

$4^\circ$ компакт $(C,\ t)$ принадлежит компактологии
$\widehat{\varkappa}$ тогда и только тогда, когда он содержится в
$S$ и является индуктивным пределом пространств $K\in \varkappa$ с $K\subset C$ в категории отделимых топологических пространств.
\end{thm}

\begin{proof}  1. Пусть $K\in \varkappa$. Покажем, что
$(K\in \widehat{\varkappa})$, то есть является компактным
подпространством пространства $(S,\ t_\varkappa)$. Пусть $Z$ ---
подмножество в $K,$ замкнутое относительно топологии, которую
индуцирует в $K$ топология $t_\varkappa$. Тогда существует
$t_\varkappa$-замкнутое подмножество $Z'\subset S$ такое, что
$Z'\cap K=Z$. Но $t_\varkappa$-замкнутость $Z'$ означает, что его
пересечение с любым компактом из $\varkappa$ замкнуто в этом
компакте. Стало быть $Z'\cap K=Z$ замкнуто в $K$, и мы показали, что каждое $t_\varkappa$-замкнутое подмножество в $K$ замкнуто также и относительно топологии $\tau_K$ компакта $K$. Тем самым топология $\tau_K$ сильнее топологии в $K$, индуцированной из $(S,\ t_\varkappa)$. Но $\tau_K$ компактна, а $t_\varkappa$ отделима (в силу регулярности $\varkappa$). Поэтому $\tau_K$ совпадает с топологией, индуцированной из $(S,\ t_\varkappa)$, и $(K,\ \tau_K)$~--- компактное подпространство пространства $(S,\
t_\varkappa).$

2. Отделимость $t_{\varkappa}$ уже была отмечена выше. Теперь
покажем, что $t_{\varkappa}$ ---$k$-топология, то есть каждое
подмножество $Z$ множества $S$ замкнуто относительно топологии
$t_{\varkappa}$ тогда и только тогда, когда $Z\cap K$ замкнуто в $K$ для любого компактного подпространства $K$ пространства $(S,\
t_{\varkappa})$. Утверждение "<тогда">\ непосредственно следует из п.~$1^\circ$. Установим "<только тогда">. Пусть $Z$~---
подпространство в $S$, замкнутое относительно топологии
$t_{\varkappa}$, а $K$ --- компактно относительно той же топологии. Тогда и $Z$, и $K$ (в силу отделимости $t_{\varkappa}$)~--- замкнутые подмножества в $(S,\ t_{\varkappa})$. Поэтому их пересечение замкнуто в $(S,\ t_{\varkappa})$ и, тем более, в $K$.

3. Пусть $i^K\colon  K\rightarrow S$, $i^K_{K'}\colon  K\rightarrow K'$ $(K,\ K'\in \varkappa,\ K\subset K')$ --- тождественные вложения. Покажем, что топологическое пространство $(S,\ t_\varkappa)$ является индуктивным пределом спектра топологических пространств $\{K,\ i^K_{K'}\colon  K\subset K',\ K,\ K'\in \varkappa\}$ относительно отображений $i^K$. Рассмотрим произвольное семейство $(f^K\colon K\rightarrow S_1)_{K\in \varkappa}$ непрерывных отображений в категории отделимых топологических пространств, согласующихся с отображениями $i^K_{K'}$. Определим отображение $f$ из $S$ в $S_1$ условием: $f(t)=f^K(t)$ для любого $K\in \varkappa$ с $t\in K$. Из того, что $\varkappa$ покрывает $S$, а отображения $f^K$ согласуются с отображениями $i^K_{K'}$, обычным образом выводится, что отображение $f$ определено корректно и однозначно. Из определения топологии $t_\varkappa$ непосредственно следует, что $f$ непрерывно. Наконец, из определения $f$ ясно, что
$f\circ i^K = f^K$ при всех $K\in \varkappa$.

4. Пусть $(C,\ t)$ --- компакт, являющийся индуктивным пределом
компактов $(K,\ \tau_K)$ из $\varkappa$ с $K\subset C$ в категории
отделимых топологических пространств. Покажем, что $(C,\ t)$ ---
подпространство пространства $(S,\ t_\varkappa)$.

Прежде всего, поскольку все вложения $j^K\colon  K\rightarrow C$
$\tau_K$ --- $t$-непрерывны, топология~$\tau_K$ есть сужение на $K$
топологии $t$ при всех $K\in \varkappa$ с $K\subset C$.
Следовательно, все компакты $(K,\ \tau_K)$ с $K\subset C$ ---
подпространства в $(C,\ t)$ и потому --- замкнутые в $(C,\ t)$
подмножества.

Далее, подмножество $Z$ множества $C$ замкнуто относительно
топологии $t$ тогда и только тогда, когда $Z\cap K$ замкнуто в $K$
относительно топологии $\tau_K=t|_K$ для всех $(K,\ \tau_K)\in
\varkappa$ с $K\subset C$. Действительно, совокупность всех таких
подмножеств в~$C$ определяет на $C$ топологию, и эта топология
обладает, как легко видеть, требуемым универсальным свойством. Тем
самым она является топологией индуктивного предела и потому
совпадает с $t$.

Теперь покажем, что топология $t$ совпадает с топологией, которую
индуцирует на $C$ топология $t_\varkappa$. Для этого, прежде всего,
заметим, что $t$ сильнее топологии, индуцированной из пространства
$(S,\ t_\varkappa)$. Действительно, пусть $Z$ --- подмножество
в~$C$, замкнутое относительно топологии, индуцированной из $(S,\
t_\varkappa)$. Тогда существует замкнутое подмножество $Z'$ в $(S,\
t_\varkappa)$ такое, что $Z=Z'\cap C$. При этом замкнутость $Z'$ в
$(S,\ t_\varkappa)$ означает, что $Z'\cap K$ замкнуто в  $(K,\
\tau_K)$ для любого $(K,\ \tau_K)\in \varkappa,$ так что, в
частности, это выполняется и для $(K,\ \tau_K)$ с $K\subset C$. Но
если $K\subset C$, то $Z'\cap K=Z'\cap C\cap K=Z\cap K$, так что и
$Z\cap K$ замкнуто в $(K,\ \tau_K)$ для любого $(K,\ \tau_K)\in
\varkappa$ с $K\cap C$. Значит, $Z$ замкнуто в $(C,\ t),$ и мы
показали, что топология $t$ сильнее той, что индуцирована из $(S,\
t_\varkappa)$. Но эта топология отделима, а $t$ компактна. Значит,
$t$ совпадает с топологией, индуцированной из пространства $(S,\
t_\varkappa)$, так что $(C,\ t)$ --- подпространство в пространстве
$(S,\ t_\varkappa)$.

Обратно, пусть $(C,\ t)$ --- подпространство в $(S,\ t_\varkappa)$.
Тогда $C$ --- замкнутое подмножество в $(S,\ t_\varkappa)$, так что
пересечение $C\cap K$ замкнуто в $(K,\ \tau_K)$ для любого $(K,\
\tau_K)\in \varkappa$. Покажем, что $(C,\ t)$  является индуктивным
пределом компактов $(K,\ \tau_K)\in \varkappa$ с $K\subset C$. Пусть
$j^K\colon  K\rightarrow C$, $i^K_{K'}\colon  K\rightarrow K'$ ---
тождественные вложения $(K,\ K'\in \varkappa,\ K\subset K')$, а
$(f^K\colon K\rightarrow S_1)_{K\in\varkappa,\ K\subset C}$ ---
семейство отображений из $K$ в $S_1$ в категории отделимых
топологических пространств, (непрерывных и) согласующихся с
отображениями $i^K_{K'}$. Определим отображение $f\colon
C\rightarrow (S,\ t_1)$, полагая $f(t)=f^K(t)$ для любого $K\in
\varkappa$ с $t\in K$ и $K\subset C$. Снова обычными образом
убеждаемся, что $f$ определено корректно и однозначно и что $f\circ
j^K=f^K$ при всех $(K,\ \tau_\varkappa)\in \varkappa$ с $K\subset
C$. Пусть $Z$ --- замкнутое подмножество а пространстве $(S_1,\
t_1)$. Убедимся, что $f^{-1}[Z]$ замкнуто в $(C,\ t)$. Как уже
отмечалось, для любого $(K,\ \tau_K)\in \varkappa$ пересечение
$C\cap K$ замкнуто в $(K,\ \tau_K)$ (даже --- в $(S,\ t_\varkappa)$)
и, наделённое топологией, индуцированной из $(K,\ \tau_\varkappa)$,
принадлежит $\varkappa$. Теперь, для любого $(K,\ \tau_\varkappa)\in
\varkappa$,
$$
  f^{-1}[Z]\cap K=f^{-1}[Z]\cap C\cap K=(f^{C\cap K})^{-1}[Z]\cap(C\cap K).
$$
Причём $(f^{C\cap K})^{-1}[Z]$ замкнуто в $C\cap K$, а $C\cap K$
замкнуто в $(K,\ \tau_\varkappa)$. Тем самым $f^{-1}[Z]\cap K$
замкнуто в $(K,\ \tau_\varkappa)$ для любого $(K,\
\tau_\varkappa)\in \varkappa$ и, тем самым, $f^{-1}[Z]$ замкнуто в
$(S,\ t_\varkappa)$. Следовательно, оно замкнуто и в $(C,\ t$, и
непрерывность $f$ установлена. Этим доказательство теоремы (4.2)
завершено.
\end{proof}

\begin{cor}
Если компактология $\varkappa$ на множестве $S$ порождается хотя бы
одной топологией, то она порождается и топологией, определенной в
предыдущей теореме.
\end{cor}

\begin{proof}  Действительно, если компактология
$\varkappa$ на множестве $S$ порождается топологией $t$, то
определённая при доказательстве теоремы (4.2) топология
$t_\varkappa$ есть $k$-усиление топологии $t$, а хорошо известно,
что компактологии, порождаемые какой-либо топологией и её
$k$-усилением, совпадают.
\end{proof}

Топологию $t_\varkappa$ естественно назвать $k$-топологией,
порождаемой компактологией~$\varkappa$. В этой терминологии
следствие (4.3) означает, что компактология порождается хотя бы
одной топологией тогда и только тогда, когда она порождается
порождённой ею $k$-топологией.

\begin{dfn}
Будем говорить, что компактология $\varkappa$ на множестве $S$
устойчива относительно образования индуктивных пределов в категории
отделимых топологических пространств, если каждый компакт $(C,\ t)$
с $C\subset S$, являющийся индуктивным пределом некоторого семейства
компактов из $\varkappa$ в категории отделимых топологических
пространств, принадлежит $\varkappa$.
\end{dfn}

Теперь из теоремы 2 и её следствия немедленно вытекает

\begin{thm}
Для того, чтобы регулярная компактология порождалась хотя бы одной
топологией, необходимо и достаточно, чтобы она была устойчива
относительно образования компактных индуктивных пределов в категории
отделимых топологических пространств.
\end{thm}

Нам не известно, существует ли аналогичная характеризация для
компактологий, порождаемых вполне регулярными топологиями. То, что
нам удалось получить в этом направлении, собрано в следующем
предложении. Доказательства аналогичны приведённым выше для
соответствующих случаев.

\begin{prop}
$1^\circ$ Если компактология $\varkappa$ на множестве $S$
порождается хотя бы одной вполне регулярной топологией, то она
порождается и слабой топологией, определяемой множеством всех
компактологических морфизмов из $(S,\ \varkappa)$ в $\mathbb{C}$.
Множество $S$, наделённое этой топологией, является индуктивным
пределом компактов $(K,\ \tau_K)\in\varkappa$ в категории
тихоновских пространств.

$2^\circ$ Если регулярная компактология на множестве $S$ порождается
хотя бы одной вполне регулярной топологией, то она устойчива
относительно образования компактных индуктивных пределов в категории
тихоновских пространств.
\end{prop}

\section{Двойственность для регулярных компактологических
пространств}

В этом параграфе доказывается принадлежащая А.~Бухвальтеру~\cite{b6}
теорема о двойственности между некоторой категорией топологических
инволютивных алгебр и категорией регулярных компактологических
пространств. Это делается по нескольким причинам. Во-первых, в
доказательствах наших последующих результатов используется не только
сама теорема двойственности, но и некоторые части приводимого нами
её доказательства. Во-вторых, строго говоря, в упомянутой работе
А.~Бухвальтера сформулирована и доказана не теорема двойственности,
а лишь две теоремы об изоморфизме, в которых полностью отсутствует
всё, что касается функторных свойств, весьма существенных для наших
последующих результатов. В-третьих, одна их этих двух теорем
сформулирована А.~Бухвальтером довольно необычно: часть условий в
ней заменена многоточием (по этой причине мы и написали в начале
параграфа "<{\em некоторой}"> категорией топологических инволютивных
алгебр вместо того, чтобы указать эту категорию явно). Какие условия
ещё имел в виду А.~Бухвальтер, ставя многоточие, трудно сказать.
Приводимое нами доказательство показывает, что никаких больше
условий не нужно, и это --- ещё одна из причин, побудивших нас
привести здесь своё доказательство теоремы А.~Бухвальтера. Наконец,
заметим, что наше доказательство хотя и близко по духу к
доказательству самого Бухвальтера, всё же отличается от него и
кажется нам более прямым и прозрачным. К тому же оно позволяет
распространить на произвольные коммутативные $K^*$-алгебры известную
теорему, согласно которой из полноты пространства $C_{co}(T)$
следует, что $T$ является $k_{\mathbb{C}}$-пространством.

\begin{thm}[А.~Бухвальтер]
Категория, дуальная к категории регулярных компактологических
пространств, эквивалентна категории $K^*$-алгебр.
\end{thm}

\begin{proof} Пусть $\mathcal{K}$ --- категория
регулярных компактологических пространств, $\mathcal{K}^*$ ---
категория $K^*$-алгебр.

а) Определим контравариантные функторы \mbox{$\mathcal{C}\colon
\mathcal{K}\rightarrow\mathcal{K}^*$} и\linebreak \mbox{$\mathcal{X}^*_c\colon
\mathcal{K}^*\rightarrow\mathcal{K}$}.

a1) Пусть $(S,\ \varkappa)$ --- регулярное компактологическое
пространство. Обозначим через $\mathcal{C}(S)$ инволютивную алгебру
всех компактологических морфизмов из $(S,\ \varkappa)$ в
$\mathbb{C}$, наделённую топологией равномерной сходимости на всех
компактах из $\varkappa$. Ясно, что $\mathcal{C}(S)$ отделима и
полна, так что является $K^*$-алгеброй. Для любого
компактологического морфизма $u\colon (S,\ \varkappa)\rightarrow
(S_1,\ \varkappa_1)$ определим отображение $\mathcal{C}(u)\colon
\mathcal{C}(S_1)\rightarrow\mathcal{C}(S)$ условием:
$\mathcal{C}(u)(\varphi_1)=\varphi_1\circ u$ для всех $\varphi_1\in
\mathcal{C}(S_1)$. Без труда проверяется, что $\mathcal{C}(u)$
определено корректно (то есть для всех $\varphi_1\in
\mathcal{C}(S_1)$ $\mathcal{C}(u)(\varphi_1)\in \mathcal{C}(S))$ и
является эрмитовым морфизмом. Для доказательства непрерывности
$\mathcal{C}(u)$ возьмём произвольный $K\in \varkappa$ и покажем,
что преднорма $p_K\circ\mathcal{C}(u)$ на алгебре $\mathcal{C}(S_1)$
непрерывна. Пусть $\varphi_1\in\mathcal{C}(S_1)$. Тогда
\begin{multline*}
  [p_K\circ\mathcal{C}(u)](\varphi_1) = p_K(\mathcal{C}(u)(\varphi_1))
  =\\=
  p_K(\varphi_1\circ u) = \sup\{|\varphi_1(u(x))|\colon x\in K\}
  =\\=
  \sup\{|\varphi_1(x_1)|\colon x_1\in u[K]\} = p_{u[K]}(\varphi_1).
\end{multline*}

Таким образом, $p_K\circ\mathcal{C}(u) = p_{u[K]}$. Но $u[K]\in
\varkappa_1$. Значит, $p_K\circ\mathcal{C}(u)$ непрерывна, что и
требовалось. Наконец, стандартные вычисления показывают, что для
любых компактологических морфизмов $u\colon (S,\
\varkappa)\rightarrow(S_1,\ \varkappa_1)$ и $u_1\colon (S_1,\
\varkappa_1)\rightarrow(S_2,\ \varkappa_2)$ имеет место равенство
$\mathcal{C}(u_1\circ u) = \mathcal{C}(u)\circ\mathcal{C}(u_1)$. Тем самым $\mathcal{C}$ --- контравариантный функтор из $\mathcal{K}$ в $\mathcal{K^*}$.

a2) Пусть $A$ --- $K^*$-алгебра, $\varkappa_A = \{\Chi^*_p (A)\colon p\in \mathcal{P}^*_c (A)\}$. Тогда, во-первых, для любой $p\in\mathcal{P}^*_c(A)$ множество $\Chi^*_p (A)$ содержится в
$\Chi^*_c(A)$ и является компактом в индуцированной топологии
(глава~2, теоремы (2.1) и (1.1)). Во-вторых, для любого
$\chi\in\Chi^*_c (A)$ $C^*$-преднорма $p\colon a\mapsto|\chi(a)|$ на $A$ непрерывна и $\chi\in\Chi^*_p(A)$, так что $\varkappa_A$~--- покрытие $\Chi^*_c(A)$. В-третьих, если $K_1,\ K_2\in \varkappa_A$, то существуют $p_1$, $p_2\in\mathcal{P}^*_c(A)$ такие, что $K_1 =
\Chi^*_{p_1}(A)$ и $K_2 = \Chi^*_{p_2}(A)$; тогда $p\colon =
\max\{p_1,\ p_2\}\in\mathcal{P}^*_c(A)$ и $K_1\cup
K_2\subset\Chi^*_p(A)\in\varkappa_A$. В-четвёртых, если $K,\
L\in\varkappa$ и $K\subset L$, то и $K$, и $L$ --- подпространства в $\Chi^*_c(A)$, так что вложение $K\rightarrow L$ непрерывно.
В-пятых, если $H\subset K\in\varkappa,$ то
$$
  \cl_KH = \cl_{\Chi^*(A)}H=\Chi^*_{p^H}(A)\in\varkappa_A.
$$
И, наконец, если $\chi_1,\ \chi_2\in\Chi^*_c(A)$ и
$\chi_1\neq\chi_2$, то существует $a\in A$, для которого
$\chi_1(a)\neq\chi_2(a)$, или, что --- то же самое,
$\gel_{A}(a)(\chi_1)\neq\gel_A(a)(\chi_2)$. При этом $\gel_A(a)$~--- непрерывное отображение из $\Chi^*_c$ в $\mathbb{C}$ и, тем более, непрерывно на всех компактах $\Chi^*_p(A)$ с $p\in \mathcal{P}^*_c(A)$. Всё это показывает, что $\varkappa_A$~--- регулярная компактология на множестве $\Chi^*_c(A)$.

Для любой $K^*$-алгебры $A$ обозначим через $\mathcal{X}^*_c(A)$
компактологическое пространство $(\Chi^*_c(A),\ \varkappa_A)$. Для любого морфизма $K^*$-алгебр $u\colon A\rightarrow A_1$ определим
$\mathcal{X}^*_c(u)$ как $\Chi^*_c(u)$. Поскольку $\Chi^*_c(u)$~--- непрерывное отображение из пространства $\Chi^*_c(A_1)$ в пространство $\Chi^*_c(A)$, оно, тем более, является компактологическим морфизмом из $\mathcal{X}^*_c(A_1)$ в
$\mathcal{X}^*_c(A)$. Функторные свойства отображений $\Chi^*_c(u)$ уже отмечались (глава 1, \S 19). Таким образом,
$\mathcal{X}^*_c$~--- контравариантный функтор из $\mathcal{K}^*$ в $\mathcal{K}$.

б) Теперь покажем, что композиции $\mathcal{X}^*_c\circ\mathcal{C}$
и $\mathcal{C}\circ\mathcal{X}^*_c$ изоморфны тождественным
функторам $I_{\mathcal{K}}$ и $I_{\mathcal{K}^*}$ соответственно.

б1) Пусть $(S,\ \varkappa)$ --- регулярное компактологическое
пространство. Наделим $S$ слабой топологией $r_\varkappa$,
определяемой всеми функциями из $\mathcal{C}(S)$. Обозначим через
$C_\varkappa(S)$ алгебру всех непрерывных комплекснозначных функций на топологическом пространстве $(S,\ r_\varkappa)$, наделённую топологией равномерной сходимости на членах семейства $\varkappa$. Ясно, что тогда $\mathcal{C}(S) = C_\varkappa(S)$. При этом преобразование Дирака
$$
  \delta_S\colon  S\rightarrow \Chi^*_c C_\varkappa(S)
$$
является гомеоморфизмом на образ (ибо $(S,\ r_\varkappa)$)~--- тихоновское пространство) и сюръективно~\cite[следствие 2.3]{b25}. Кроме того, для любого $K\in\varkappa$ множество $\delta_S[K] =
\Chi^*_{p_K}C_\varkappa(S)$ принадлежит $\varkappa_{C_\varkappa(S)}$
(то есть $\varkappa_A$ с $A = C_\varkappa(S)$), причём иных
компактов в $\varkappa_{C_\varkappa(s)}$ нет. Тем самым, $\delta_S$ оказывается изоморфизмом компактологических пространств из пространства $(S,\ \varkappa)$ в пространство
$\mathcal{X}^*_c(\mathcal{C}(S,\ \varkappa))$. То, что для любого
морфизма $u\colon (S,\ \varkappa)\rightarrow (S_1,\ \varkappa_1)$
выполняется равенство $\delta_{S'_1}\circ u =
\mathcal{X}^*_c\mathcal{C}(u)\circ\delta_S$~--- очевидно (ибо
сводится к аналогичному утверждению относительно тихоновских
пространств, для которых это уже отмечалось (глава 1, \S22)).

б2) Пусть $A$~--- $K^*$-алгебра. Поскольку все члены
$\varkappa_A$~--- компактные подмножества в $\Chi^*_c(A)$, алгебра $C\Chi^*_c(A)$ является подалгеброй алгебры
$\mathcal{C}\mathcal{X}^*_c(A)$ и преобразование Гельфанда можно
рассматривать как морфизм из $A$ в $\mathcal{C}\mathcal{X}^*_c(A)$.
Поскольку $A$ полупроста (глава 1, предложение $(LC^*)$,
п.~$5^\circ$) $\gel_A$ инъективно (там же, п.~$4^\circ$). При этом $\gel_A[A]$ --- самосопряжённая, отделяющая точки и содержащая
константы подалгебра в $\mathcal{C}\mathcal{X}^*_c(A)$. По теореме Вейерштрасса-Стоуна $\gel_A[A]$ плотна в
$\mathcal{C}\mathcal{X}^*_c(A)$. Но и топология алгебры $A$, и
топология алгебры $\mathcal{C}\mathcal{X}^*_c(A)$ есть топология
равномерной сходимости на всех множествах $\Chi^*_p(A)$ c $p\in
\mathcal{P}^*_c(A)$: относительно $\mathcal{C}\mathcal{X}^*_c(A)$
это следует из определений $\mathcal{C}$ и $\mathcal{X}^*_c$, а
относительно алгебры $A$~--- из предложения (1.1) главы 2.
Поэтому~$\gel_A$~--- гомеоморфизм на образ. В силу полноты
алгебры~$A$ подалгебра~$\gel_A[A]$ замкнута в
$\mathcal{C}\mathcal{X}^*_c(A)$. Значит, $\gel_A[A] =
\mathcal{C}\mathcal{X}^*_c(A)$ и $\gel_A$ --- изоморфизм. Функторные свойства $\gel$ уже отмечались (глава 1). Теорема доказана.
\end{proof}

\begin{cor}
Пусть $A$ --- отделимая $LC^*$-алгебра. Тогда:

$1^\circ$ Если $A$ полна, то каждая функция
$\Chi^*_c(A)\rightarrow\mathcal{C}$, непрерывная на каждом
$K\in\varkappa_A$, непрерывна на $\Chi^*_c(A)$. В частности, тогда $\Chi^*_c(A)$ есть $k_\mathbb{C}$-пространство.

$2^\circ$ Пополнение $\widehat{A}$ алгебры $A$ реализуется алгеброй $\mathcal{C}(\mathcal{X}^*_c(A))$ всех функций
$\Chi^*_c(A)\rightarrow\mathbb{C},$ непрерывных на каждом
$K\in\varkappa_A$.

$3^\circ$ Топология пространства $\Chi^*_c(\widehat{A})$ совпадает со слабой топологией, определяемой всеми функциями
$\Chi^*_c(A)\rightarrow\mathbb{C},$ непрерывными на каждом $K\in
\varkappa_A$

$4^\circ$ Для того, чтобы отображение сужения $r\colon
\Chi^*_c(\widehat{A}\rightarrow\Chi^*_c(A))$ было гомеоморфизмом,
необходимо и достаточно, чтобы на пространстве $\Chi^*_c(A)$ были
непрерывны все функции со значениями в $\mathbb{C}$, непрерывные на каждом $K\in\varkappa_A$. В частности, для этого необходимо (a если каждый компакт из пространства $\Chi^*_c(A)$ равностепенно
непрерывен, то и достаточно), чтобы $\Chi^*_c(A)$ было
$k_\Cnum$-пространством.
\end{cor}

\begin{proof} В последней части доказательства теоремы~(2.1) мы установили, что
$$
  \gel_A[A]=\mathcal{C}\mathcal{X}^*_c(A).
$$
Кроме того, ясно, что
$$
  \gel_A[A]\subset C\Chi^*_c(A)
\quad\text{и}\quad
  C\Chi^*_c(A)\subset\mathcal{C}\mathcal{X}^*_c(A).
$$
Таким образом, для любой $K^*$-алгебры $A$
$$
  C\Chi^*_c(A)=\mathcal{C}\mathcal{X}^*_c(A).
$$
Это и утверждается в п.~$1^\circ$.

2. Пополнение $\widehat{A}$ $LC^*$-алгебра $A$ есть $K^*$-алгебра. В силу теоремы (2.1) $\widehat{A}\cong\mathcal{C}\mathcal{X}^*_c(\widehat{A})$. Это и
есть п.~$2^\circ$.

Пп.~$3^\circ$, $4^\circ$ теперь очевидны.
\end{proof}

{\sc Замечания.} 1. П.~$1^\circ$ следствия (2.2) есть обобщение
половины теоремы 1 из~\cite{b49}.

2. П.~$2^\circ$ следствия (2.2) является аналогом гротендиковского описания пополнения отделимого локально выпуклого пространства (см., например,~\cite{b37}).

3. Вопрос о том, является ли отображение $r$ гомеоморфизмом, был
поставлен А.~Маллиосом в~\cite{b24}. В. Дитрих в~\cite{b11} построил пример, дающий отрицательный ответ на этот вопрос. Благодаря следствию (2.2) всё становится очевидным: если~$T$~--- тихоновское пространство, не являющееся $k_\Cnum$-пространством, то $\Chi^*_c(C(T))$ гомеоморфно~$T$ и, значит, не является
$k_\Cnum$-пространством; стало быть, в этом случае отображение
сужения из $\Chi^*_c(\widehat{C(T)})$ на $\Chi^*_c(C(T))$ не
является гомеоморфизмом.


\section[Теоремы двойственности для некоторых классов $K^*$-алгебр]{Теоремы двойственности для некоторых классов $K^*$-алгебр и $k_{\mathbb{C}}$-пространств}

\begin{lemma}
Пусть $(S,\ \varkappa)$ --- регулярное компактологическое
пространство. Для того, чтобы компактология $\varkappa$ порождалась
хотя бы одной вполне регулярной топологией, необходимо и достаточно,
чтобы на $K^*$-алгебре $\mathcal{C}(S)$ каждая спектральная
$C^*$-преднорма была непрерывна.
\end{lemma}

\begin{proof} "<Необходимость">. Пусть $\varkappa$
порождается хотя бы одной вполне регулярной топологией. Тогда
(предложение (1.6)) она порождается и слабой топологией
$r_\varkappa$, определяемой всеми функциями из $\mathcal{C}(S)$, то
есть топологией пространства $\Chi^*_c(\mathcal{C}(S))$. Если теперь
$p$ --- спектральная $C^*$-преднорма на алгебре $\mathcal{C}(S)$, то
$\Chi^*_p\mathcal{C}(S)\subset\Chi^*_c\mathcal{C}(S)$ и, по
теореме~(2.1) $\Chi^*_c\mathcal{C}(S) = S$, так что
$\Chi^*_p\mathcal{C}(S)$~--- компактное подмножество пространства
$(S,\ r_\varkappa)$. Но все компакты из $(S, r_\varkappa)$
принадлежат $\varkappa$ (ибо~$\varkappa$ порождена $r_\varkappa$), а
все преднормы $p_K$ с $K\in \varkappa$ на алгебре $\mathcal{C}(S)$
непрерывны. Значит,~$p$ непрерывна (ибо $p = p_K$ с $K =
\Chi^*_p\mathcal{C}(S)$) (глава 1, предложение (1.1), б), что и
требовалось.

"<Достаточность">. Пусть каждая спектральная $C^*$-преднорма на
алгебре $\mathcal{C}(S)$ непрерывна. Тогда, в частности, для каждого
компакта $K$ из пространства $(S,\ r_\varkappa)$ $C^*$-преднорма
$p_K$ на $\mathcal{C}(S)$ непрерывна. Но на $\mathcal{C}(S)$
непрерывны только те $C^*$-пред\-нормы, которые отвечают компактам
из $\varkappa$. Значит, каждое компактное подпространство
пространства $(S,\ r_\varkappa)$ принадлежит $\varkappa$. Тем самым,
компактология $\varkappa$ порождена топологией $r_\varkappa$, что и
требовалось.
\end{proof}

\begin{thm}
Категория, дуальная к категории вполне регулярных
$k_\mathbb{C}$-про\-странств, эквивалентна категории всех
$K^*$-алгебр, на которых каждая спектральная $C^*$-преднорма
непрерывна.
\end{thm}

\begin{proof} В силу теоремы (2.1) и леммы (3.1)
категория $K^*$-алгебр, на которых каждая спектральная
$C^*$-преднорма непрерывна, эквивалентна категории
компактологических пространств, порождаемых вполне регулярными
отделимыми топологическими пространствами. Но категория таких
компактологических пространств эквивалентна (даже --- изоморфна)
категории вполне регулярных $k_\mathbb{C}$-пространств. Отсюда и
следует теорема.
\end{proof}

Согласно~\cite{b38} комплексная топологическая алгебра $A$ с
совместно непрерывным умножением называется $B$-полной, если каждый
непрерывный почти открытый морфизм из $A$ в такую же топологическую
алгебру $A_1$ открыт. (Напомним, что линейное отображение $f$
локально выпуклого пространства $E$ в локально выпуклое пространство
$F$ называется почти открытым, если для любой окрестности нуля $U$ в
$E$ замыкание множества $f[U]$ в $F$ является окрестностью нуля в
$F$.) В~\cite{b38} доказано, что алгебра $C_{co}(T)$ $B$-полна тогда
и только тогда, когда пространство $T$ является $k$-пространством.
Отсюда и из теоремы $(3.2)$ следует

\begin{thm}
Категория, дуальная к категории вполне регулярных $k$-прост\-ранств,
эквивалентна категории $B$-полных $K^*$-алгебр, на которых каждая
спектральная $C^*$-преднорма непрерывна.
\end{thm}

Теперь вспомним, что, согласно теореме (4.2) из главы 2,
$LC^*$-алгебра $A$ $C^*$-бо\-че\-ч\-на тогда и только тогда, когда
множество $\Chi^*_c(A)$ $k$-замкнуто в $\Chi^*(A)$ и топология
алгебры $A$ есть топология равномерной сходимости на всех компактах
из $\Chi^*_c(A)$. В силу этого на $C^*$-бочечных $LC^*$-алгебрах все
спектральные $C^*$-преднормы непрерывны. Кроме того (теорема (4.3)
главы 1), алгебра $C_{co}(T)$ $C^*$-бочечна тогда и только тогда,
когда $T$ является $\mu$-пространством. Назовем топологическое
пространство $\mu k_{\mathbb{C}}$-пространством, если оно является
$\mu$-пространством и $k_\mathbb{C}$-пространством одновременно.
Coчетая всё сказанное с теоремой (3.2), приходим к следующему
результату.

\begin{thm}
Категория, дуальная к категории $C^*$-бочечных $K^*$-алгебр,
эквивалентна категории $\mu k_{\mathbb{C}}$-пространств.
\end{thm}

Согласно теореме Нахбина и Сироты о борнологических пространствах
непрерывных функций (\cite{b34},~\cite{b42}; см. также~\cite[теорема 2.6-1]{b1},) пространство $C_{co}(T)$ борнологично тогда и только
тогда, когда пространство $T$ вещественнокомпактно. Поэтому:

\begin{thm}
Категория, дуальная к категории борнологических $K^*$-алгебр,
эквивалентна категории вещественнокомпактных
$k_{\mathbb{C}}$-пространств.
\end{thm}

Соединяя (3.4) с (3.2), а (3.5) с (3.3), получаем ещё две теоремы.

\begin{thm}
Категория, дуальная к категории $B$-полных $C^*$-бочечных
\mbox{$K^*$-ал}\-гебр, эквивалентна категории $\mu k$-пространств.
\end{thm}

\begin{thm}
Категория, дуальная к категории $B$-полных борнологических
$K^*$-алгебр, эквивалентна категории вещественнокомпактных
$k$-пространств.
\end{thm}

Сочетая результаты этого параграфа с известными результатами о
$C_{co}(T)$, можно получить ещё целый ряд теорем двойственности и их
следствий. Ограничимся двумя примерами.

Известно, что пространство $C_{co}(T)$ метризуемо тогда и только
тогда, когда $T$ хемикомпактно (то есть представляется в виде
объединения счётного семейства компактов $K_i$ таких, что каждое
компактное подмножество пространства $T$ содержится в одном из
$K_i$) и является $k$-пространством. Отсюда вытекает

\begin{thm}
Категория, дуальная к категории метризуемых $K^*$-алгебр,
эквивалентна категории вполне регулярных хемикомпактных
$k$-пространств.
\end{thm}

Но полные метризуемые локально выпуклые пространства
бочечны~\cite[глава IV, терема 2]{b37}. Получаем

\begin{cor}
Каждое вполне регулярное хемикомпактное $k$-пространство является
$\mu$-пространством.
\end{cor}


\section{$lc^*$-алгебры и двойственность для локально компактных
пространств}

\begin{dfn}
$K^*$-алгебру $A$ будем называть $lc^*$-алгеброй, если она
удовлетворяет условию:
\begin{itemize}
\item[($lc^*$)] для любого $\chi\in\Chi^*_c(A)$ найдутся $a\in A$ и
непрерывная $C^*$-преднорма $p$ на $A$ такие, что $\chi(a)=1$ и
$a\cdot N(p)=\{0\}$.
\end{itemize}
\end{dfn}

\begin{predo}
Обращаем внимание на различие между понятиями $LC^*$-алгебры и
$lc^*$-алгебры.
\end{predo}

Полную подкатегорию категории $\mathcal{K}^*$ $K^*$-алгебр,
образованную $lc^*$-алгебрами, будем называть категорией
$lc^*$-алгебр и обозначать через $\mathcal{L}c^*$.

\begin{thm}
Категория, дуальная к категории отделимых локально компактных
пространств, эквивалентная категории $lc^*$-алгебр.
\end{thm}

\begin{proof} Временно назовем регулярное компактологическое
пространство $(S,\ \varkappa)$ строго регулярным, если для любого
$t\in S$ найдутся $K\in \varkappa$ и $\varphi\in\mathcal{C}(S)$
такие, что $\varphi(t) = 1$ и $\varphi = 0$ вне $K$. Покажем, что
категория $\widetilde{\mathcal{L}c}$ строго регулярных
компактологических пространств эквивалентна категории $\mathcal{L}c$
отделимых локально компактных пространств.

Определим функторы $k\colon
\mathcal{L}c\rightarrow\widetilde{\mathcal{L}c}$ и $t\colon
\widetilde{\mathcal{L}c}\rightarrow\mathcal{L}c$, полагая $k(T,\
\tau) = (T,\ \varkappa_\tau)$ и $t(S,\ \varkappa) = (S,\
r_\varkappa)$, где $\varkappa_\tau$ --- класс всех компактов
топологического пространства $(T,\ \tau)$, а $r_\varkappa$ ---
слабая топология на $T$, определяемая всеми функциями
$\varphi\in\mathcal{C}(S)$. На морфизмах $k$ и $t$ действуют как
тождественные отображения. Конечно, нужно проверить, что $k$ и $t$
действительно является функторами между указанными категориями. Что
$k$ есть функтор из $\mathcal{L}c$ в $\widetilde{\mathcal{L}c}$, ---
очевидно. Кроме того, ясно, что если $(S,\ \varkappa)$ --- строго
регулярное компактологическое пространство, то топология
$r_\varkappa$ локально компактна. Остается показать, что любой
компактологический морфизм $u\colon (S,\ \varkappa)\rightarrow
(S_1,\ \varkappa_1)$ между строго регулярными пространствами
является непрерывным отображением из $(S,\ r_\varkappa)$ в $(S_1,\
r_{\varkappa_1})$. Но поскольку $r_\varkappa$ и $r_{\varkappa_1}$
--- слабые топологии, определяемые семействами $\mathcal{C}(S)$ и
$\mathcal{C}(S_1)$, для проверки непрерывности $u$ достаточно
установить, что для любой функции $\varphi_1\in\mathcal{C}(S_1)$
композиция $\varphi_1\circ u$ принадлежит $\mathcal{C}(S)$, а это
очевидно.

Теперь покажем, что $t\circ k = I_{\mathcal{L}c}$ и $k\circ t =
I_{\widetilde{\mathcal{L}c}}$. Ясно, что для любого $(T,\ \tau)\in
\mathcal{L}c$ имеем $C(T) = \mathcal{C}(T,\ \varkappa_\tau)$. Далее,
в силу полной регулярности локально компактных пространств, $\tau$~--- слабейшая из топологий, относительно которых непрерывны все $\varphi\in C(T)$. Поэтому $r_{\varkappa_\tau}=\tau$ и $(t\circ
k)(T,\ \tau) = t(T,\ \varkappa_\tau) = (T,\ r_{\varkappa_\tau}) =
(T,\ \tau)$. Это и означает, что $t\circ k = I_{\mathcal{L}c}$.

Пусть $(S,\ \varkappa)\in\widetilde{\mathcal{L}c}$. Покажем, что
$\varkappa = \varkappa_{r_\varkappa}$, то есть компакты в $(S,\
r_\varkappa)$ суть члены $\varkappa$ и только они.

а) Пусть $K\in\varkappa$. Для любой $\varphi\in\mathcal{C}(S)$
композиция $\varphi\circ I^K$ функции $\varphi$ с тождественным
вложением $I^K\colon K\rightarrow S$ совпадает с сужением $\varphi$
на $K$ и потому непрерывна. Значит, $i^K$ является непрерывным
отображением из $(K,\ \tau)$ в $(S,\ r_\varkappa)$. Но $(K,\
\tau_K)$ компактно, а $(S,\ r_\varkappa)$, в силу регулярности $(S,\
\varkappa)$, отделимо. Поэтому $i^K$ --- гомеоморфизм на образ, и
$(K,\ \tau_K)$~--- компактное подпространство в $(S,\ r_\varkappa)$.

б) Пусть $(K,\ \tau_K)$ --- компактное подпространство в $(S,\
r_\varkappa)$. Тогда для любого $x\in K$ найдется $K_x\in \varkappa$
такой, что $K_x$ является $r_\varkappa$-окрестностью точки $x$
(поскольку компактология $\varkappa$ строго регулярна). Внутренности
таких $K_x$ в совокупности покрывают компакт $K$, поэтому существует
конечное число $x_1,\ \ldots,\ x_n$ точек $x$, для которых
$$
  K\subset\bigcup\limits_{i=1}^n K_{x_i}.
$$

Поскольку $\varkappa$ фильтруется вправо по отношению $\subset$,
существует $L\in\varkappa$, для которого
$$
  \bigcup\limits_{i=1}^n K_{x_i}\subset L.
$$
При этом и $(K,\ \tau_K)$, и $(L,\ \tau_L)$ --- подпространства в
$(S,\ r_\varkappa)$, значит, $(K,\ \tau_K)$ --- подпространство в
$(L,\ \tau_L)$ и $K = \cl_LK\in\varkappa$. Этим доказано, что
$\varkappa_{r_\varkappa} = \varkappa$. Теперь $(k\circ t)(S,\
\varkappa)=(S,\ \varkappa_{\tau_\varkappa}) = (S,\ \varkappa)$. Тем
самым, $k\circ t = I_{\widetilde{Lc}}$, что и требовалось.

Итак, категория отделимых локальных локально компактных пространств
эквивалентна категории строго регулярных компактологических
пространств. Значит, теорема будет доказана, если мы установим, что
категория, дуальная к категории $\widetilde{\mathcal{L}c}$ строго
регулярных компактологических пространств, эквивалентна категории
$\Lcal c^*$ $lc^*$-алгебр. Для этого воспользуемся теоремой (2.1).
Согласно доказательству этой теоремы контравариантные функторы
$\mathcal{C}$ и $\mathcal{X}^*_c$ осуществляют двойственность между
категориями $\mathcal{K}$ регулярных компактологических пространств
и $\mathbb{}\mathcal{K}^*$ $K^*$-алгебр. Нам остается только
показать, что строго регулярным пространствам отвечают при этом
$lc^*$-алгебры, а $lc^*$-алгебрам --- строго регулярные
пространства.

Пусть $(S,\ \varkappa)\in \ob\widetilde{\Lcal c}$. Ясно, что
$\mathcal{C}(S,\ \varkappa)$ --- $K^*$-алгебра. Проверим выполнение
условия $(lc^*)$. Пусть $\chi$ --- произвольный элемент множества
$\Chi^*\mathcal{C}(S,\ \varkappa)$. В силу теоремы (2.1) существует
$t\in S$ такой, что $\chi = \delta_S(t)$. А в силу строгой
регулярности пространства $(S,\ \varkappa)$ существует $K\in
\varkappa$ и такая $\varphi\in\mathcal{C}(S,\ \varkappa)$, что
$\varphi(t)=1$ и $\varphi = 0$ вне $K$. Поскольку каждая $\psi\in
N(p_K)$ равна 0 на $K$, отсюда получаем, что $\chi(\varphi) =
\delta_S(t)(\varphi) = \varphi(t) = 1$ и $\varphi\cdot N(p_K) =
\{0\}$. Это и есть условие $(lc^*)$ для $\mathcal{C}(S,\
\varkappa)$.

Пусть $A$ $lc^*$-алгебра и $\chi\in\Chi^*_c(A)$. Тогда существуют
$a\in A$ и $p\in \Pcal^*_c(A)$ такие, что $\chi(a)=1$ и $a\cdot N(p)
= \{0\}$. Покажем, что $\chi\in\Chi^*_p(A)$. Допустим противное.
Тогда найдётся $\varphi\in\mathcal{C}(\mathcal{X}^*_c(A))$ такая,
что $\varphi(\chi) = 1$ и $ \varphi[\Chi^*_p(A)] = \{0\}$.

Поскольку $\gel_A\colon
A\rightarrow\mathcal{C}\mathcal{X}^*_c(A)$~--- изоморфизм
(теорема~(2.1)), существует $b\in A$ такое, что $\varphi=\gel_A(b)$.
При этом
$$
  p(b)=\sup\limits_{\chi\in X^*_p(A)}|\chi(b)| = \sup\limits_{\chi\in
  X^*_p(A)}|\gel_A(b)(\chi)| = \sup\limits_{\chi\in
  X^*_p(A)}|\varphi(\chi)|=0,
$$
то есть $b\in N(p)$. Но $\chi(a)\neq0$ и
$\chi(b)=\varphi(\chi)=1\neq0$. Значит, $\chi(ab)\neq0$, так что и
$ab\neq0$. Получается, что $a\cdot N(p)\ne\{0\}$. Однако это
противоречит выбору $p$. Тем самым, наше предложение неверно, и
$\chi\in \Chi^*_p(A)$. Теперь имеем $\gel_A(a)(\chi)=1$ и
$\gel_A(a)=0$ вне $\Chi^*_p(A)$, а это и означает, что
$\mathcal{X}^*_c(A)\in \ob\widetilde{\mathcal{L}c}$. Теорема
доказана.
\end{proof}

Сочетая только что доказанную теорему с результатами предыдущего
параграфа, можно получить ещё целый ряд теорем. Сформулируем
некоторые из них.

\begin{thm}
Категория, дуальная к категории $C^*$-бочечных $lc^*$-алгебр
эквивалентна категории локально компактных $\mu$-пространств.
\end{thm}

\begin{thm}
Категория, дуальная к категории борнологических $lc^*$-алгебр,
эквивалентна категории вещественнокомпактных локально компактных
про\-ст\-ранств.
\end{thm}

\begin{thm}
Категория, дуальная к категории метризуемых $lc^*$-алгебр,
эквивалентна категории хемикомпактных локально компактных (или, что
--- то же самое, $\sigma$-компактных локально компактных)
пространств.
\end{thm}

\begin{thm}
Каждое $\sigma$-компактное локально компактное пространство является
$\mu$-пространством.
\end{thm}

\chapter{Тензорные произведения $LC^*$-алгебр}

\section{Проективное тензорное произведение $LC^*$-алгебр}

Пусть $E$, $F$ --- локально выпуклые пространства.
Напомним~\cite[Т.~2]{b16},~\cite{b37}, что на тензорном произведении
$E\otimes F$ векторных пространств $E$, $F$ (без топологии)
существует сильнейшая локально выпуклая топология, относительно
которой каноническое билинейное отображение
$$
  \otimes\colon  E\times F \to E\otimes F
$$
непрерывно. Она называется проективной локально выпуклой топологией
тензорного произведения и обозначается через $\pi$. Наделённое этой
топологией пространство $E\otimes F$ называется проективным
тензорным произведением локально выпуклых пространств $E$, $F$ и
обозначается через $E\mathbin{\underset{\pi}{\otimes}} F$. Пополнение пространства $E\mathbin{\underset{\pi}{\otimes}} F$ называется полным проективным тензорным произведением локально выпуклых пространств
$E$, $F$ и обозначается через $E\mathbin{\hat{\underset{\pi}{\otimes}}} F$.

Если $A_1$, $A_2$ --- локально выпуклые алгебры с совместно
непрерывным умножением, то пространства
$A_1\mathbin{\underset{\pi}{\otimes}}A_2$ и $A_1\mathbin{\hat{\underset{\pi}{\otimes}}} A_2$ естественным образом превращаются в топологические алгебры, в которых умножение определяется условием 
$$
  (x_1\otimes x_2)\cdot(y_1\otimes y_2) = x_1y_1\otimes x_2y_2
$$
для элементарных тензоров и затем канонически распространяется на
остальные элементы. Если $A_1$, $A_2$ локально мультипликативно
выпуклы, то $A_1\mathbin{\underset{\pi}{\otimes}}A_2$ и
$A_1\mathbin{\hat{\underset{\pi}{\otimes}}}A_2$ --- тоже~\cite{b24}. Наконец,
наличие на алгебрах $A_1$, $A_2$ непрерывной инволюции позволяет
наделить непрерывной инволюцией и алгебры
$A_1\mathbin{\underset{\pi}{\otimes}}A_2$, $A_1\mathbin{\hat{\underset{\pi}{\otimes}}}
A_2$: полагаем
$$
  (x_1\otimes x_2)^*=x^*_1\otimes x^*_2
$$
для элементарных тензоров и канонически продолжаем на остальные
элементы.

\begin{lemma}
Пусть $A_1$, $A_2$ --- $LC^*$-алгебры. Тогда все непрерывные
характеры на алгебрах $A_1\mathbin{\underset{\pi}{\otimes}}A_2$ и
$A_1\mathbin{\hat{\underset{\pi}{\otimes}}}A_2$ --- эрмитовы.
\end{lemma}

\begin{proof} Пусть $\chi$ --- непрерывный характер
алгебры $A_1\mathbin{\underset{\pi}{\otimes}}A_2$. Тогда отображения
$$\chi_1\colon
\left\{%
\begin{array}{ll}
    A_1\rightarrow\mathbb{C} \\
    x_1\mapsto\chi(x_1\otimes e_2), \\
\end{array}%
\right. \chi_2\colon \left\{%
\begin{array}{ll}
    A_2\rightarrow\mathbb{C} \\
    x_2\mapsto\chi(e_1\otimes x_2),\\
\end{array}%
\right.$$ являются непрерывными характерами алгебр $A_1$, $A_2$
соответственно. При этом для любых $x_1\in A_1$, $x_2\in A_2$ имеем:
$$
  \chi(x_1 \otimes x_2)=\chi(x_1\otimes e_2\cdot e_1\otimes x_2)=
  \chi(x_1\otimes e_2)\cdot\chi(e_1\otimes x_2)=\chi_1(x_1)\chi_2(x_2).
$$
Отсюда, в силу линейности $\chi$, для любого $z=\sum x_{1i}\otimes
x_{2i}$ получаем
$$
  \chi(z)=\sum\chi_1(x_1)\cdot\chi_2(x_2).
$$
Но на $LC^*$-алгебрах все непрерывные характеры --- эрмитовы (гл. 1,
$\S20$, предл. $(LC^*)$, п.~$2^\circ$). Поэтому $\chi_1$,
$\chi_2$~--- эрмитовы, а тогда и $\chi$ --- тоже эрмитов. Этим
доказано первое утверждение. Второе следует из первого в силу
непрерывности инволюции.
\end{proof}

Пусть
$$
  m_{\mathbb{C}}\colon \mathbb{C}\times\mathbb{C}\rightarrow\mathbb{C}
$$
--- операция умножения комплексных чисел. Это --- непрерывное
билинейное отображение, следовательно, ему канонически соответствует
непрерывное линейное отображение
$$
  \widetilde{m}_{\mathbb{C}}\colon \mathbb{C}\mathbin{\underset{\pi}{\otimes}}
  \mathbb{C}\rightarrow\mathbb{C},
$$
для которого $\widetilde{m}_{\mathbb{C}}(z\otimes t)=z\cdot t$.
Поскольку и $\mathbb{C}\otimes\mathbb{C}$, и $\mathbb{C}$ являются
одномерными топологическими инволютивными алгебрами над полем
$\mathbb{C}$, $\widetilde{m}_{\mathbb{C}}$ есть изоморфизм
топологических инволютивных алгебр. Значит,
$\hat{\mathbb{C}\mathbin{\underset{\pi}{\otimes}}\mathbb{C}}$ совпадает с
$\mathbb{C}\mathbin{\underset{\pi}{\otimes}}\mathbb{C}$, а
$\widetilde{m}_{\mathbb{C}}$ --- с его продолжением на
$\hat{\mathbb{C}\mathbin{\underset{\pi}{\otimes}}\mathbb{C}}$.

Пусть $A_1$, $A_2$ --- локально выпуклые алгебры. Каждая пара
$(\chi_1, \chi_2)\in\Chi_c(A_1)\times\Chi_c(A_2)$ порождает
непрерывные характеры
$$
  A_1\mathbin{\underset{\pi}{\otimes}}A_2\overset{\chi_1\otimes\chi_2}%
  {\longrightarrow}\mathbb{C}\otimes\mathbb{C}\overset%
  {\widetilde{m}_{\mathbb{C}}}{\longrightarrow}\mathbb{C},
$$
$$
  A_1\mathbin{\hat{\underset{\pi}{\otimes}}}A_2\overset%
  {\chi_1\mathbin{\hat{\underset{\pi}{\otimes}}}\chi_2}%
  {\longrightarrow}\mathbb{C}\otimes\mathbb{C}\overset%
  {\widetilde{m}_{\mathbb{C}}}{\longrightarrow}\mathbb{C}.
$$

\begin{prop}
$1^\circ$ Пусть $A_1$, $A_2$ --- отделимые $LC^*$-алгебры. Тогда
отображение
$$
  J \colon
    \left\{%
    \begin{array}{ll}
    \Chi^*_c(A_1)\times\Chi^*_c(A_2)\rightarrow\Chi^*_c(A_1\mathbin{\underset{\pi}{\otimes}} A_2)  \\
    (\chi_1,\chi_2)\mapsto\widetilde{m}_{\mathbb{C}}\circ(\chi_1\otimes\chi_2) \\
    \end{array}%
    \right.
$$
есть гомеоморфизм.

$2^\circ$ Пусть $A_1$, $A_2$ --- $K^*$-алгебры. Тогда:

а) отображение
$$
  \hat{J}\colon
  \left\{%
  \begin{array}{ll}
    \Chi^*_c(A_1)\times\Chi^*_c(A_2)\rightarrow\Chi^*_c(A_1\mathbin{\hat{\underset{\pi}{\otimes}}} A_2)  \\
    (\chi_1,\chi_2)\mapsto\overline{m}_{\mathbb{C}}\circ(\chi_1\mathbin{\hat{\underset{\pi}{\otimes}}}\chi_2) \\
  \end{array}%
  \right.
$$
есть гомеоморфизм;

б) отображение сужения
$$
  r\colon \Chi^*_c(A_1\mathbin{\hat{\underset{\pi}{\otimes}}} A_2)%
  \rightarrow\Chi^*_c(A_1\mathbin{\underset{\pi}{\otimes}}A_2)
$$
есть гомеоморфизм.
\end{prop}

\begin{proof} В силу леммы (1.1)
\begin{align*}
  \Chi^*_c(A_1\mathbin{\underset{\pi}{\otimes}}
  A_2)&=\Chi_c(A_1\mathbin{\underset{\pi}{\otimes}}A_2)\\
  \Chi^*_c(A_1\mathbin{\hat{\underset{\pi}{\otimes}}}A_2)&=
  \Chi_c(A_1\mathbin{\hat{\underset{\pi}{\otimes}}}A_2)
\end{align*}
По теореме Маллиоса~\cite{b24} $J$ есть гомеоморфизм
$\Chi_c(A_1)\times\Chi_c(A_2)$ на $\Chi_c(A_1\mathbin{\underset{\pi}{\otimes}}
A_2)$, а по теореме Лопушанского~\cite{b21} $\hat{J}$ есть
гомеоморфизм $\Chi_c(A_1)\times\Chi_c(A_2)$ на 
$\Chi_c(A_1 \mathbin{\hat{\underset{\pi}{\otimes}}} A_2)$. Наконец, ясно, что $r=J\circ\hat{J}^{-1}$. Предложение доказано.
\end{proof}

Пусть $A_1$, $A_2$ --- $LC^*$-алгебры. В конце $\S2$ мы покажем, что
на $A_1\mathbin{\underset{\pi}{\otimes}}A_2$ существует сильнейшая
$LC^*$-топология, относительно которой каноническое отображение
$\otimes$ непрерывно. Ею является $LC^*$-топология, ассоциированная
с топологией $\pi$, то есть топология, определяемая совокупностью
всех непрерывных $C^*$-преднорм на $A_1\mathbin{\underset{\pi}{\otimes}}A_2$.

\begin{dfn}
Алгебру $A_1\otimes A_2$, наделённую $LC^*$-топологией,
ассоциированной с топологией $\pi$, будем называть проективным
тензорным произведением $LC^*$-алгебр $A_1$, $A_2$ и обозначать
через $A_1\mathbin{\underset{\ \pi^*}{\otimes}} A_2$. Алгебру, получаемую из
$A_1\mathbin{\underset{\ \pi^*}{\otimes}} A_2$ отделимым пополнением, будем
называть полным проективным произведением $LC^*$-алгебр $A_1$, $A_2$
и обозначать через $A_1\mathbin{{\underset{\ \pi^*}{\bar\otimes}}} A_2$.
\end{dfn}

Выделим для дальнейших ссылок:

\begin{prop}
Пусть $A_1$, $A_2$ --- $LC^*$-алгебры. Тогда топология $\pi^*$ на
\mbox{$A_1\otimes A_2$} является сильнейшей из тех $LC^*$-топологий,
относительно которых отображение $\otimes$ непрерывно.
\end{prop}

Важно знать, в каких случаях $A_1\otimes A_2$ можно считать
подалгеброй в $A_1{\underset{\ \pi^*}{\bar\otimes}}A_2$. Ответ даёт

\begin{prop}
Пусть $A_1$, $A_2$ --- $LC^*$-алгебры. Тогда следующие утверждения
равносильны:

(а) $A_1$ и $A_2$ отделимы;

(б) $A_1\mathbin{\underset{\ \pi^*}{\otimes}} A_2$ отделима;

(в) $A_1\otimes A_2$ является подалгеброй в
$A_1\mathbin{\underset{\ \pi^*}{\bar\otimes}}A_2$.
\end{prop}

\begin{proof} Поскольку
$A_1\mathbin{\underset{\ \pi^*}{\bar\otimes}}A_2$ получается из
$A_1\mathbin{\underset{\ \pi^*}{\otimes}}A_2$ {\em отделимым} пополнением,
равносильность (б) и (в) очевидна. Покажем, что (a) $\Eq$ (б).

(а) $\Rightarrow$ (б). Пусть $A_1$, $A_2$ отделимы и, значит, строго
$*$-полупросты (гл. 1, $\S20$, предл. ($LC^*$), п.~$5^\circ$).
Покажем, что и $A_1\mathbin{\underset{\ \pi^*}{\otimes}}A_2$ строго
$*$-полупроста (и, значит, отделима). Пусть
$$
  z = \sum_{i=1}^n a_{1i}\otimes a_{2i} \text{ ---}
$$
отличный от нуля элемент в $A_1\otimes A_2$. Найдём
непрерывный эрмитов характер $\chi$ на $A_1\otimes A_2$, для
которого $\chi(z)\ne 0$. Не уменьшая общности, можем считать, что
векторы $a_{21}$, \ldots, $a_{2n}$ линейно независимы. Далее,
поскольку $z\ne 0$, среди $a_{11}$, \ldots, $a_{1n}$ есть отличный
от нуля, скажем, $a_{11}$. В силу строгой $*$-полупростоты алгебры
$A_1$ существует $\chi_1\in\Chi^*_c(A_1)$ такой, что $\chi_1(a_{11})
\ne 0$. Поскольку векторы $a_{21}$, \ldots, $a_{2n}$ линейно
независимы,
$$
  z_2 = \sum\chi_1(a_{1i})a_{2i}
$$
есть отличный от нуля элемент алгебры $A_2$. В силу строгой
$*$-полупростоты $A_2$ существует $\chi_2\in\Chi^*_c(A_2)$ такой,
что $\chi_2(z_2)\ne 0$. Положив $\chi = J(\chi_1,\ \chi_2),$ получим
$$
  \chi(z) = J(\chi_1,\ \chi_2)(z) = \sum\chi_1(a_{1i})\chi_2(a_{2i}) =
  \chi_2(z_2)\neq0.
$$
Это и требовалось.

(б) $\Rightarrow$ (a). Пусть $A_1\mathbin{\underset{\ \pi^*}{\otimes}} A_2$
отделима. Покажем, что и $A_1$ отделима. Пусть $a_1\in A_1$ и
$a_1\neq0$. Тогда и $a_1\otimes e_2\neq0$. В силу строгой
$*$-полупростоты $A_1\mathbin{\underset{\ \pi^*}{\otimes}} A_2$ существует
$\chi\in \Chi^*_c(A_1\mathbin{\underset{\ \pi^*}{\otimes}}A_2)$ такой, что
$\chi(a_1\otimes e_2)\neq0$. Из определения (1.3) ясно, что
$$
  \Chi^*_c(A_1\mathbin{\underset{\ \pi^*}{\otimes}}A_2) =
  \Chi^*_c(A_1\mathbin{\underset{\pi}{\otimes}}A_2).
$$
Отсюда, в силу (1.2), обнаруживаем существование $(\chi_1,\
\chi_2)\in\Chi^*_c(A_1)\times\Chi^*_c(A_2)$ с $J(\chi_1,\chi_2) =
\chi$ и получаем
$$
  \chi_1(a_1) = \chi_1(a_1)\chi_2(e_2) = J(\chi_1,\chi_2)(a_1\otimes e_2)
  = \chi(a_1\otimes e_2)\ne 0.
$$
Тем самым, $A_1$ строго $*$-полупроста и потому --- отделима. Для
$A_2$ рассуждаем точно так же.
\end{proof}

\begin{prop}
Пусть $A_1$, $A_2$, $B_1$, $B_2$ --- $LC^*$-алгебры, $u_1\colon
A_1\rightarrow B_1$, $u_2\colon A_2\rightarrow B_2$ --- непрерывные
эрмитовы морфизмы. Тогда существуют однозначно определенные
непрерывные эрмитовы морфизмы
$$
  u_1\otimes u_2\colon A_1\mathbin{\underset{\ \pi^*}{\otimes}}%
  A_2\longrightarrow B_1\mathbin{\underset{\ \pi^*}{\otimes}}B_2,
$$
$$
  u_1\mathbin{\bar\otimes}u_2\colon A_1\mathbin{\underset{\ \pi^*}%
  {\bar\otimes}}A_2\longrightarrow B_1\mathbin{\underset{\ \pi^*}{\bar\otimes}}B_2,
$$
для которых при всех $a_1\in A_1$, $a_2\in A_2$ выполняются
равенства
$$
  (u_1\mathbin{\bar\otimes}u_2)(a_1\otimes a_2) = u_1(a_1)u_2(a_2),
$$
$$
  (u_1\mathbin{\bar\otimes}u_2)(a_1\otimes a_2) = u_1(a_1)u_2(a_2).
$$
\end{prop}

\begin{proof} Поскольку $u_1$ и $u_2$ являются
непрерывными линейными отображениями, из свойств проективного
тензорного произведения локально выпуклых пространств вытекает
существование непрерывного линейного отображения
$$
  u_1\otimes u_2\colon A_1\mathbin{\underset{\pi}{\otimes}}A_2\longrightarrow B_1
  \mathbin{\underset{\pi}{\otimes}}B_2,
$$
удовлетворяющего первому равенству. Тривиально проверяется, что
это~--- эрмитов морфизм. Покажем, что он является непрерывным
отображением из $A_1\mathbin{\underset{\ \pi^*}{\otimes}}A_2$ в
$B_1\mathbin{\underset{\ \pi^*}{\otimes}}B_2$. Пусть $q$ --- $\pi^*$-непрерывная
$C^*$-преднорма на $B_1\otimes B_2$. Тогда она и $\pi$-непрерывна.
Значит, композиция $q\circ(u_1\otimes u_2)$ является
$\pi$-непрерывной $C^*$-преднормой на $A_1\otimes A_2$. В силу
определения топологии $\pi^*$ эта преднорма $\pi^*$-непрерывна.
Итак, для любой $\pi^*$-непрерывной $C^*$-преднормы $q$ на
$B_1\otimes B_2$ композиция $q\circ(u_1\otimes u_2)$ является
$\pi^*$-непрерывной $C^*$-преднормой на $A_1\otimes A_2$. Это и
означает, что $u_1\otimes u_2$ является непрерывным отображением из
$A_1\mathbin{\underset{\ \pi^*}{\otimes}}A_2$ в
$B_1\mathbin{\underset{\ \pi^*}{\otimes}}B_2$. Первое утверждение доказано.
Второе выводится из первого с помощью универсального свойства
отделимых пополнений.
\end{proof}

\begin{prop}
Пусть $A$ --- $lc^*$-алгебра, соотв., $K^*$-алгебра. Тогда
существует ровно один непрерывный эрмитов морфизм
$\widetilde{m}_A\colon A\underset{\ \pi^*}{\otimes}A\rightarrow A$,
соотв., $\bar{m}_A\colon
A\mathbin{\underset{\ \pi^*}{\bar\otimes}} A\rightarrow A$ такой, что для всех
$x,\ y \in A$ выполняется равенство $\widetilde{m}_A(x\otimes
y)=x\cdot y$, соотв., $\hat{m}_A(x\otimes y)=x\cdot y$.
\end{prop}

\begin{proof} Поскольку $A$ --- $LC^*$-алгебра, операция умножения
$$
  m_A\colon  A\times A\rightarrow A
$$
является непрерывным билинейным отображением (гл.1, $\S20$, предл.
$(LC^*)$, п.~$1^\circ$). Поэтому существует ровно одно непрерывное
линейное отображение
$$
  \widetilde{m}_A\colon A\mathbin{\underset{\pi}{\otimes}} A \to A,
$$
удовлетворяющее первому условию. Тривиально проверяется, что это ---
эрмитов морфизм, после чего доказательство завершается так же, как
доказательство предложения (1.6).
\end{proof}

\begin{prop}
$1^\circ$ Для любых $LC^*$-алгебр $A_1, A_2, A$ отображение
$$
  \gamma^{A_1,A_2}_{A}\colon
  \left\{%
   \begin{array}{ll}
    \mathcal{L}\mathcal{C}^*(A_1, A)\times\mathcal{L}\mathcal{C}^*(A_2, A)\rightarrow\mathcal{L}\mathcal{C}^*(A_1\otimes A_2, A)\\
    (u_1, u_2)\mapsto\widetilde{m}_A\circ(u_1\otimes u_2)\\
   \end{array}%
  \right.
$$
является биекцией.

$2^\circ$ Для любых $K^*$-алгебр $A_1$, $A_2$, $A$ отображение

$$
  \bar{\gamma}^{A_1,A_2}_{A}\colon \left\{%
  \begin{array}{ll}
    \mathcal{K}^*(A_1, A)\times\mathcal{K}^*(A_2, A)\rightarrow\mathcal{K}^*(A_1\mathbin{\underset{\ \pi^*}{\bar\otimes}}A_2, A)\\
    (u_1, u_2)\mapsto\overline{m}_A\circ(u_1\mathbin{\bar\otimes} u_2)\\
  \end{array}%
  \right.
$$
является биекцией.
\end{prop}

\begin{proof} 1. Инъективность устанавливается так же,
как и в алгебраическом случае (без топологий). Докажем
сюръективность. Пусть $u\colon A_1 \mathbin{\underset{\ \pi^*}{\otimes}} A_2 \to
A$ --- произвольный непрерывный эрмитов морфизм. Рассмотрим
отображения
$$
  i_1\colon \left\{%
  \begin{array}{ll}
    A_1\rightarrow A_1\mathbin{\underset{\ \pi^*}{\otimes}}A_2 \\
    x_1\mapsto x_1\otimes e_2 \\
  \end{array}%
  \right.,\qquad i_2\colon\left\{%
  \begin{array}{ll}
    A_2\rightarrow A_1\mathbin{\underset{\ \pi^*}{\otimes}}A_2 \\
    x_2\mapsto e_1\otimes x_2 \\
  \end{array}%
  \right..
$$
Первое есть композиция гомеоморфизма $x_1\mapsto(x_1,\ e_2)$ из
$A_1$ в $A_1\times A_2$ и отображения~$\otimes$, второе~---
композиция гомеоморфизма $x_2\mapsto(e_1,\ x_2)$ из $A_2$ в
$A_1\times A_2$ и того же отображения. Поскольку $\otimes$~---
непрерывное отображение, $i_1$ и $i_2$ непрерывны. Кроме того, ясно,
что они являются эрмитовыми морфизмами. Все это показывает, что
композиции $u_1 = u\circ i_1$ и $u_2 = u\circ i_2$ являются
непрерывными эрмитовыми морфизмами из $A_1$, соотв., $A_2$, в $A$.
Стандартная проверка показывает, что $u = \gamma(u_1,\ u_2)$.
Доказательство п.~$1^\circ$ завершено.

2. Пусть $A_1$, $A_2$~--- $K^*$-алгебры. Тогда они отделимы, так что
$A_1\mathbin{\underset{\ \pi^*}{\otimes}} A_2$ является плотной подалгеброй в $A_1\mathbin{\underset{\ \pi^*}{\bar\otimes}}A_2$ (см. предл. (1.5)). Если при этом и $A$ является $K^*$-алгеброй, то каждый непрерывный эрмитов морфизм из  $A_1\mathbin{\underset{\ \pi^*}{\otimes}}A_2 $ в $A$
однозначно продолжается до непрерывного  эрмитова морфизма из
$A_1\mathbin{\underset{\ \pi^*}{\bar\otimes}} A_2$ в $A$. Иначе говоря, в этом случае $\mathcal{K}^*(A_1\mathbin{\underset{\ \pi^*}{\bar{\otimes}}}A_2,
A)=\mathcal{LC}^*(A_1\mathbin{\underset{\ \pi^*}{\otimes}}A_2, A)$. Ясно также, что $\mathcal{K}^*(A_i, A)=\mathcal{LC}^*(A_i, A),$ $i=1,\ 2$. Применяя п.~$1^\circ$, получаем требуемое.
\end{proof}

\begin{cor}
$1^\circ$ $A_1\mathbin{\underset{\ \pi^*}{\otimes}}A_2$ есть прямая сумма $A_1$
и $A_2$ в категории $\mathcal{L}\mathcal{C}^*$.

$2^\circ$ $A_1\mathbin{\underset{\ \pi^*}{\bar{\otimes}}}A_2$ есть прямая сумма $A_1$ и $A_2$ в категории $\mathcal{K^*}$.
\end{cor}

\begin{proof} 1. Пусть
$u_k\in\mathcal{L}\mathcal{C}^*(A_k,\ A)$, $k = 1,\ 2$. Положим $u =
\gamma(u_1,\ u_2)$ и покажем, что $u\circ i_k = u_k$, $k = 1,\ 2$.
Имеем:
$$
  (u\circ i_1)(x_1) = \gamma(u_1,\ u_2)(x_1\otimes e_2) = u_1(x_1)\cdot
  u_2(e_2) = u_1(x_1).
$$
Это показывает, что $u\circ i_1=u_1$. Точно так
же проверяется и второе равенство.

2. Пусть $u_k\in\mathcal{K}^*(A_k,\ A)$, $k = 1,\ 2$. Тогда
$A_1\mathbin{\underset{\ \pi^*}{\otimes}}A_2$ является плотной подалгеброй в
$A_1\mathbin{\underset{\ \pi^*}{\hat{\otimes}}}A_2,$ a $\hat{m}_A$~---
продолжением $\tilde{m}_A$. Остается применить п.~$1^\circ$.
\end{proof}

\section{Инъективное тензорное произведение $LC^*$-алгебр}

Сначала~--- напоминания. Пусть $E$~--- локально выпуклое
пространство с сопряжённым $E'$. Полярой множества $U\subset E$
соотв., $V\subset E'$ называется множество
$$
  U^\circ = \{l\in E'\colon(\forall x\in U)(|\langle l,\ x\rangle|\leqslant 1)\},
$$
соотв.,
$$
  V^\circ = \{x\in E\colon (\forall l\in V)(|\langle l,\ x\rangle|)\leqslant 1\}.
$$

Пусть $E_1$, $E_2$ --локально выпуклые пространства с сопряжёнными
$E'_1$, $E'_2$. Каждая пара $(l_1,\ l_2)\in E'_1\times E'_2$
порождает линейную форму
$$
  E_1\otimes E_2\overset{l_1\otimes\, l_2}{\longrightarrow}%
  \Cbb\otimes\Cbb\overset{\widetilde{m}_\Cbb}%
  {\longrightarrow}\Cbb,
$$
которую обычно, допуская вольность речи, обозначают просто через
$l_1\otimes l_2$. Продолжая в том же духе, для любых $V_1\subset
E'_1$, $V_2\subset E'_2$ полагают
$$
  V_1\otimes V_2 = \{l_1\otimes l_2\colon l_1\in V_1,\ l_2\in V_2\}.
$$
(Разумеется, в подобных случаях нужно соблюдать осторожность и не
принимать $V_1\otimes V_2$ за тензорное произведение.) В этих
обозначениях инъективное тензорное произведение локально выпуклых
пространств определяется следующим образом.

Пусть $E_1$, $E_2$~--- отделимые локально выпуклые пространства, а
$U_1$, $U_2$ пробегают какие-нибудь фундаментальные системы
окрестностей нуля в $E_1$, $E_2$ соответственно. Тогда множества
$(U^\circ_1\otimes U^\circ_2)^\circ$ образуют фундаментальную
систему замкнутых абсолютно выпуклых окрестностей нуля некоторой
локально выпуклой топологии на $E_1\otimes E_2$. Эту топологию
называют инъективной топологией тензорного произведения локально
выпуклых пространств $E_1$, $E_2$ и обозначают через $\varepsilon$.
Наделённое ею пространство $E_1\otimes E_2$ называют инъективным
тензорным произведением локально выпуклых пространств $E_1$, $E_2$ и
обозначают через $E_1\underset{\varepsilon}{\otimes} E_2$. Его
пополнение называют полным инъективным тензорным произведением
пространств $E_1, E_2$ и обозначают через
$E_1\hat{\underset{\varepsilon}{\otimes}}E_2$.

Пусть $p_1$, $p_2$~--- непрерывные преднормы на $E_1$, $E_2$
соответственно, $U_i = \{x_i\in E_i\colon p_i(x_i)\leqslant 1\}$,
$i=1,2$. Для любого $z = \sum\limits^n_{j=1}x_{1j}\otimes x_{2j}$ и
любого $\omega = \omega_1\otimes\omega_2\in U^\circ_1\otimes
U^\circ_2$ имеем
$$
  \langle w,\ z\rangle =
  \sum\limits^n_{j=1}\langle w_1,\ x_{1j}\rangle\cdot
  \langle w_2,\ x_{2j}\rangle.
$$
В силу этого условие
$$
  (\forall w\in U^{\circ}_1\otimes U^\circ_2)(|\langle w,\ z\rangle|\leqslant 1)
$$
оказывается равносильным тому, что
$$
  \sup\left\{\abs{\sum\langle w_1,\ x_{1j}\rangle\cdot\langle w_2, x_{2j}\rangle}
  \colon w_1\in U^\circ_1,\ w_2\in U^\circ_2\right\}\leqslant 1.
$$
Полагая
$$
  (p_1\underset{\varepsilon}{\otimes}p_2)(z) =
  \sup\left\{\abs{\sum\langle w_1,\ x_{1j}\rangle\cdot
  \langle w_2,\ x_{2j}\rangle}\colon
   w_1\in U^\circ_1,\ w_2\in U^\circ_2\right\},
$$
получаем:
$$
  (U^\circ_1\otimes U^\circ_2)^\circ = \{z\in E_1\otimes E_2\colon
  (p_1\underset{\varepsilon}{\otimes}p_2)(z)\leqslant 1\}.
$$

Таким образом, топология $\varepsilon$ на $E_1\otimes E_2$
порождается совокупностью преднорм вида
$p_1\underset{\varepsilon}{\otimes}p_2$, где $p_1$, соотв., $p_2$
пробегает какую-нибудь совокупность преднорм, задающих топологию на
$E_1$, соотв., $E_2$.

Следующий результат даёт описание $\varepsilon$-произведения
$C^*$-преднорм, более удобное для применений в теории $LC^*$-алгебр.

\begin{lemma}
Пусть $p_1$, $p_2$~--- непрерывные $C^*$-преднормы на
$LC^*$-алгебрах $A_1$, $A_2$ соответственно и $K_1 =
\Chi^*_{p_1}(A_1)$, $K_2 = \Chi^*_{p_2}(A_2)$. Тогда
$$
  p_1\underset{\varepsilon}{\otimes}p_2 = p^{J[K_1\times K_2]}.
$$
\end{lemma}

\begin{proof} Напомним, что для любого $a_i\in A_i$
имеем
$$
  p^{K_i}(a_i) = p_i(a_i) = \sup\{|\langle l_i, a_i\rangle |\colon l_i\in
  V^\circ_{p_i}\},
$$
где
$$
  p^{K_i}(a_i) = \sup\{|\chi_i(a_i)|\colon \chi_i\in K_i\}
$$
и
$$
  V_{p_i} = \{a_i\in A_i\colon  p_i(a_i)\leqslant 1\},\  i = 1, 2,
$$
(см. гл.~2, предл. (1.1) и~\cite[Т.~1]{b16}). Используя это, для
любого
$$
  z = \sum\limits^n_{j=1}a_{1j}\otimes a_{2j}\in A_1\otimes A_2
$$
получаем (мы пишем \glqq$\sum$\grqq вместо
\glqq$\sum\limits^n_{j=1}$\grqq):
\begin{align*}
  (p_1\underset{\varepsilon}{\otimes}p_2)(z)\leqslant 1
  &\Eq (\forall l_1\in V^\circ_{p_1})(\forall l_2\in
  V^\circ_{p_2})\left(\abs{\sum\langle l_1,\ a_{1j}\rangle \langle l_2,
  a_{2j}\rangle}\leqslant 1\right)\\
  &\Eq (\forall l_1\in V^\circ_{p_1})(\forall l_2\in
  V^\circ_{p_2})\left(\abs{\langle l_2,\ \sum\langle l_1,\
  a_{1j}\rangle a_{2j}\rangle}\leqslant 1\right)\\
  &\Eq (\forall l_1\in V^\circ_{p_1})\left(p_2\left(\sum
  \langle l_1,\ a_{1j}\rangle a_{2j}\right)\leqslant 1\right)\\
  &\Eq (\forall l_1\in V^\circ_{p_1})\left(p^{K_2}\left(\sum\langle l_1,
  a_{1j}\rangle a_{2j}\right)\leqslant 1\right)\\
  &\Eq (\forall l_1\in V^\circ_{p_1})(\forall\chi_2\in K_2)\left(\abs{\langle\chi_2,\ \sum\langle l_1, a_{1j}\rangle a_{2j}\rangle}\leqslant 1\right)\\
  &\Eq (\forall\chi_2\in K_2)(\forall l_1\in V^\circ_{p_1})\left(\abs{\sum\langle l_1, a_{1j}\rangle \langle\chi_2, a_{2j}\rangle}\leqslant 1\right)\\
  &\Eq (\forall\chi_2\in K_2)(\forall l_1\in V^\circ_{p_1})\left(\abs{\langle l_1,\ \sum\langle\chi_2,\ a_2j\rangle a_{1j}\rangle} \leqslant 1\right)\\
  &\Eq (\forall\chi_2\in K_2)\left(p_1\left(\sum\langle\chi_2,\ a_{2j}\rangle a_{1j}\right)\leqslant 1\right)\\
  &\Eq (\forall\chi_2\in K_2)\left(p^{K_1}\left(\sum\langle\chi_2,a_{2j}\rangle a_{1j}\right)\leqslant 1\right)\\
  &\Eq (\forall\chi_2\in K_2)(\forall \chi_1\in K_1)\left(\abs{\sum\langle\chi_1, a_1j\rangle \langle\chi_2, a_{2j}\rangle}\leqslant 1\right)\\
  &\Eq p^{J[K_1\times K_2]}(z)\leqslant 1.
\end{align*}
Это и означает, что $p_1\underset{\varepsilon}{\otimes}p_2 =
p^{J[K_1\times K_2]}$.
\end{proof}

Вспоминая, что каждая преднорма вида $p^S$ является
$C^*$-преднормой, получаем

\begin{prop}
Если $A_1$, $A_2$~--- $LC^*$-алгебра, то и
$A_1\underset{\varepsilon}{\otimes}A_2$~--- тоже.
\end{prop}

Теперь, $K_1$, $K_2$ из леммы (2.1) компактны, а $J$~--- непрерывно.
Значит $J[K_1\times K_2]$ есть компактное и потому замкнутое
подмножество в $\Chi^*(A_1\underset{\varepsilon}{\otimes}A_2)$.
Применяя предложение (1.6) из гл.2, выводим

\begin{prop}
$\Chi^*_{p_1\mathbin{\underset{\varepsilon}{\otimes}}\, p_2}(A_1\otimes
A_2)=J[\Chi^*_{p_1}(A_1)\times\Chi^*_{p_2}(A_2)].$
\end{prop}

Это показывает, что непрерывными характерами на алгебре
$A_1\mathbin{\underset{\varepsilon}{\otimes}} A_2$ являются в точности элементы множеств $J[K_1\times K_2]$ с $K_1\in \varkappa_1$,
$K_2\in\varkappa_2$, a членами
$\varkappa_{A_1\mathbin{\underset{\varepsilon}{\otimes}}A_2}$~--- в точности
замкнутые подмножества тех же множеств. Иначе говоря,
отображение~$J$ является компактологическим изоморфизмом
пространства $\Chi^*_c(A_1)\times\Chi^*_c(A_2)$ на пространство
$\Chi^*_{c}(A_1\mathbin{\underset{\varepsilon}{\otimes}}A_2).$ Получаем

\begin{prop}
Для любых $LC^*$-алгебр $A_1$, $A_2$ имеет место компактологический
изоморфизм $\Chi^*_{c}(A_1)\times\Chi^*_c(A_2)\cong
\Chi^*_c(A_1\mathbin{\underset{\varepsilon}{\otimes}}A_2)$.
\end{prop}

Поскольку каждый непрерывный эрмитов характер и каждая
$C^*$-преднорма однозначно продолжаются с
$A_1\mathbin{\underset{\varepsilon}{\otimes}}A_2$ на
$A_1\mathbin{\underset{\varepsilon}{\hat{\otimes}}}A_2$, отображение
сужения
$r\colon \mathcal{\Chi}^*_c(A_1\mathbin{\underset{\varepsilon}{\hat{\otimes}}}A_2)
\rightarrow\mathcal{\Chi}^*_c(A_1\mathbin{\underset{\varepsilon}{\otimes}}A_2)$
является компактологическим изоморфизмом. Компонируя его с $J$,
получаем компактологический изоморфизм
$\Chi^*_{c}(A_1)\times\Chi^*_c(A_2)$ на
$\mathcal{X}^*_c(A_1\mathbin{\underset{\varepsilon}{\hat{\otimes}}}A_2)$.
Применяя к нему функтор $\mathcal{C}$, приходим к следующему
результату.

\begin{thm}
Для любых отделимых $LC^*$-алгебр $A_1$, $A_2$ имеет место
топологический эрмитов изоморфизм
$$
  A_1\mathbin{\underset{\varepsilon}{\hat{\otimes}}}A_2\cong
  \mathcal{C}(\Chi^*_{c}(A_1)\times\Chi^*_c(A_2)).
$$
\end{thm}

\begin{proof}
Действительно, поскольку $A_1\mathbin{\underset{\varepsilon}{\hat{\otimes}}}A_2$~---
$K^*$-алгебра, имеет место топологический эрмитов изоморфизм
$A_1\mathbin{\underset{\varepsilon}{\hat{\otimes}}}A_2\cong\mathcal{C}
\mathcal{X}^*_c(A_1\mathbin{\underset{\varepsilon}{\hat{\otimes}}}A_2)$
(гл. 3, теорема (2.1)), и остается применить замечания,
предшествующие теореме.
\end{proof}

Доказательство той же теоремы (2.1) из главы 3 показывает, что
функторы $\mathcal{C}$ и $\mathcal{X}^*_c$ образую двойственность
между $\mathcal{K}$ и $\mathcal{K}^*$. Но в таком случае функтор
$\mathcal{C}$ преобразует произведения в прямые суммы, и
$\mathcal{C}(\mathcal{X}^*_c(A_1)\times\mathcal{X}^*_c(A_2))$ есть
прямая сумма $\mathcal{C}\mathcal{X}^*_c(A_1)$ и
$\mathcal{C}\mathcal{X}^*_c(A_2)$ в $\mathcal{K}^*$. При этом
$$
\mathcal{C}\mathcal{X}^*_c(A_i)\cong A_i, i=1, 2,
$$
так что
$\mathcal{C}(\mathcal{X}^*_c(A_1)\times\mathcal{X}^*_c(A_2))$ есть
также прямая сумма $A_1$ и $A_2$ в $\mathcal{K}^*.$ Применяя теорему
(2.5), заключаем, что
$A_1\mathbin{\underset{\varepsilon}{\hat{\otimes}}}A_2$ есть прямая сумма
$A_1$ и $A_2$ в $\mathcal{K}^*$. Вспоминая, что
$A_1\mathbin{\underset{\ \pi^* }{\bar\otimes}}A_2$ тоже является прямой суммой $A_1$ и $A_2$ в $\mathcal{K}^*$, получаем

\begin{cor}
$1^\circ$ $A_1\mathbin{\underset{\varepsilon}{\hat{\otimes}}}A_2$ есть прямая сумма $A_1$ и $A_2$ в $\mathcal{K}^*$.

$2^\circ$ Для любых $K^*$-алгебр $A_1$, $A_2$ алгебры
$A_1\mathbin{\underset{\ \pi^*}{\bar\otimes}} A_2$ и
$A_1\mathbin{\underset{\varepsilon}{\bar\otimes}} A_2$ топологически
эрмитово изоморфны.
\end{cor}

\begin{cor}
Пусть $A_1$, $A_2$~--- отделимые $LC^*$-алгебры. Тогда существует
только одна $LC^*$-топология $\tau\geqslant\varepsilon$ на
$A_1\otimes A_2$, относительно которой каноническое билинейное
отображение $\otimes$ непрерывно; это топология $\varepsilon$.
\end{cor}

В заключение покажем, что топология $\pi$ на тензорном произведении
$LC^*$-алгебр (даже~--- $K^*$-алгебр) не обязательно является
$LC^*$-топологией. Для этого нам понадобится теорема об инъективном
тензорном произведении алгебр непрерывных функций на тихоновских
пространствах. Чтобы её сформулировать, введем следующие
обозначения.

Пусть $T_1$, $T_2$~--- тихоновские пространства. Будем обозначать
через $T_1\boxtimes T_2$ $k_\Cbb$-пространство,
ассоциированное с тихоновским произведением пространств $T_1$ и
$T_2$, то есть множество $T_1\times T_2$, наделённое слабой
топологией, определяемой всеми функциями $T_1\times
T_2\longrightarrow\mathbb{C}$, непрерывными на каждом компакте вида
$K_1\times K_2$, где $K_1$, $K_2$~--- компактное подмножество в
$T_1$, $T_2$ соответственно.

\begin{thm}
Для любых тихоновских пространств $T_1$, $T_2$ имеет место
топологический эрмитов изоморфизм $$
C(T_1)\mathbin{\underset{\varepsilon}{\hat{\otimes}}} C(T_2)\cong
C(T_1\boxtimes T_2).
$$
\end{thm}

\begin{proof} В силу следствия (2.2) из гл. 3 имеет место
топологический эрмитов изоморфизм
$$
  \widehat{C(T)}\cong \mathcal{C}(\widetilde{T}),
$$
где $\widetilde{T}$~--- ассоциированное с $T$ компактологическое
пространство, а значок `\mbox{\;${\hat{}}$\;}' означает
пополнение. Применяя это замечание, теорему (2.1) из гл. 3 и
предыдущие результаты этого параграфа, получаем следующую цепочку
топологических эрмитовых изоморфизмов:
\begin{align*}
  C(T_1)\mathbin{\underset{\varepsilon}{\hat{\otimes}}} C(T_2)
&\cong
  \mathcal{C}\mathcal{X}^*_c(C(T_1)\mathbin{\underset{\varepsilon}{\hat{\otimes}}}C(T_2))\\
&\cong
  \mathcal{C}(\mathcal{X}^*_c(C(T_1)))\times
  \mathcal{X}^*_c(C(T_2))\\
&\cong
  \mathcal{C}(\mathcal{X}^*_c(\widehat{C(T_1)})\times
  \mathcal{X}^*_c(\widehat{C(T_2)})\\
&\cong
  \mathcal{C}(\widetilde{T_1}\times\widetilde{T_2})\cong C(T_1\boxtimes T_2).
\end{align*}
Этим теорема доказана.
\end{proof}

В. Дитрих в~\cite{b11} построил вполне регулярные $k$-пространства
$D$ и $H$ такие, что отображение сужения
$$
  r_1\colon \Chi^*_c(C(D)\mathbin{\underset{\varepsilon}{\hat{\otimes}}} C(H))
  \rightarrow\Chi^*_c(C(D)\mathbin{\underset{\varepsilon}{\otimes}} C(H))
$$
не является гомеоморфизмом. В то же время в силу предложения (1.2)
отображение сужения
$$
  r\colon \Chi^*_c(C(D)\mathbin{\hat{\underset{\pi}{\otimes}}} C(H))
  \rightarrow\Chi^*_c(C(D)\mathbin{\underset{\pi}{\otimes}} C(H))
$$
является гомеоморфизмом. Если бы топологии $\pi$ и $\varepsilon$ на
тензорном произведении $C(D)\otimes C(H)$ совпадали, эти два
утверждения не могли бы выполняться одновременно. Значит, топологии
$\pi$ и $\varepsilon$ на $C(D)\otimes C(H)$ различны. При этом
никакой $LC^*$-топологии $\tau > \varepsilon$, относительно которой
каноническое отображение $\otimes$ было бы непрерывно, на
$C(D)\otimes C(H)$ нет. Значит, топология $\pi$ не является
$LC^*$-топологией. Поскольку $C(D)$ и $C(H)$ являются
$K^*$-алгебрами, этим доказана

\begin{thm}
Существуют $K^*$-алгебры $A_1$, $A_1$, для которых топология $\pi$
на $A_1\otimes A_2$ не является $LC^*$-топологией.
\end{thm}

В силу теоремы (2.8) пространство
$\Chi^*_c(C(D)\mathbin{\underset{\varepsilon}{\hat{\otimes}}} C(H))$
гомеоморфно $D \boxtimes H$. Нетрудно показать, что пространство
$\Chi^*_c(C(D)\mathbin{\underset{\varepsilon}{\otimes}}C(H))$ гомеоморфно $D\times H$. Тем самым $D$ и $H$ являются вполне регулярными $k$-пространствами, тихоновское произведение которых не является $k_{\Cbb}$-пространством.

\chapter{Двойственность для полугрупп}

\section{$K^*$-бигебры}

$K^*$-бигеброй (соотв., $K^*$-алгеброй Хопфа) мы будем называть
$K^*$-алгебру $A$, наделённую непрерывным эрмитовым морфизмом
$\mu\colon A\rightarrow A\widehat{\underset{\varepsilon}{\otimes}}
A$ (и непрерывными эрмитовыми морфизмами $\varepsilon\colon
A\rightarrow \mathbb{C}$, $\sigma\colon  A\rightarrow A$), для
которых коммутативна диаграмма
$$
  \begin{diagram}
  \node{A}\arrow{s,t}{\mu}\arrow{e,t}{\mu}
  \node{A\ehatotimes A} \arrow{s,b}{\mu\ehatotimes I_A}\\
  \node{A\ehatotimes A} \arrow{e,b}{I_A\ehatotimes \mu}
  \node{A\ehatotimes A\ehatotimes A}
  \end{diagram}\eqno{(A^*)}
$$
(и диаграммы
$$
  \begin{diagram}
  \node[2]{A}\arrow{sw,t}{a\mapsto 1\otimes a}\arrow{s,b}{\mu}
  \arrow{se,t}{a\mapsto a\otimes 1}\\
  \node{\Cbb\ehatotimes A}
  \node{A\ehatotimes A} \arrow{e,b}{\eps\ehatotimes I_A}
  \arrow{w,b}{I_A\ehatotimes\eps}
  \node{A\ehatotimes\Cbb,}
  \end{diagram}\eqno{(E^*)}
$$
$$
  \begin{diagram}
  \node{A}\node{\Cbb}\arrow{w,t}{}
  \node{A}\arrow{w,t}{\eps}\arrow{s,t}{\mu}\arrow{e,t}{\eps}
  \node{\Cbb}\arrow{e}{}\node{A}\\
  \node{A\ehatotimes A}\arrow{n,t}{m}
  \node[2]{A\ehatotimes A}
  \arrow[2]{w,b}{I_A\ehatotimes\sigma}\arrow[2]{e,b}{I_A\ehatotimes\sigma}
  \node[2]{A\ehatotimes A,}\arrow{n,b}{m}
  \end{diagram}\eqno{(S^*)}
$$
где $m$~--- линейное отображение, канонически соответствующее
умножению на $A$, а стрелка \glqq $\mathbb{C}\rightarrow A$\grqq\ изображает единственный морфизм из $\mathbb{C}$ в $A$). Морфизмы $\mu$, $\varepsilon$ и $\sigma$ называются, соответственно коумножением, коединицей и кообращением (последний ещё называется сопряжением). Коединица и кообращение, если они существуют, определяются по $A$ и $\mu$ однозначно. По этой причине, называя $K^*$-алгебру Хопфа, мы будем обычно указывать лишь $K^*$-алгебру и коумножение.

Морфизмом из $K^*$-бигебры $(A_1,\ \mu_1)$ в $K^*$-бигебру $(A_2,
\mu_2)$ будем называть каждый непрерывный эрмитов морфизм $u:
A_1\rightarrow A_2$, согласующийся с коумножениями $\mu_1$ и
$\mu_2$, то есть такой, для которого диаграмма
$$
  \begin{diagram}
  \node{A_1}\arrow{s,t}{\mu_1}\arrow{e,t}{u}
  \node{A_2} \arrow{s,b}{\mu_2}\\
  \node{A_1\ehatotimes A_1} \arrow{e,b}{u_1\ehatotimes u_2}
  \node{A_2\ehatotimes A_2}
  \end{diagram}
$$
коммутативна. Если при этом $(A_1,\ \mu_1)$ и $(A_2,\ \mu_2)$~---
$K^*$-алгебры Хопфа, то $u$ согласуется также с коединицей и
кообращением. По этой причине под морфизмами $K^*$-алгебра Хопфа
понимаются просто морфизмы соответствующих $K^*$-бигебр.

\section{Двойственность для компактологических полугрупп}

Полугруппы, соотв., группы над категорией $\mathcal{K}$ регулярных
компактологических пространств будем называть компактологическими
полугруппами, соотв., группами.

Пусть $(S,\ m)$~--- компактологическая полугруппа. Применяя к
умножению $m: S\times S\rightarrow S$ кофунктор $\mathcal{C},$
получим непрерывный эрмитов морфизм
$$
  \mathcal{C}(m)\colon \mathcal{C}(S)\longrightarrow \mathcal{C}(S\times S),
$$
который, напомним, действует по правилу
$$
  \mathcal{C}(m)(\varphi)=\varphi\circ m
$$
для всех $\varphi\in\mathcal{C}(S)$. Обозначим через $\Phi$
канонический изоморфизм $K^*$-алгебры
\mbox{$\mathcal{C}(S)\mathbin{\underset{\varepsilon}{\hat\otimes}}\mathcal{C}(S)$}
на $K^*$-алгебру $\mathcal{C}(S\times S)$, определяемый условием
$$
  \Phi(\varphi\otimes\psi) = \varphi\cdot\psi
$$
для элементарных тензоров $\varphi\otimes\psi$ и продолжаемый
линейно и непрерывно на остальные элементы алгебры
$\mathcal{C}(S)\mathbin{\underset{\varepsilon}{\hat\otimes}}\mathcal{C}(S)$.
Положим
$$
  \mu_{\mathcal{C}(S)}=\Phi^{-1}\circ\mathcal{C}(m).
$$

\begin{lemma}
$1^\circ$ Для любой компактологической полугруппы соответственно,
группы $(S,\ m)$ пара $(\mathcal{C}(S),\ \mu_{\mathcal{C}(S)})$ есть
$K^*$-бигебра (соотв., $K^*$-алгебра Хопфа).

$2^\circ$ Для любого морфизма $u$ компактологических полугрупп
$\mathcal{C}(u)$ есть морфизм $K^*$-бигебр.
\end{lemma}

\begin{proof} Это тотчас же следует из того, что $\mathcal{C}$ есть
контравариантный функтор из категории $\mathcal{K}$ в категорию
$\mathcal{K}^*$ $K^*$-алгебр и того, что
$$
  \Phi\colon \mathcal{C}(S)\mathbin{\underset{\varepsilon}{\hat\otimes}}
  \mathcal{C}(S)\longrightarrow \mathcal{C}(S)
$$
есть изоморфизм $K^*$-алгебр.
\end{proof}

Пусть $(A,\ \mu)$~--- $K^*$-бигебра. Применяя к коумножению
$\mu\colon A\rightarrow A\mathbin{\underset{\varepsilon}{\hat\otimes}}
A$ кофунктор $\mathcal{X}^*_c$, получаем компактологический морфизм
$$
  \Chi^*_c(\mu)\colon \mathcal{X}^*_c(A\mathbin{\underset{\varepsilon}
  {\hat\otimes}}A)\longrightarrow\mathcal{X}^*_c(A)
$$
(гл.~3, $\S2$). Кроме того, существует компактологический изоморфизм
$$
  \widehat{J}\colon \mathcal{X}^*_c(A)\times\mathcal{X}^*_c(A)
  \longrightarrow\mathcal{X}^*_c(A\mathbin{\underset{\varepsilon}{\hat\otimes}}A)
$$
(гл.~4, $\S1$), который, напомним, каждой паре $(\chi_1,\ \chi_2)$
непрерывных эрмитовых характеров на $A$ сопоставляет непрерывный
эрмитов характер $\widehat{J}(\chi_1,\ \chi_2)$ на
$A\mathbin{\underset{\varepsilon}{\hat\otimes}} A$, определяемый
условием
$$
  \widehat{J}(\chi_1,\ \chi_2)(a\otimes b)=\chi_1(a)\chi(b)
$$
на элементарных тензорах $a\otimes b$ и продолжаемый непрерывно и
линейно на остальные элементы алгебры
$A\underset{\varepsilon}{\otimes}A$. Полагая
$$
  m_{\mathcal{X}^*_c(A)} = \Chi^*_c(\mu)\circ\widehat{J},
$$
получаем компактологический морфизм
$$
  m_{\mathcal{X}^*_c(A)}\colon \mathcal{X}^*_c(A)\times\mathcal{X}^*_c(A)
  \longrightarrow\mathcal{X}^*_c(A).
$$

\begin{lemma}
$1^\circ$ Для любой $K^*$-бигебры (соотв., $K^*$-алгебры Хопфа)
$(A,\ \mu)$ пара $(\mathcal{X}^*_c(A),\ m_{\mathcal{X}^*_c(A)})$
есть компактологическая полугруппа (соотв., группа).

$2^\circ$ Для любого морфизма $K^*$-бигебр $u$ $\Chi^*_c(u)$ есть
морфизм компактологических полугрупп.
\end{lemma}

\begin{proof} Это тотчас же следует из того, что $\mathcal{X}^*_c$
есть контравариантный функтор из $\mathcal{K}^*$ в $\mathcal{K}$
(гл.~3, $\S2$) и отмеченного выше свойства
отображения~$\widehat{J}.$
\end{proof}

Для любой $K^*$-алгебры $A$ преобразование Гельфанда  $\gel_A$
является изоморфизмом $A$ на $\mathcal{C}(\mathcal{X}^*_c(A))$ в
категории $K^*$-алгебр (гл.~3, $\S2$).

\begin{lemma}
Пусть $(A,\ \mu)$~--- $K^*$-бигебра. Тогда преобразование Гельфанда
$\gel_A$ является изоморфизмом $(A,\ \mu)$ на
$(\mathcal{C}(\mathcal{X}^*_c(A)),\ \mu_{\mathcal{X}^*_c(A)})$ в
категории $K^*$-бигебр.
\end{lemma}

\begin{proof} Нужно установить коммутативность диаграммы, то есть
показать, что для любого $a\in A$ имеет место равенство
$$
  \mu_{\mathcal{C}\mathcal{X}^*_c(A)}(\gel_A(a)) =
  (\gel_A\underset{\varepsilon}{\otimes}\gel_A)(\mu(a)).
$$

Пусть $a\in A$. Поскольку $\mu(a)\in
A\underset{\varepsilon}{\otimes}A$, существуют $a_{i\nu}$,
$b_{i\nu}$ такие, что
$$
  \mu(a)=\lim\limits_{\nu}\sum\limits^{n_{\nu}}_{i = 1}a_{i\nu}
  \otimes b_{i\nu}.
$$
Используя это, а также определение, непрерывность и
линейность $\gel_A\mathbin{\underset{\varepsilon}{\hat\otimes}}\gel_A$,
получаем:
$$
  (\gel_A\mathbin{\underset{\varepsilon}{\hat\otimes}}\gel_A)
  (\mu(a))=(\gel_A\mathbin{\underset{\varepsilon}{\hat\otimes}}
  \gel_A)(\lim\limits_{\nu}\sum\limits^{n_{\nu}}_{i=1}a_{i\nu}
  \otimes b_{i\nu})=
$$
$$
  \lim\limits_{\nu}(\gel_A\mathbin{\underset{\varepsilon}{\hat\otimes}}
  \gel_A)(\sum\limits^{n_{\nu}}_{i=1}a_{i\nu}\otimes b_{i\nu})=
  \lim\limits_{\nu}\sum\gel_A(a_{i\nu})\otimes\gel_A(b_{i\nu}).
$$
Отсюда, используя определение, непрерывность и линейность $\Phi$ и
$\widehat{J}$, для любых $\chi_1$, $\chi_2\in\Chi^*_c$ выводим:
\begin{align*}
  \Phi(\lim\limits_{\nu}\sum\limits^{n_{\nu}}_{i=1}\gel_A(a_{i\nu})\otimes\gel_A
  (b_{i\nu}))(\chi_1,\ \chi_2)
  &=\lim\limits_{\nu}\sum\limits^{n_{\nu}}_{i=1}\Phi(\gel_A(a_{i\nu})\otimes\gel_A
   (b_{i\nu}))(\chi_1, \chi_2)=\\
  &=\lim\limits_{\nu}\sum\limits^{n_{\nu}}_{i=1}\gel_A(a_{i\nu})(\chi_1)
  \cdot\gel_A(b_{i\nu})(\chi_2)=\\
  &=\lim\limits_{\nu}\sum\limits^{n_{\nu}}_{i=1}\chi_1(a_{i\nu})
  \cdot\chi_2(b_{i\nu}) =\\
  &= \lim\limits_{\nu}\sum\limits^{n_{\nu}}_{i=1}
  \widehat{J}(\chi_1,\chi_2)(a_{i\nu}\otimes b_{i\nu})=\\
  &=\widehat{J}(\chi_1,\chi_2)(\lim\limits_{\nu}
  \sum\limits^{n_{\nu}}_{i=1}a_{i\nu}\otimes b_{i\nu}) =\\
  &= \widehat{J}(\chi_1,\chi_2)(\mu(a)).
\end{align*}
C другой стороны, применяя определение и свойства $\mathcal{C}(m)$,
$\gel_A(a)$, $m_{\mathcal{X}^*_c(A)}$ и $\Chi^*_c(\mu)$, находим:
\begin{align*}
  (\Phi\circ\mu_{\mathcal{C}\mathcal{X}^*_c(A)}\circ\gel_A)(a)(\chi_1,\chi_2)
  &= \mathcal{C}(m_{\mathcal{X}^*_c(A)})(\gel_A(a))(\chi_1,\chi_2)=\\
  &= (\gel_A(a)\circ m_{\mathcal{X}^*_c(A)})(\chi_1,\chi_2)=\\
  &= \gel_A(a)(m_{\mathcal{X}^*_c(A)}(\chi_1,\chi_2))=\\
  &= m_{\mathcal{X}^*_c(A)}(\chi_1,\chi_2)(a)=\\
  &= (\Chi^*_c(\mu)\circ\widehat{J})(\chi_1,\chi_2)(a)=\\
  &= \Chi^*_c(\mu)(\widehat{J}(\chi_1,\chi_2))(a)=\\
  &= (\widehat{J}(\chi_1,\chi_2)\circ\mu)(a)=\\
  &= \widehat{J}(\chi_1,\chi_2)(\mu(a)).
\end{align*}
В силу произвольности $\chi_1$, $\chi_2$ и $a$ это означает, что
$\Phi\circ(\gel_A\mathbin{\underset{\varepsilon}{\hat\otimes}}\gel_A)
\circ\mu=\Phi\circ\mu_{\mathcal{C}\mathcal{X}^*_c(A)}\circ\gel_A$.
Опуская изоморфизм $\Phi$, получаем требуемое.
\end{proof}

Для любого компактологического пространства $S$ преобразование
Дирака
$$
  \delta_S\colon S\longrightarrow\mathcal{X}^*_c(\mathcal{C}(S))
$$
является изоморфизмом компактологических пространств (гл. 3, $\S2$,
доказательство теоремы (2.1)).

\begin{lemma} Для любой компактологической полугруппы $(S,\ m)$
преобразование Дирака $\delta_S$ является изоморфизмом $(S,m)$ на
$(\mathcal{X}^*_c\mathcal{C}(S), m_{\mathcal{X}^*_c\mathcal{C}(S)})$
в категории компактологических полугрупп.
\end{lemma}

\begin{proof} Нужно установить коммутативность диаграммы, то есть
показать, что для любых $(x_1, x_2)\in S\times S$ и
$\varphi\in\mathcal{C}(S)$ имеет место равенство
$$
  (\delta_S\circ m)(x_1, x_2)(\varphi) =
  (m_{\mathcal{X}^*_c\mathcal{C}(S)}\circ(\delta_S\times\delta_S))
  (x_1, x_2)(\varphi).
$$

Пусть $(x_1,\ x_2)\in S\times S$ и $\varphi\in\mathcal{C}(S)$.
Используя определения и свойства $m_{\mathcal{X}^*_c\mathcal{C}(S)}$
и $\Chi^*_c(u)$, находим:
\begin{align*}
  [m_{\mathcal{X}^*_c\mathcal{C}(S)}\circ(\delta_S\times\delta_S)](x_1,
  x_2)(\varphi) &=[\Chi^*_c(\Phi^{-1}\circ\mathcal{C}(m))\circ\widehat{J}\circ
  (\delta_S\times\delta_S)](x_1, x_2)(\varphi) =\\
  &=[\Chi^*_c(\Phi^{-1}\circ\mathcal{C}(m))\circ\widehat{J}](\delta_S(x_1),
  \delta_S(x_2))(\varphi)= \\ &=\Chi^*_c(\Phi^{-1}\circ\mathcal{C}(m))(\widehat{J}
  (\delta_S(x_1),\delta_S(x_2)))(\varphi)=\\
  &=[\widehat{J}(\delta_S(x_1),\delta_S(x_2))\circ\Phi^{-1}\circ\mathcal{C}(m)]
   (\varphi)= \\ 
  &= \widehat{J}(\delta_S(x_1),\delta_S(x_2))(\Phi^{-1}(\mathcal{C}(m)
   (\varphi)))= \\
  &=\widehat{J}(\delta_S(x_1),\delta_S(x_2))(\Phi^{-1}(\varphi\circ
  m)).
\end{align*}

Теперь, $\Phi^{-1}(\varphi\circ m)\in\mathcal{C}(S)
\mathbin{\underset{\varepsilon}{\hat\otimes}}\mathcal{C}(S)$, поэтому
существуют $\varphi_{i\nu}$, $\psi_{i\nu}\in\mathcal{C}(S)$ такие,
что
$$
  \Phi^{-1}(\varphi\circ m) = \lim\limits_{\nu} \sum
  \limits^{n_\nu}_{i=1} \varphi_{i\nu}\otimes\psi_{i\nu}.
$$
Отсюда, используя непрерывность и линейность $\widehat{J}(\chi_1,\
\chi_2)$, $\Phi^{-1}$ и $\Phi$, получаем:
\begin{align*}
  [m_{\mathcal{X}^*_c\mathcal{C}(S)}\circ(\delta_S\times\delta_S)](x_1,
  x_2)(\varphi)
  &= \widehat{J}(\delta_S(x_1),\delta_S(x_2))(\lim\limits_{\nu}
     \sum\limits^{n_\nu}_{i=1}\varphi_{i\nu}\otimes\psi_{i\nu})=\\
  &= \lim\limits_{\nu} \sum\limits^{n_\nu}_{i=1} \widehat{J}(\delta_S(x_1),
     \delta_S(x_2))(\varphi_{i\nu}\otimes\psi_{i\nu})=\\
  &= \lim\limits_{\nu}\sum\limits^{n_\nu}_{i=1}
     \varphi_{i\nu}(x_1)\cdot\psi_{i\nu}(x_2)=\\
  &= \lim\limits_{\nu}\sum\limits^{n_\nu}_{i=1}
     \Phi(\varphi_{i\nu}\otimes\psi_{i\nu})(x_1, x_2)=\\
  &= \Phi(\lim\limits_{\nu}\sum\limits^{n_\nu}_{i=1}
     \varphi_{i\nu}\otimes\psi_{i\nu})(x_1, x_2)=\\
  &= \Phi(\Phi^{-1}(\varphi\circ m))(x_1, x_2)=\\
  &= (\varphi\circ m)(x_1, x_2)=\\
  &= \varphi(m(x_1, x_2))=\delta_S(m(x_1, x_2))(\varphi) =\\
  &= (\delta_S\circ m)(x_1, x_2)(\varphi).
\end{align*}
Это и требовалось.
\end{proof}

Сопоставляя леммы (2.1)--(2.4), видим, что нами доказана следующая
теорема.

\begin{thm}
Контравариантные функторы $\mathcal{C}$ и $\mathcal{X}^*_c$ образуют
двойственность между категорией компактологических полугрупп
(соотв., групп) и категорией $K^*$-бигебр (соотв., $K^*$-алгебр
Хопфа).
\end{thm}

\section{Теоремы двойственности для топологических полугрупп}

Обозначим через $\mathcal{K}_{\mathbb{C}}$ категорию
$k_{\mathbb{C}}$-пространств и непрерывных отображений. Это~---
категория с конечными произведениями. Произведением
$k_{\mathbb{C}}$-пространств $T_1$, $T_2$ в категории
$\mathcal{K}_{\mathbb{C}}$ является не тихоновское произведение
$T_1\times T_2$, которое может не быть
$k_{\mathbb{C}}$-пространством (гл.~4, $\S 2$), а пространство
$T_1\boxtimes T_2$, введенное в $\S 2$ гл.~4. Напомним, что по
определению это есть декартово произведение $T_1\times T_2$ множеств
$T_1$ и $T_2$, наделённое слабой топологией, определяемой всеми
функциями $T_1\times T_2\longrightarrow \Cbb$, непрерывными на
каждом компакте вида $K_1\times K_2$, где $K_1$, $K_2$~---
компактное подмножество в $T_1$, $T_2$ соответственно.

Полугруппы над категорией $\mathcal{K}_{\mathbb{C}}$ мы будем
называть $k_{\mathbb{C}}$-полугруппами. Иными словами,
$k_{\mathbb{C}}$-полугруппа~--- это $k_{\mathbb{C}}$-пространство~ $T$, наделённое непрерывным отображением
$$
  m\colon T\boxtimes T\longrightarrow T,
$$
для которого коммутативна диаграмма

Для любого $k_{\mathbb{C}}$-пространства $T$ обозначим через
$\widetilde{T}$ ассоциированное компактологическое пространство, то
есть множество $T$, наделённое компактологией, состоящей в точности
из всех компактных подмножеств топологического пространства $T$.
Ясно, что отображение $u: T_1\rightarrow T_2$ непрерывно тогда и
только тогда, когда оно является компактологическим морфизмом из
$\widetilde{T}_1$ в $\widetilde{T}_2$. (Действительно, поскольку
$T_1$~--- $k_{\mathbb{C}}$-пространство, это верно для всех
отображений из $T_1$ в $\mathbb{C}$, а поскольку $T_2$ вполне
регулярно, это верно и для отображений из $T_1$ в $T_2$.) Тем самым,
категория $\mathcal{K}_{\mathbb{C}}$ эквивалентна (даже изоморфна)
полной подкатегории $\widetilde{\mathcal{K}}_{\mathbb{C}}$ категории
$\mathcal{K}$, образованной компактологическими пространствами вида
$\widetilde{T}$.

Далее, легко видеть, что
$$
  \widetilde{T_1\boxtimes T_2}=\widetilde{T}_1\times\widetilde{T}_2.
$$
В силу этого $\widetilde{\mathcal{K}}_{\mathbb{C}}$ есть
подкатегория, замкнутая относительно произведений в $\mathcal{K}$.
Стало быть: а) категория полугрупп над $\mathcal{K}_{\mathbb{C}}$
эквивалентна категории полугрупп над
$\widetilde{\mathcal{K}}_{\mathbb{C}}$ и б) полугруппами над
$\widetilde{\mathcal{K}}_{\mathbb{C}}$ являются просто полугруппы
$(S,\ m)$ над $\mathcal{K}$, у которых $S\in
\ob\widetilde{\mathcal{K}}_{\mathbb{C}}$. Поэтому эквивалентны
категории, дуальные к категории полугрупп над
$\mathcal{K}_{\mathbb{C}}$ и к категории полугрупп над
$\widetilde{\mathcal{K}}_{\mathbb{C}}$. Эта же категория
эквивалентна категории тех $K^*$-бигебр $(A,\ \mu)$, у которых
$A\cong\mathcal{C}(S)$ для некоторого $S\in
\ob\widetilde{\mathcal{K}}_{\mathbb{C}}$ и, значит,
$A\cong\mathcal{C}_{co}(S)$ для некоторого $S\in
\ob\mathcal{K}_{\mathbb{C}}$. Вспоминая, что, в силу теоремы~(3.2) из
гл.~3, $A\cong\mathcal{C}_{co}(S)$ для некоторого $S\in
\ob\mathcal{K}_{\mathbb{C}}$ тогда и только тогда, когда на
$K^*$-алгебре $A$ каждая спектральная $C^*$-преднорма непрерывна,
видим, что нами доказана

\begin{thm}
Категория, дуальная к категории $k_\Cbb$-полугрупп, (соотв., групп)
эквивалентна категории $K^*$-бигебр (соотв., $K^*$-алгебр Хопфа), на
которых все спектральные $C^*$-преднормы непрерывны.
\end{thm}

Теперь обратимся к подкатегориям категории $\mu k_\Cbb$-полугрупп.

\begin{lemma}
$1^\circ$ Пусть $T_1$, $T_2$~--- $\mu$-пространства. Тогда и
$T_1\boxtimes T_2$~--- $\mu$-пространство.

$2^\circ$ Пусть $A_1$, $A_2$~--- $C^*$-бочечные $K^*$-алгебры. Тогда
и $A_1\mathbin{\underset{\varepsilon}{\hat\otimes}}A_2$~---
$C^*$-бочечная $K^*$-алгебра.
\end{lemma}

\begin{proof} 1. Пусть $B$~--- замкнутое ограниченное подмножество в
$T_1\boxtimes T_2$. Тогда оно ограниченно и в $T_1\times T_2$.
Поскольку $T_1\times T_2$~--- $\mu$-пространство, замыкание
$\overline{B}$ множества $B$ в $T_1\times T_2$ в компактно в
$T_1\times T_2$. Но запас компактов у $T_1\boxtimes T_2$ тот же, что
и у $T_1\times T_2$. Значит, $\overline{B}$ компактно и в
$T_1\boxtimes T_2$. В силу этого $B$ оказывается замкнутым
подмножеством компактного множества и потому само компактно.

2. Пусть $A_1, A_2$~--- $C^*$-бочечные $K^*$-алгебры. Тогда (гл.~3,
теорема (3.4)) $\Chi^*_c(A_1)$ и $\Chi^*_c(A_2)$~---
$\mu$-пространства и $A_1\cong C_{co}(\Chi^*_c(A_1))$, $A_2\cong
C_{co}(\Chi^*_c(A_2))$. Применяя теорему (2.8) из гл.~4, получаем:
$$A_1\mathbin{\underset{\varepsilon}{\hat\otimes}}A_2\cong C_{co}(\Chi^*_c(A_1))
\mathbin{\underset{\varepsilon}{\hat\otimes}}C_{co}(\Chi^*_c(A_2))\cong
C_{co}(\Chi^*_c(A_1)\boxtimes\Chi^*_c(A_2)).$$ В силу п. $1^\circ$
$\Chi^*_c(A_1)\boxtimes\Chi^*_c(A_2)$ есть
$\mu k_\Cbb$-пространство. Значит,
$C_{co}(\Chi^*_c(A_1)\boxtimes\Chi^*_c(A_2))$, а вместо с ней~--- и
$A_1\mathbin{\underset{\varepsilon}{\hat\otimes}}A_2$ есть
$C^*$-бочечная $K^*$-алгебра (снова применяем теорему (3.4) из
гл.~3).
\end{proof}

Условимся называть $\mu k_\Cbb$-полугруппами те $k_\Cbb$-полугруппы,
подлежащее пространство которых является $\mu$-пространством. Еще
раз применяя теорему (3.4) из гл.~3 получаем:

\begin{thm}
Категория, дуальная к категории $\mu k_\Cbb$-полугрупп (соотв.,
групп), эквивалентна категории $C^*$-бочечных $K^*$-бигебр (соотв.,
$K^*$-алгебр Хопфа).
\end{thm}

\begin{thm}
$1^\circ$ Пусть $T_1$, $T_2$~--- $Q$-пространства. Тогда и
$T_1\boxtimes T_2$ есть $Q$-пространство.

$2^\circ$ Пусть $A_1$, $A_2$~--- борнологические $K^*$-алгебры.
Тогда и $A_1\mathbin{\underset{\varepsilon}{\hat\otimes}}A_2$ есть
борнологическая $K^*$-алгебра.
\end{thm}

\begin{proof} 1. Известно, что произведение $Q$-пространств является
$Q$-пространством. Значит, $T_1\times T_2$ есть $Q$-пространство,
так что $\Chi^*C(T_1\times T_2) = T_1\times T_2$ (гл.~2,
теорема~(1.9)). Но $C(T_1\times T_2)\subset C(T_1\boxtimes T_2)$,
поэтому
$$
  T_1\boxtimes T_2\subset\Chi^*C(T_1\boxtimes T_2)
  \subset\Chi^*C(T_1\times T_2) = T_1\times T_2.
$$
Тем самым, $\Chi^*C(T_1\boxtimes T_2)=T_1\boxtimes T_2$ и
$T_1\boxtimes T_2$ есть $Q$-пространство.

2. Пусть $A_1$, $A_2$~--- борнологические $K^*$-алгебры. Тогда
$\Chi^*_c(A_1)$ и $\Chi^*_c(A_2)$ являются $Q$-пространствами,
причем $A_1\cong C_{co}\Chi^*_c(A_1)$ и $A_2\cong
C_{co}\Chi^*_c(A_2)$ (теорема (3.5) гл.~3). Поэтому
$$
  A_1\mathbin{\underset{\varepsilon}{\hat\otimes}}A_2\cong
  C_{co}(\Chi^*_c(A_1)\boxtimes\Chi^*_c(A_2)),
$$
где $\Chi^*_c(A_1)\boxtimes\Chi^*_c(A_2)$, в силу сказанного, есть
$Q$-пространство. В силу той же теоремы (3.5) из гл.~3 алгебра
$C_{co}(\Chi^*_c(A_1)\boxtimes\Chi^*_c(A_2))$, а вместе с ней~--- и
$A_1\mathbin{\underset{\varepsilon}{\hat\otimes}}A_2$ является
борнологическими $K^*$-алгебрами.
\end{proof}

Еще раз применяя теорему (3.5) из гл.~3, приходим к следующему
результату.

\begin{thm}
Категория, дуальная к категории вещественно-компактных
$k_\Cbb$-полугрупп (соотв., групп), эквивалентна категории
борнологических $K^*$-бигебр (соотв., $K^*$-алгебр Хопфа).
\end{thm}

Теперь обратимся к локально компактным полугруппам. Как известно,
локально компактные пространства являются $k$-пространствами и, тем
самым, $k_\Cbb$-пространст\-вами, причём тихоновское
произведение локально компактных пространств локально компактно. В
силу этого: а) категория локально компактных пространств является
полной подкатегорией категории $k_\Cbb$-пространств и
б) полугруппами над категорией локально компактных пространств
является локально компактные полугруппы в обычном смысле слова.
Значит, категория, дуальная к категории локально компактных
полугрупп, эквивалентна категории тех $K^*$-бигебр, у которых
$A\cong C_{co}(T)$ для некоторых локально компактного пространства
$T$, то есть тех, у которых $A$ есть $lc^*$-алгебра (гл.~3, теорема~(4.2)). Называя такие $K^*$-бигебры $lc^*$-бигебрами, получаем основной результат этой главы.

\begin{thm}
Категория, дуальная к категории отделимых локально компактных
полугрупп (соотв., групп), эквивалентна категории $lc^*$-бигебр
(соотв., $lc^*$-алгебр Хопфа).
\end{thm}

Хорошо известно, что  подлежащее пространство локально компактной
группы паракомпактно, а все паракомпактные пространства являются
$\mu$-пространствами. Отсюда~---
\begin{cor}
Каждая $lc^*$-алгебра Хопфа $C^*$-бочечна.
\end{cor}

Сочетая теорему (3.7) с теоремами (3.4) и (3.6), получаем еще две
теоремы двойственности.

\begin{thm}
Категория, дуальная к категории $\mu$-топологических локально
компактных полугрупп, эквивалентна категории $C^*$-бочечных
$lc^*$-бигебр.
\end{thm}

\begin{thm} Категория, дуальная к категории вещественно-компактных
локально компактных полугрупп (соотв., групп), эквивалентна
категории борнологических $lc^*$-бигебр (соотв., $lc^*$-алгебр
Хопфа).
\end{thm}

\section{Понтрягинская двойственность}

В этом параграфе мы покажем, как выглядит двойственность Понтрягина
для коммутативных локально компактных групп в рамках
пространственной выше теории.

Как и в случае обычных бигебр, элемент $a$ $K^*$-бигебры $(A,\ \mu)$
будем называть примитивным, если он отличен от нуля и имеет место
равенство:
$$
  \mu(a)=a\otimes a.
$$
Купер и Михор в~\cite{b20} и Купер в~\cite{b19} вводят еще понятие
строго примитивного элемента. А именно, элемент $a$ $K^*$-алгебры
Хопфа $(A,\ \mu,\ \varepsilon,\ \sigma)$ можно было бы, следуя
Куперу и Михору, называть строго примитивным, если он примитивен и:

а) $\varepsilon(a)=1$ ;

б) $\sigma(a)=a^{-1}$.

Никакой надобности в этом, однако, нет.

\begin{prop}
В любой $K^*$-алгебре Хопфа $(A,\ \mu,\ \varepsilon,\ \sigma)$ каждый примитивный элемент строго примитивен.
\end{prop}

\begin{proof} Пусть $a$~--- примитивный элемент $K^*$-алгебры Хопфа $(A,\ \mu,\ \varepsilon,\ \sigma)$.

а) В силу коммутативности диаграммы $(E^*)$ из $\S1$ имеем:
$$
  (\varepsilon\mathbin{\underset{\varepsilon}{\hat\otimes}}I_A)\mu(a)=1\otimes a.
$$
При этом
$$
  (\varepsilon\mathbin{\underset{\varepsilon}{\hat\otimes}}I_A)\mu(a)=
  (\varepsilon\mathbin{\underset{\varepsilon}{\hat\otimes}}I_A)(a\otimes
  a)=\varepsilon(a)\otimes a,
$$
где $\varepsilon(a)$~--- число. Значит
$$
  \varepsilon(a)\otimes a=\varepsilon(a)(1\otimes a),
$$
так что
$$
  \varepsilon(a)(1\otimes a)=1\otimes a
$$
Поскольку $a\neq 0$, $1\otimes a\neq 0$. Следовательно,
$\varepsilon(a)=1$.

б) В силу коммутативности диаграммы $(S^*)$ из $\S1$ имеем:
$$
  [m\circ(\sigma\otimes I_A)\circ\mu](a)=\varepsilon(a)\cdot
  e=1\cdot e=e.
$$
C другой стороны,
$$
  [m\circ(\sigma\otimes I_A)\circ\mu](a)=m((\sigma\otimes
  I_A)(a\otimes a)=m(\sigma(a)\otimes a)=\sigma(a)\cdot a.
$$
Таким образом, $\sigma(a)\cdot a = e$. Аналогично проверяется, что и
$a\cdot\sigma(a) = e$. Вместе эти равенства означают, что $\sigma(a)
= a^{-1}$. Предложение доказано.
\end{proof}

\begin{prop} Элемент $a$ $K^*$-алгебры Хопфа $(A,\ \mu)$
примитивен тогда и только тогда, когда $\gel_A(a)$ является
(непрерывным) морфизмом топологической группы $\Chi^*_c(A)$ в
мультипликативную группу $\mathbb{C}^*$ отличных от нуля комплексных
чисел.
\end{prop}

\begin{proof} Напомним, что умножение $\ast$ на полугруппе
$\Chi^*_c(A)$ определяется как композиция
$$\Chi^*_C(A)\times\Chi^*_c(A)\overset{\widehat{J}}{\longrightarrow}
\Chi^*_c(A\mathbin{\underset{\varepsilon}{\hat\otimes}} A)\overset{\Chi^*_c(\mu)}{\longrightarrow}\Chi^*_c(A).$$

Используя это, для любого $a\in A$ и любых $\chi_1,
\chi_2\in\Chi^*_c(A)$ получаем:
$$\gel_A(a)(\chi_1\ast\chi_2)=\gel_A(a)[(\Chi^*_c(\mu)\circ\widehat{J})(\chi_1,
\chi_2)]=[\Chi^*_c(\mu)(\widehat{J}(\chi_1,\chi_2))](a)=$$
$$=[\widehat{J}(\chi_1, \chi_2)\circ\mu](a)=\widehat{J}(\chi_1,
\chi)(\mu(a)),$$ и, с другой стороны,
$$\gel_A(a)(\chi_1)\cdot\gel_A(a)(\chi_2)=\chi_1(a)\cdot\chi_2(a)=\widehat{J}(\chi_1, \chi_2)(a\otimes a).$$

(Кроме того, известно, что все $\gel_A(a)$ непрерывны.)

Пусть теперь $a$ примитивен. Тогда, в силу предложения (4.1), он
обратим, поэтому функция $\gel_A(a)$ всюду отлична от нуля, то есть
является отображением в $\mathbb{C}^*$. Из равенства $\mu(x) =
x\otimes x$ вытекает, что последние члены выписанных выше двух
цепочек равенства равны. Значит, равны и их первые члены, а это и
означает, что $\gel_A(a)$ является морфизмом полугрупп.

Обратно, пусть $\gel_A(a)$~--- морфизм из полугруппы $\Chi^*_c(A)$ в
полугруппу $\mathbb{C}^*$. Тогда при всех $\chi_1$, $\chi_2$ из
$\Chi^*_c(A)$ равны первые члены тех же цепочек равенств и, значит,
их последние члены. Таким образом, при всех $(\chi_1,\
\chi_2)\in\Chi^*_c(A)\times\Chi^*_c(A)$
$$
  \widehat{J}(\chi_1,\ \chi_2)(\mu(a))=\widehat{J}(\chi_1,\ \chi_2)(a\otimes a).
$$
Поскольку отображение $\widehat{J}$ сюръективно, это означает, что
для любого непрерывного эрмитова характера $\chi$ на
$A\mathbin{\underset{\varepsilon}{\hat\otimes}} A$ имеет место равенство
$\chi(\mu(a))=\chi(a\otimes a).$ При этом алгебра
$A\mathbin{\underset{\varepsilon}{\hat\otimes}} A$ строго $*$-полупроста
(глава 4, $\S2$). Значит, $\mu(a)=a\otimes a$, что и требовалось.
\end{proof}

Элемент $a$ топологической инволютивной алгебры $A$ будем называть
ограниченным, если $\{ \chi(a): \chi\in\Chi^*_c(A)\}$ есть
ограниченное множество в $\mathbb{C}$.

\begin{prop} $1^\circ$ Пусть $(A,\ \mu)$~--- $lc^*$-алгебра
Хопфа, $a\in A$. Функция $\gel_A(a)$ является характером дуальной к
$(A,\ \mu)$ локально компактной группы $\Chi^*_c(A)$ тогда и только
тогда, когда $a$ есть ограниченный примитивный элемент в $(A,\
\mu)$.

$2^\circ$ Пусть $G$-отделимая локально компактная группа,
$\varphi\in C(G)$. Функция $\varphi$ является характером группы $G$
тогда и только тогда, когда она является ограниченным примитивным
элементом дуальной к $G$ алгебры Хопфа $C_{co}(G).$
\end{prop}

\begin{proof} 1. В силу предложения (5.2) $\gel_A(a)$ является
морфизмом из $\Chi^*_c(A)$ в $\mathbb{C}^*$. Ограниченность $a$ есть
не что иное, как ограниченность функции $\gel_A(a)$. Из неё же
стандартным образом выводится, что для любого $\chi\in\Chi^*_c(A)$
$|\gel_A(a)(\chi)|=1$. Тем самым $\gel_A(a)$ оказывается морфизмом
из $\Chi^*_c(A)$ в единичную окружность, то есть характером.

2. Это тотчас же следует из части $1^\circ$ и того, что группы $G$ и
$\Chi^*_c(C_{co}(G))$ топологически изоморфны ($\S3$).
\end{proof}

Несложная проверка показывает, что множество $P(A)$ всех и множество
$P_b(A)$ всех ограниченных примитивных элементов $K^*$-алгебры Хопфа
$(A,\ \mu)$ являются подгруппами группы всех обратимых элементов
алгебры $A$. Поскольку $K^*$-алгебры отделимы и умножение на них
совместно непрерывно, множества $P(A)$ и $P_b(A)$ с индуцированными
из $A$ топологиями являются отделимыми топологическими группами.
Беря в качестве $(A,\ \mu)$ алгебру Хопфа $C_{co}(G)$, дуальную к
локально компактной группе~$G$, и вспоминая, что группа
$\widehat{G}$, дуальная по Понтрягину к группе~$G$, есть группа всех
непрерывных характеров группы~$G$, наделённая, как и $C_{co}(G)$,
компактно-открытой топологией, приходим к следующему результату.

\begin{prop} Группа $\widehat{G}$, дуальная по Понтрягину к
отделимой локально компактной группе $G$, есть в точности группа
всех ограниченных примитивных элементов дуальной к $G$
$lc^*$-алгебры Хопфа $C_{co}(G)$.
\end{prop}

%
%
\bibliographystyle{amsplain}

\end{document}